\documentclass[12pt,a4paper,leqno]{article}

\usepackage{latexsym}
\usepackage{amssymb,exscale}
\usepackage[centertags]{amsmath}
\usepackage{amsthm}
\usepackage[all]{xy}
\usepackage{graphicx}
\usepackage{multicol}
\usepackage{makeidx}
\usepackage{url}
\newcommand{\qqed}{}
\newcommand{\no}[1]{}

\numberwithin{equation}{subsection}

\swapnumbers
\theoremstyle{definition}
\newtheorem{theorem}[equation]{Theorem}
\newtheorem{lemma}[equation]{Lemma}
\newtheorem{corollary}[equation]{Corollary}
\newtheorem{definition}[equation]{Definition}
\newtheorem{example}[equation]{Example}
\newtheorem{remark}[equation]{Remark}
\newtheorem{question}[equation]{Question}

\renewcommand{\phi}{\varphi}

\newcommand{\I}{{\rm i}}
\newcommand{\D}{\mathrm{d}}
\newcommand{\E}{\mathrm{e}}

\newcommand{\ti}{\tilde}

\renewcommand{\(}{\bigl(}
\renewcommand{\)}{\bigr)\vphantom{)}}

\newcommand{\ip}[2]{\langle#1,#2\rangle}

\newcommand{\equi}{\;\Longleftrightarrow\;}
\newcommand{\imply}{\;\;\;\Longrightarrow\;\;\;}
\newcommand{\impl}{\;\Longrightarrow\;}
\newcommand{\imp}{$ \Longrightarrow $ }

\newcommand{\sgn}{\mathrm{sgn}}
\newcommand{\Exp}{\mathrm{Exp}}

\newcommand{\pd}{\partial}

\newcommand{\jetblack}{\text{\textup{jetblack}}}
\newcommand{\countable}{\text{\textup{countable}}}
\newcommand{\perfect}{\text{\textup{perfect}}}
\newcommand{\stable}{\text{\textup{stable}}}
\newcommand{\discrete}{\text{\textup{discrete}}}

\newcommand{\all}{\operatorname{all}}

\newcommand{\Proj}{\operatorname{Proj}}
\newcommand{\MALG}{\operatorname{MALG}}
\newcommand{\Lim}{\operatorname{Lim}}
\newcommand{\mes}{\operatorname{mes}}
\newcommand{\Var}{\operatorname{Var}}
\newcommand{\Cov}{\operatorname{Cov}}
\newcommand{\Corr}{\operatorname{Corr}}
\newcommand{\N}{\operatorname{N}}
\newcommand{\Poisson}{\operatorname{Poisson}}
\newcommand{\Geom}{\operatorname{Geom}}
\newcommand{\entier}{\operatorname{entier}}

\newcommand{\One}{\mathbf1}

\newcommand{\bH}{\mathbf H}
\newcommand{\bN}{\mathbf N}

\newcommand{\eps}{\varepsilon}
\newcommand{\si}{\sigma}
\newcommand{\ga}{\gamma}
\newcommand{\om}{\omega}
\newcommand{\Om}{\Omega}

\newcommand{\De}{\Delta}
\newcommand{\al}{\alpha}
\newcommand{\be}{\beta}
\newcommand{\cC}{\mathcal C}
\newcommand{\cM}{\mathcal M}
\newcommand{\cN}{\mathcal N}
\newcommand{\Ec}{\mathcal E}
\newcommand{\F}{\mathcal F}
\newcommand{\A}{\mathcal A}
\newcommand{\B}{\mathcal B}
\renewcommand{\P}{\mathcal P}
\newcommand{\la}{\lambda}

\newcommand{\Cfinite}{\mathcal{C}_{\text{finite}}}
\newcommand{\dist}{\operatorname{dist}}
\newcommand{\const}{{\mathrm{const}}}
\newcommand{\Ex}{\mathbb E\,}
\renewcommand{\Pr}[1]{\mathbb P\,\(\,#1\,\)}
\newcommand{\Prob}{\mathbb P}
\newcommand{\R}{\mathbb R}
\newcommand{\C}{\mathbb C}
\newcommand{\Q}{\mathbb Q}
\newcommand{\Z}{\mathbb Z}
\newcommand{\T}{\mathbb T}
\newcommand{\cE}[2]{\mathbb{E}\,\(\,#1\,\big|\,#2\,\)}
\newcommand{\cP}[2]{\mathbb{P}\,\(\,#1\,\big|\,#2\,\)}
\newcommand{\cVar}[2]{\operatorname{Var}\,\(\,#1\,\big|\,#2\,\)}

\newcommand{\sif}{$\sigma$\nobreakdash-field}

\newcommand{\tendsk}{\xrightarrow[k\to\infty]{}}
\newcommand{\tendsi}{\xrightarrow[i\to\infty]{}}
\newcommand{\tendseps}{\xrightarrow[\eps\to0]{}}
\newcommand{\tendst}{\xrightarrow[t\to0]{}}

\newenvironment{myitemize}{\begin{list}{$\bullet$}
{\setlength{\topsep}{1mm}
\setlength{\partopsep}{0mm}
\setlength{\itemsep}{1mm}
\setlength{\parsep}{0mm}
\setlength{\parskip}{0mm}}}
{\end{list}}

\def\emailwww#1#2{\par\quad {\tt #1}\par\quad {\tt #2}\medskip}
\makeatletter
  {\end{multicols}}%
\makeatother
\makeatletter
\renewcommand*\l@section[2]{%
  \ifnum \c@tocdepth >\z@
    \addpenalty\@secpenalty
    \addvspace{0.25em \@plus\p@}%
    \setlength\@tempdima{2.5em}%
    \begingroup
      \parindent \z@ \rightskip \@pnumwidth
      \parfillskip -\@pnumwidth
      \leavevmode \bfseries
      \advance\leftskip\@tempdima
      \hskip -\leftskip
      #1\nobreak\hfil \nobreak\hb@xt@\@pnumwidth{\hss #2}\par
    \endgroup
  \fi}
\renewcommand*\numberline[1]{\hb@xt@\@tempdima{\hfil#1\hskip1em}}
\makeatother

\begin{document}

\title{Scaling Limit, Noise, Stability}

\author{Boris Tsirelson}

\date{}
\maketitle

\begin{abstract}
Linear functions of many independent random variables lead to
classical noises (white, Poisson, and their combinations) in the
scaling limit. Some singular stochastic flows and some models of
oriented percolation involve very nonlinear functions and lead to
nonclassical noises. Two examples are examined, Warren's `noise made
by a Poisson snake' and the author's `Brownian web as a black
noise'. Classical noises are stable, nonclassical are not. A new
framework for the scaling limit is proposed. Old and new results are
presented about noises, stability, and spectral measures.
\end{abstract}

\setcounter{tocdepth}{1}
\tableofcontents

\vfill

\section*{Introduction}
\addcontentsline{toc}{section}{Introduction}
Functions of $ n $ independent random variables and limiting
procedures for $ n \to \infty $ are a tenor of probability theory.

Classical limit theorems investigate linear functions, such as $ f
(\xi_1,\dots,\xi_n) = ( \xi_1 + \dots + \xi_n ) / \sqrt n $. The
well-known limiting procedure (a classical
example of scaling limit) leads to the Brownian motion. Its
derivative, the white noise, is not a continuum of independent random
variables, but rather an infinitely divisible `reservoir of
independence', a classical example of a continuous product of
probability spaces.

Percolation theory investigates some very special nonlinear functions
of independent two-valued random variables, either in the limit of an
infinite discrete lattice, or in the scaling limit. The latter is now
making spectacular progress. The corresponding `reservoir of
independence' is already constructed for oriented percolation (which
is much simpler). That is a modern, nonclassical example of a
continuous product of probability spaces.

An essential distinction between classical and nonclassical
continuous products of probability spaces is revealed by the concept
of stability/sensitivity, framed for the discrete case by computer
scientists and (in parallel) for the continuous case by
probabilists. Everything is stable if and only if the setup is
classical.

Some readers prefer discrete models, and treat continuous models as
a mean of describing asymptotic behavior. Such readers may skip
Sects.\ \ref{sec:6.2}, \ref{sec:6.3}, \ref{sec:8.2},
\ref{sec:8.3}, \ref{sec:8.4}. Other readers are interested only
in continuous models. They may restrict themselves to Sects.\
\ref{sec:3.4}, \ref{sec:3.5}, \ref{sec:4.9}, \ref{sec:5.2},
\ref{sec:6}, \ref{sec:7}, \ref{sec:8}.

\emph{Scaling limit.} A new framework for the scaling limit is
proposed in Sects.\ \ref{sec:1.3}, \ref{sec:2},
\ref{sec:3.1}--\ref{sec:3.3}.

\emph{Noise.} The idea of a continuous product of probability spaces
is formalized by the notions of `continuous factorization'
(Sect.\ \ref{sec:3.4}) and `noise' (Sect.\ \ref{sec:3.5}). (Some other
types of continuous product are considered in \cite{Ts99},
\cite{Ts02}.) For two nonclassical examples of noise see
Sects.\ \ref{sec:4}, \ref{sec:7}.

\emph{Stability.} Stability (and sensitivity) is studied in
Sects.\ \ref{sec:5}, \ref{sec:6.1}, \ref{sec:6.4}. For an interplay
between discrete and continuous forms of stability/sensitivity, see
especially Sects.\ \ref{sec:5.3}, \ref{sec:6.4}.

The spectral theory of noises, presented in Sects.\ \ref{sec:3.3},
\ref{sec:3.4} and used in Sects.\ \ref{sec:5}, \ref{sec:6}, generalizes
both the Fourier transform on the discrete group $ \Z_2^n $ (the
Fourier-Walsh transform) and the It\^o decomposition into multiple
stochastic integrals. For the scaling limit of spectral measures, see
Sect.\ \ref{sec:3.3}.

Throughout, either by assumption or by construction, all probability
spaces will be Lebesgue-Rokhlin spaces; that is, isomorphic $ \bmod \,
0 $ to an interval with Lebesgue measure, or a discrete (finite or
countable) measure space, or a combination of both.

\section{A First Look}
\label{sec:1}
\subsection{Two toy models}
\label{sec:1.2}

The most interesting thing is a scaling limit as a transition from a
lattice model to a continuous model. A transition from a finite
sequence to an infinite sequence is much simpler, but still
nontrivial, as we'll see on simple toy models.

Classical theorems about independent increments are exhaustive, but a
small twist may surprise us. I demonstrate the twist on two models,
`discrete' and `continuous'. The `continuous' model is a Brownian
motion on the circle. The `discrete' model takes on two values $ \pm 1
$ only, and increments are treated multiplicatively: $ X(t) / X(s) $
instead of the usual $ X(t) - X(s) $. Or equivalently, the `discrete'
process takes on its values in the two-element group $ \Z_2 $; using
additive notation we have $ \Z_2 = \{ 0, 1 \} $, $ 1+1=0 $, increments
being $ X(t) - X(s) $. In any case, the twist stipulates values in a
compact group (the circle, $ \Z_2 $, etc.), in contrast to the
classical theory, where values are in $ \R $ (or another linear
space). Also, the classical theory assumes continuity (in
probability), while our twist does not. The `continuous' process (in
spite of its name) is discontinuous at a single instant $ t=0 $. The
`discrete' process is discontinuous at $ t = \frac1n $, $ n =
1,2,\dots $, and also at $ t = 0 $; it is constant on $ [\frac1{n+1},
\frac1n) $ for every $ n $.

\begin{example}\label{1b1}
Introduce an infinite sequence of random signs $ \tau_1, \tau_2, \dots
$; that is,
\begin{gather*}
\Pr{ \tau_k = -1 } = \Pr{ \tau_k = +1 } = \frac12 \quad \text{for each
 $ k $,} \\
\tau_1, \tau_2, \dots \quad \text{are independent.}
\end{gather*}
For each $ n $ we define a stochastic process $ X_n(\cdot) $, driven
by $ \tau_1, \dots, \tau_n $, as follows:
\[
\begin{gathered} X_n (t) = \prod_{k:1/n\le 1/k\le t} \tau_k \,
. \end{gathered} \qquad
\begin{gathered}\includegraphics{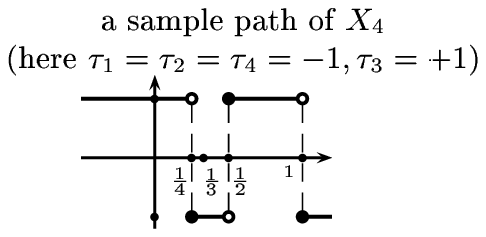}\end{gathered}
\]
For $ n \to \infty $, finite-dimensional distributions of $ X_n $
converge to those of a process $ X(\cdot) $. Namely, $ X $ consists of
countably many random signs, situated on intervals $ [\frac1{k+1},
\frac1k) $. Almost surely, $ X $ has no limit at $ 0+ $. We have
\begin{equation}\label{1b2}
\frac{ X(t) }{ X(s) } = \prod_{k:s<1/k\le t} \tau_k
\end{equation}
whenever $ 0 < s < t < \infty $. However, \eqref{1b2} does not hold
when $ s <
0 < t $. Here, the product contains infinitely many factors and
diverges almost surely; nevertheless, the increment $ X(t) / X(s) $ is
well-defined. Each $ X_n $ satisfies \eqref{1b2} for all $ s,t $
(including $ s < 0 < t $; of course, $ k \le n $), but $ X $ does
not. Still, $ X $ is an
independent increment process (multiplicatively); that is, $ X(t_2)
/X(t_1), \dots, X(t_n) /X(t_{n-1}) $ are independent whenever $
-\infty < t_1 < \dots < t_n < \infty $. However, we cannot describe
the whole $ X $ by a countable collection of its independent
increments. The infinite sequence of $ \tau_k = X(\frac1k+) /
X(\frac1k-) $ does not suffice since, say, $ X(1) $ is independent of
$ (\tau_1,\tau_2,\dots) $. Indeed, the global sign change $ x(\cdot)
\mapsto -x(\cdot) $ is a measure-preserving transformation that leaves
all $ \tau_k $ invariant. The conditional distribution of $ X(\cdot) $
given $ \tau_1, \tau_2, \dots $ is concentrated at two functions of
opposite global sign. It may seem that we should add to $ (\tau_1,
\tau_2, \dots) $ one more random sign $ \tau_\infty $ independent of $
(\tau_1, \tau_2, \dots) $ such that $ X(\frac1k) $ is a measurable
function of $ \tau_k, \tau_{k+1}, \dots $ and $ \tau_\infty
$. However, it is impossible. Indeed, $ X(1) = \tau_1 \dots \tau_k
X(\frac1k) $. Assuming $ X(\frac1k) = f_k (\tau_k, \tau_{k+1}, \dots;
\tau_\infty) $ we get $ f_1 ( \tau_1, \tau_2, \dots; \tau_\infty ) =
\tau_1 \dots \tau_{k-1} f_k (\tau_k, \tau_{k+1}, \dots; \tau_\infty )
$ for all $ k $. It follows that $ f_1 (\tau_1, \tau_2, \dots;
\tau_\infty) $ is orthogonal to all functions of the form $ g (\tau_1,
\dots, \tau_n ) h (\tau_\infty) $ for all $ n $, and thus, to a dense
(in $ L_2 $) set of functions of $ \tau_1, \tau_2, \dots; \tau_\infty
$; a contradiction.

So, for each $ n $ the process $ X_n $ is driven by $ (\tau_k) $, but
the limiting process $ X $ is not.

\end{example}

\begin{example}\label{1b3}

\begin{sloppypar}
(See also \cite{ES}.)
We turn to the other, the `continuous' model. For any $ \eps \in (0,1) $
we introduce a (complex-valued) stochastic process
\[
Y_\eps (t) = \begin{cases}
 \exp \( \I B(\ln t) - \I B(\ln\eps) \) &\text{for $ t \ge \eps $}, \\
 1 &\text{otherwise},
\end{cases}
\]
where $ B(\cdot) $ is the usual Brownian motion; or rather, $ \( B(t)
\)_{t\in[0,\infty)} $ and $ \( B(-t) \)_{t\in[0,\infty)} $ are two
independent copies of the usual Brownian motion. Multiplicative
increments $ Y_\eps (t_2) / Y_\eps (t_1), \dots, Y_\eps(t_n) /Y_\eps
(t_{n-1}) $ are independent whenever $ -\infty < t_1 < \dots < t_n <
\infty $, and the distribution of $ Y_\eps (t) /Y_\eps (s) $ does not
depend on $ \eps $ as far as $ \eps < s < t $ (in fact, the
distribution depends on $ t/s $ only). The distribution of $ Y_\eps
(1) $ converges for $ \eps \to 0 $ to the uniform distribution on the
circle $ |z|=1 $. The same for each $ Y_\eps (t) $. It follows easily
that, when $ \eps \to 0 $, finite dimensional distributions of $
Y_\eps $ converge to those of some process $ Y $. For every $ t > 0 $,
$ Y(t) $ is distributed uniformly on the circle; $ Y $ is an
independent increment process (multiplicatively), and $ Y(t) = 1 $ for
$ t \le 0 $. Almost surely, $ Y(\cdot) $ is continuous on $ (0,\infty)
$, but has no limit at $ 0+ $. We may define $ B(\cdot) $ by
\begin{gather*}
Y(t) = Y(1) \exp \( \I B(\ln t) \) \quad \text{for } t \in \R \, , \\
B(\cdot) \quad \text{is continuous on $ \R $} \, .
\end{gather*}
Then $ B $ is the usual Brownian motion, and
\[
\frac{ Y(t) }{ Y(s) } = \frac{ \exp(\I B(\ln t)) }{ \exp(\I B(\ln s)) }
\quad \text{for } 0 < s < t < \infty \, .
\]
However, $ Y(1) $ is independent of $ B(\cdot) $. Indeed, the global
phase change $ y(\cdot) \mapsto e^{i\al} y(\cdot) $ is a measure
preserving transformation that leaves $ B(\cdot) $ invariant. The
conditional distribution of $ Y(\cdot) $ given $ B(\cdot) $ is
concentrated on a continuum of functions that differ by a global phase
(distributed uniformly on the circle). Similarly to the
`discrete' example, we cannot introduce a random variable $ B(-\infty)
$ independent of $ B(\cdot) $, such that $ Y(t) $ is a function of $
B(-\infty) $ and increments of $ B(r) $ for $ -\infty < r < \ln t $.
\end{sloppypar}

So, for each $ \eps $, the process $ Y_\eps $ is driven by the
Brownian motion, but the limiting process $ Y $ is not.

\end{example}

Both toy models are singular at a given instant $ t = 0 $.
Interestingly, continuous stationary processes can demonstrate such 
strange behavior, distributed in time! (See Sects.\ \ref{sec:4},
\ref{sec:7}).

\subsection{Our limiting procedures}
\label{sec:1.3}

Imagine a sequence of elementary probabilistic models such that the $
n $-th model is driven by a finite sequence $ (\tau_1,\dots,\tau_n) $
of random signs (independent, as before). A limiting procedure may
lead to a model driven by an infinite sequence $ (\tau_1,\tau_2,\dots)
$ of random signs. However, it may also lead to something else, as
shown in \ref{sec:1.2}. This is an opportunity to ask ourselves:
what do we mean by a limiting procedure?

The $ n $-th model is naturally described by the finite probability
space $ \Om_n = \{ -1, +1 \}^n $ with the uniform measure. A
prerequisite to any limiting procedure is some structure able to join
these $ \Om_n $ somehow. It may be a sequence of
`observables',\index{observable} that is, functions on the disjoint
union,
\[
f_k : ( \Om_1 \uplus \Om_2 \uplus \dots ) \to \R \, .
\]

\begin{example}\label{1c1}

Let $ f_k (\tau_1,\dots,\tau_n) = \tau_k $ for $ n \ge k
$. Though $ f_k $ is defined only on $ \Om_k \uplus \Om_{k+1} \uplus
\dots $, it is enough. For every $ k $, the joint distribution of
$ f_1, \dots, f_k $ on $ \Om_n $ has a limit for $ n \to \infty $
(moreover, the distribution does not depend on $ n $, as far as $ n
\ge k $). The
limiting procedure should extend each $ f_k $ to a new probability
space $ \Om $ such that the joint distribution of $ f_1, \dots, f_k $
on $ \Om_n $ converges for $ n \to \infty $ to their joint
distribution on $ \Om $. Clearly, we may take the space of infinite
sequences $ \Om = \{ -1, +1 \}^\infty $ with the product measure, and
let $ f_k $ be the $ k $-th coordinate function.

\end{example}

\begin{example}

Still $ f_k (\tau_1,\dots,\tau_n) = \tau_k $ (for $ n
\ge k \ge 1 $), but in addition, the product $ f_0 (\tau_1,\dots,\tau_n) =
\tau_1 \dots \tau_n $ is included. For every $ k $, the joint
distribution of $ f_0, f_1, \dots, f_k $ on $ \Om_n $ has a limit for
$ n \to \infty $; in fact, the distribution does not depend on $ n $,
as far as $ n > k $ (this time, not just $ n \ge k $). Thus, in the
limit, $ f_0, f_1, f_2,
\dots $ become independent random signs. The functional dependence $
f_0 = f_1 f_2 \dots $ holds for each $ n $, but disappears in the
limit. We still may take $ \Om = \{-1,+1\}^\infty $, however, $ f_0 $
becomes a new coordinate.

\end{example}

This is instructive; the limiting model depends on the class of
`observables'.

\begin{example}\label{1c3}

Let $ f_k (\tau_1,\dots,\tau_n) = \tau_k \dots \tau_n $ for
$ n \ge k \ge 1 $. In the limit, $ f_k $ become
independent random signs. We may define $ \tau_k $ in the limiting
model by $ \tau_k = f_k / f_{k+1} $; however, we cannot express $ f_k
$ in terms of $ \tau_k $. Clearly, it is the same as the `discrete'
toy model of \ref{sec:1.2}.

\end{example}

The second and third examples are isomorphic. Indeed, renaming $ f_k $
of the third example as $ g_k $ (and retaining $ f_k $ of the second
example) we have
\[
g_k = \frac{ f_0 }{ f_1 \dots f_{k-1} } \, ; \qquad
f_k = \frac{ g_k }{ g_{k+1} } \text{ for } k > 0 \, , \quad \text{and}
\quad f_0 = g_1 \, ;
\]
these relations hold for every $ n $ (provided that the same $ \Om_n =
\{-1,+1\}^n $ is used for both examples) and naturally, give us an
isomorphism between the two limiting models.

That is also instructive; some changes of the class of `observables'
are essential, some are not.

It means that the sequence $ (f_k) $ is not really the structure
responsible for the limiting procedure. Rather, $ f_k $ are generators
of the relevant structure. The second and third examples differ only
by the choice of generators for the same structure. In contrast, the
first example uses a different structure. So, what is the mysterious
structure?

I can describe the structure in two equivalent ways. Here is the first
description. In the commutative Banach algebra $ l_\infty ( \Om_1
\uplus \Om_2 \uplus \dots ) $ of all bounded functions on the disjoint
union, we select a subset $ C $ (its elements will be called
observables) such that
\begin{equation}\label{1cf}
C \text{ is a separable closed subalgebra of } l_\infty ( \Om_1 \uplus
\Om_2 \uplus \dots ) \text{ containing the unit.}
\end{equation}
In other words,
\begin{equation}\label{1cg}
\begin{gathered}
C \text{ contains a sequence dense in the uniform topology;} \\
f_n \in C, \, f_n \to f \text{ uniformly} \imply f \in C \, ; \\
f,g \in C, \, a,b \in \R \imply af+bg \in C \, ; \\
\One \in C \, ; \\
f,g \in C \imply fg \in C
\end{gathered}
\end{equation}
(here $ \One $ stands for the unity, $ \One(\om) = 1 $ for all $ \om
$). Or equivalently,
\begin{equation}\label{1ch}
\begin{gathered}
C \text{ contains a sequence dense in the uniform topology;} \\
f_n \in C, \, f_n \to f \text{ uniformly} \imply f \in C \, ; \\
f,g \in C, \, \phi : \R^2 \to \R \text{ continuous} \imply \phi(f,g)
\in C \, . \\
\end{gathered}
\end{equation}
Indeed, on one hand, both $ af+bg $ and $ fg $ (and $ \One $) are
special cases of $ \phi(f,g) $. On the other hand, every continuous
function on a bounded subset of $ \R^2 $ can be uniformly approximated
by polynomials. The same holds for $ \phi (f_1,\dots,f_n) $ where $
f_1,\dots,f_n \in C $, and $ \phi : \R^n \to \R $ is a continuous
function. Another equivalent set of conditions is also well-known:
\begin{equation}\label{1ci}
\begin{gathered}
C \text{ contains a sequence dense in the uniform topology;} \\
f_n \in C, \, f_n \to f \text{ uniformly} \imply f \in C \, ; \\
f,g \in C, \, a,b \in \R \imply af+bg \in C \, ; \\
\One \in C \, ; \\
f \in C \imply |f| \in C \, ;
\end{gathered}
\end{equation}
here $ |f| $ is the pointwise absolute value, $ |f| (\om) = |f(\om)|
$.

The smallest set $ C $ satisfying these (equivalent) conditions
\eqref{1cf}--\eqref{1ci} and containing all given functions $ f_k $
is, by definition, generated by these $ f_k $.

Recall that $ C $ consists of functions defined on the disjoint union
of finite probability spaces $ \Om_n $; a probability measure $ P_n $
is given on each $ \Om_n $. The following condition is relevant:
\begin{equation}\label{condC}
\lim_{n\to\infty} \int_{\Om_n} f \, dP_n \text{ exists for every } f
\in C \, .
\end{equation}
Assume that $ C $ is generated by given functions $ f_k $. Then the
property \eqref{condC} of $ C $ is equivalent to such a property of
functions $ f_k $:
\begin{equation}\label{condf}
\parbox{10cm}{%
  For each $ k $, the joint distribution of $ f_1,\dots,f_k $ on $
  \Om_n $ weakly converges, when $ n \to \infty $.}
\end{equation}
Proof: \eqref{condf} means convergence of $ \int \phi(f_1,\dots,f_k)
\, dP_n $ for every continuous function $ \phi : \R^k \to \R
$. However, functions of the form $ f = \phi(f_1,\dots,f_k) $ (for all
$ k,\phi $) belong to $ C $ and are dense in $ C $.

We see that \eqref{condf} does not depend on the choice of generators
$ f_k $ of a given $ C $.

The second (equivalent) description of our structure is the `joint
compactification'\index{joint compactification} of $ \Om_1, \Om_2,
\dots $ I mean a pair $ (K,\al) $ such that
\begin{equation}\label{1cm}
\begin{gathered}
K \text{ is a metrizable compact topological space,} \\
\al : (\Om_1 \uplus \Om_2 \uplus \dots) \to K \text{ is a map,} \\
\text{the image } \al (\Om_1 \uplus \Om_2 \uplus \dots) \text{ is
dense in } K .
\end{gathered}
\end{equation}
Every joint compactification $ (K,\al) $ determines a set $ C $
satisfying \eqref{1cf}. Namely,
\[
C = \al^{-1} \( C(K) \) \, ;
\]
that is, observables $ f \in C $ are, by definition, functions of the
form
\[
f = g \circ \al, \text{ that is, } f(\om) = g(\al(\om)), \quad g \in
C(K) \, .
\]
The Banach algebra $ C $ is basically the same as the Banach algebra
$ C(K) $ of all continuous functions on $ K $.

Every $ C $ satisfying \eqref{1cf} corresponds to some joint
compactification. Proof: $ C $ is generated by some $ f_k $ such that
$ |f_k(\om)| \le 1 $ for all $ k,\om $. We introduce
\[
\begin{gathered}
\al(\om) = \( f_1(\om), f_2(\om), \dots \) \in [-1,1]^\infty \, , \\
K \text{ is the closure of } \al (\Om_1 \uplus \Om_2 \uplus \dots)
\text{ in } [-1,1]^\infty \, ;
\end{gathered}
\]
clearly, $ (K,\al) $ is a joint compactification. Coordinate functions
on $ K $ generate $ C(K) $, therefore $ f_k $ generate $ \al^{-1} \(
C(K) \) $, hence $ \al^{-1} \( C(K) \) = C $.

Finiteness of each $ \Om_n $ is not essential. The same holds for
arbitrary probability spaces $ (\Om_n,\F_n,P_n) $. Of course, instead
of $ l_\infty ( \Om_1 \uplus \Om_2 \uplus \dots ) $ we use $ L_\infty
( \Om_1 \uplus \Om_2 \uplus \dots ) $, and the map $ \al : (\Om_1
\uplus \Om_2 \uplus \dots) \to K $ must be measurable. It sends the
given measure $ P_n $ on $ \Om_n $ into a measure $ \al(P_n) $
(denoted also by $ P_n \circ \al^{-1} $) on $ K $. If measures $
\al(P_n) $ weakly converge, we get the limiting model $ (\Om,P) $ by
taking $ \Om = K $ and $ P = \lim_{n\to\infty} \al(P_n) $.

\subsection{Examples of high symmetry}
\label{sec:1.4}

\begin{example}\label{1d1}

Let $ \Om_n $ be the set of all permutations $ \om :
\{1,\dots,n\} \to \{1,\dots,n\} $, each permutation having the same
probability ($ 1/n! $);
\begin{gather*}
f : ( \Om_1 \uplus \Om_2 \uplus \dots ) \to \R \text{ is defined by}
 \\
f(\om) = | \{ k : \om(k) = k \} | \, ;
\end{gather*}
that is, the number of fixed points of a random permutation. Though $
f $ is not bounded, which happens quite often, in order to embed it
into the framework of \ref{sec:1.3}, we make it bounded by some
homeomorphism
from $ \R $ to a bounded interval (say, $ \om \mapsto \arctan f(\om)
$). The distribution of $ f(\cdot) $ on $ \Om_n $ converges (for $ n
\to \infty $) to the Poisson distribution $ P(1) $. Thus, the limiting
model exists; however, it is scanty: just $ P(1) $.

We may enrich the model by introducing
\[
f_u (\om) = | \{ k < un : \om(k) = k \} | \, ;
\]
for instance, $ f_{0.5} (\cdot) $ is the number of fixed points among
the first half of $ \{ 1, \dots, n \} $. The parameter $ u $ could run
over $ [0,1] $, but we need a countable set of functions; thus we
restrict $ u $ to, say, rational points of $ [0,1] $. Now the limiting
model is the Poisson process.

Each finite model here is invariant under permutations. Functions $
f_u $ seem to break the invariance, but the latter survives in their
increments, and turns in the limit into invariance of the Poisson
process (or rather, its derivative, the point process) under all
measure preserving transformations of $ [0,1] $.

Note also that \emph{independent} increments in the limit emerge from
\emph{dependent} increments in finite models.

We feel that all these $ f_u (\cdot) $ catch only a small part of the
information contained in the permutation. You may think about more
information, say, cycles of length $ 1,2,\dots $ (and what about
length $ n/2 \, $?)

\end{example}

\begin{example}\label{1d2}

Let $ \Om_n $ be the set of all graphs over $ \{ 1,\dots,n \}
$. That is, each $ \om \in \Om_n $ is a subset of the set $
\binom{\{1,\dots,n\}}{2} $ of all unordered pairs (treated as edges,
while $ 1,\dots,n $ are vertices); the probability of $ \om $ is $
p_n^{|\om|} (1-p_n)^{n(n-1)/2-|\om|} $, where $ |\om| $ is the number
of edges. That is, every edge is present with probability $ p_n $,
independently of others. Define $ f(\om) $ as the number of isolated
vertices. The limiting model exists if (and only if) there exists a
limit $ \lim_n n(1-p_n)^{n-1} = \la \in [0,\infty) $;\footnote{%
 Formally, the limiting model exists also for $ \la=\infty $, since
 the range of $ f $ is compactified.}
the Poisson distribution $ P(\la) $ exhausts the limiting model.

A Poisson process may be obtained in the same way as before.

You may also count small connected components which are more
complicated than single points.

Note that the finite model contains a lot of independence (namely, $
n(n-1)/2 $ independent random variables); the limiting model (Poisson
process) also contains a lot of independence (namely, independent
increments). However, we feel that independence is not inherited;
rather, the independence of finite models is lost in the
limiting procedure, and a new independence emerges.

\end{example}

\begin{example}\label{1d3}

Let $ \Om_n = \{-1,+1\}^n $ with uniform measure, and $ f_n : (
\Om_1 \uplus \Om_2 \uplus \dots ) \to \R $ be defined by
\[
f_u (\om) = \frac1{\sqrt n} \sum_{k<un} \tau_k (\om) \, ;
\]
as before, $ \tau_1, \dots, \tau_n $ are the coordinates, that is, $
\om = \( \tau_1(\om), \dots, \tau_n (\om) \) $ and $ u $ runs over
rational points of $ [0,1] $. The limiting model is the Brownian
motion, of course.

Similarly to \ref{1d1}, each finite model is invariant under
permutations. The invariance survives in increments of functions $ f_k
$, and in the limit, the white noise (the derivative of the Brownian
motion) is invariant under all measure preserving transformations of $
[0,1] $.

\end{example}

A general argument of \ref{sec:6.3} will show that a high
symmetry model cannot lead to a nonclassical scaling limit.

\subsection{Example of low symmetry}
\label{sec:1.5}

Example \ref{1d3} may be rewritten via the composition of random maps
\[
\begin{gathered}
 \al_-, \al_+ : \Z \to \Z \, , \\
 \al_- (k) = k-1 \, , \quad \al_+ (k) = k+1 \, ; \\
 \al_\om = \al_{\tau_n(\om)} \circ \dots \al_{\tau_1(\om)} \, ;
\end{gathered} \qquad
\begin{gathered}\includegraphics{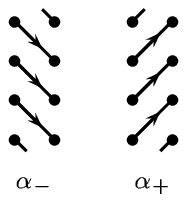}\end{gathered}
\]
thus, $ \al_\om (k) = k + \tau_1(\om) + \dots + \tau_n(\om) $, and we
may define $ f_1 (\om) = \frac1{\sqrt n} \al_\om(0) $, which conforms
to \ref{1d3}. Similarly, $ f_u (\om) = \frac1{\sqrt n} \al_{\om,u} (0)
$, where $ \al_{\om,u} $ is the composition of $ \al_{\tau_k(\om)} $
for $ k \le un $. The order does not matter, since $ \al_-, \al_+ $
commute, that is, $ \al_- \circ \al_+ = \al_+ \circ \al_- $. It is
interesting to try a pair of noncommuting maps.

\begin{example}\label{1e1}

(See Warren \cite{Wa1}.) Define
\[
\begin{gathered}
 \al_-, \al_+ : \Z+\frac12 \to \Z+\frac12 \, , \\
 \begin{gathered} \al_- (x) = x-1 \, , \\ \al_+ (x) = x+1 \end{gathered}
 \quad \text{for } x \in \(
 \Z + \tfrac12 \) \cap (0,\infty) \, , \\ 
 \al_- (-x) = - \al_- (x) \, , \quad \al_+ (-x) = - \al_+ (x) \, .
\end{gathered}
\qquad
\begin{gathered}\includegraphics{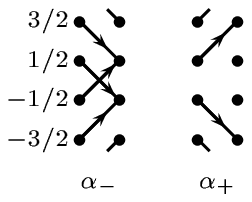}\end{gathered}
\]
These are not invertible functions; $ \al_- $ is not injective, $
\al_+ $ is not surjective. Well, we do not need to invert them, but
need their compositions:
\[
\begin{gathered}
\al_\om = \al_{\tau_n(\om)} \circ \dots \circ \al_{\tau_1(\om)} \, .
\end{gathered}
\qquad
\begin{gathered}\includegraphics{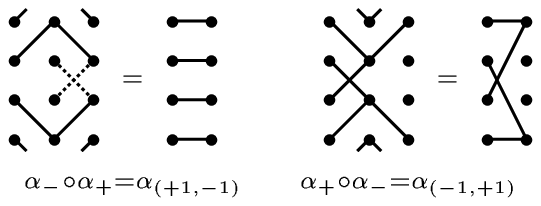}\end{gathered}
\]
All compositions belong to a two-parameter set of functions $ h_{a,b}
$,
\begin{gather*}
\begin{gathered}
\al_\om (x) = h_{a,b} (x) = \begin{cases}
 x+a &\text{for $ x \ge b $}, \\
 x-a &\text{for $ x \le -b $}, \\
 (-1)^{b-x} (a+b) &\text{for $ -b \le x \le b $};
 \end{cases}
\end{gathered}
\;
\begin{gathered}\includegraphics{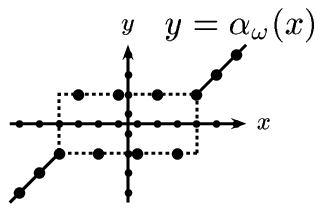}\end{gathered}
\\
b, a+b \in \( \Z + \tfrac12 \) \cap (0,\infty) = \{ \tfrac12,
 \tfrac32, \tfrac52, \dots \} \, .
\end{gather*}
Indeed, $ \al_- = h_{-1,1.5} $, $ \al_+ = h_{1,0.5} $, and $
h_{a_2,b_2} \circ h_{a_1,b_1} = h_{a,b} $ where $ a = a_1 +
a_2 $, $ b = \max ( b_1, b_2-a_1 ) $. Thus, $ \al_\om =
h_{\al(\om),b(\om)} $, and we define
\begin{gather*}
f_1 : ( \Om_1 \uplus \Om_2 \uplus \dots ) \to \R^2 \times \{-1,+1\} \,
 , \\
f_1 (\om) = \bigg( \frac{ a(\om) }{ \sqrt n }, \frac{ b(\om) }{ \sqrt
 n }, (-1)^{b(\om)-0.5} \bigg) \, .
\end{gather*}
However, the function is neither bounded nor real-valued; in order to
fit into the framework of \ref{sec:1.3} we take, say, $ \arctan
\( a(\om) / \sqrt
n \) $, $ \arctan \( b(\om) / \sqrt n \) $, and $ (-1)^{b(\om)-0.5}
$. The latter is essential if, say, $ \frac1{\sqrt n} \al_\om (0.5) $
is treated as an `observable'; indeed, $ \frac1{\sqrt n} \al_\om
(0.5) = (-1)^{b(\om)-0.5} \frac1{\sqrt n} (a(\om)+b(\om)) $.
The limiting model exists, and is quite interesting. (See also
\ref{sec:8.3}.) As before, a random process appears by considering the
composition over $ k < un $.

Here, finite models are not invariant under permutations of their
independent random variables (since the maps do not commute), and the
limiting model appears not to be invariant under measure preserving
transformations of $ [0,1] $.

Independence present in finite models survives in the limit, provided
that the limit is described by a two-parameter random process; we'll
return to this point in \ref{sec:4.3}.

\end{example}

\subsection{Trees, not cubes}
\label{sec:1.6}

\begin{example}\label{1f1}
A particle moves on the sphere $ S^2 $. Initially it is at a given
point $ x_0 \in S^2 $. Then it jumps by $ \eps $ in a random
direction. That is, $ X_0 = x_0 $, while the next random variable $
X_1 $ is distributed uniformly on the circle $ \{ x \in S^2 :
|x_0-x| = \eps \} $. Then it jumps again to $ X_2 $ such that $
|X_1-X_2| = \eps $, and so on. We have a Markov chain $ (X_k) $
in discrete time (and continuous space). Let $ \Om_\eps $ be the
corresponding probability space; it may be the space of sequences $
(x_0,x_1,x_2,\dots) $ satisfying $ |x_k-x_{k+1}| = \eps $, or
something else, but in any case $ X_k : \Om_\eps \to S^2 $. We choose
$ \eps_n \to 0 $ (say, $ \eps_n = 1/n $), take $ \Om_n = \Om_{\eps_n}
$ and define $ f_u : ( \Om_1 \uplus \Om_2 \uplus \dots ) \to S^2 $ by
\[
f_u (\om) = X_k(\om) \quad \text{for } \eps_n^2 k \le u < \eps_n^2
(k+1) \, , \quad \om \in \Om_n \, .
\]
Of course, the limiting model is the Brownian motion on the sphere $
S^2 $.

In contrast to previous examples, here $ \Om_n $ is not a product; the
$ n $-th model does not consist of independent random
variables. But, though we can parameterize these Markov transitions by
independent random variables, there is a lot of freedom in
doing so; none of the parameterizations may be called
canonical. The same holds for the limiting model. The Brownian motion
on $ S^2 $ can be driven by the Brownian motion on $ R^2 $ according
to some stochastic differential equation, but the latter involves a lot
of freedom.

\end{example}

\begin{example}\label{1f2}

(See \cite{ST}.)
Consider the random walk on such an oriented graph:
\[
\begin{gathered}\includegraphics[scale=0.8]{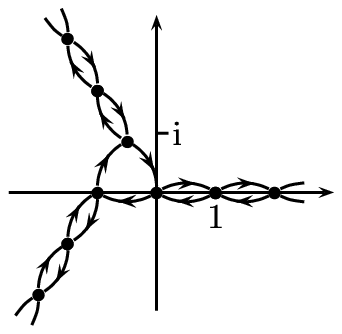}\end{gathered}
\]
A particle starts at $ 0 $ and chooses at random (with probabilities $
1/2 $, $ 1/2 $) one of the two outgoing edges, and so on (you see,
exactly two edges go out of any vertex). Such $ (Z_0,Z_1,\dots) $ is
known as the simplest spider walk. It is a complex-valued
martingale. The set $ \Om_n $ of all $ n $-step trajectories contains
$ 2^n $ elements and carries its natural structure of a binary
tree. (It can be mapped to the binary cube $ \{-1,+1\}^n $ in many
ways.) We define $ f_u : ( \Om_1 \uplus \Om_2 \uplus \dots ) \to \C $
by
\[
f_u (\om) = \frac1{\sqrt n} Z_k (\om) \quad \text{for } k \le nu < k+1
\, , \quad \om \in \Om_n \, .
\]
The limiting model is a continuous complex-valued martingale whose
values belong to the union of three rays.
\[
\begin{gathered}\includegraphics{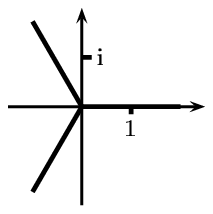}\end{gathered}
\]
The process is known as Walsh's Brownian motion, a special case of
the so-called spider martingale.

\end{example}

\subsection{Sub-\protect\sif s}
\label{sec:1.7}

Every example considered till now follows the pattern of
\ref{sec:1.3}; a joint
compactification of probability spaces $ \Om_n $, and the limiting $
\Om $. Moreover, $ \Om_n $ is usually related to a set $ T_n $ (a
parameter space, interpreted as time or space), and $ \Om $ to a joint
compactification $ T $ of these $ T_n $.
\[
\begin{array}{|c|c|c|}
\hline
\text{Example} & T_n & T \\
\hline
\text{\ref{1b1}} & \{ 1, \frac12, \dots, \frac1n \} &
 \{ 1, \frac12, \frac13, \dots \} \cup \{0\} \\
\text{\ref{1b3}} & [\eps_n,1] & [0,1] \\
\text{\ref{1d1}, \ref{1d2}, \ref{1d3}, \ref{1e1}, \ref{1f1},
\ref{1f2}} & \{ \frac1n, \frac2n, \dots, 1 \} & [0,1] \\
\hline
\end{array}
\]
Examples \ref{1b1}, \ref{1b3}, \ref{1d3} deal (for a finite $ n $)
with independent increment processes, taking on their values in a
group, namely, \ref{1d3}: $ \R $ (additive); \ref{1b1}: $ \{-1,+1\} $
(multiplicative), \ref{1b3}: the circle $ \{ z \in \C : |z|=1 \} $
(multiplicative). Every $ t \in T_n $ splits the process into two
parts, the past and the future; in order to keep them independent, we
define them via increments, not values.\footnote{%
 In fact, the process of \ref{1b1} has also independent values
 (not only increments); but that is irrelevant.}
In terms of random signs $ \tau_k $ (for \ref{1b1}, \ref{1d3}) it
means simply $ \{-1,+1\}^n = \{-1,+1\}^k \times \{-1,+1\}^{n-k} $;
here $ k $ depends on $ t $. The same idea (of independent parts) is
formalized by sub-\sif s $ \F_{0,t} $ (the past) and $ \F_{t,1} $ (the
future) on our probability space ($ \Om_m $ or $ \Om $). Say, for the
Brownian motion \ref{1d3}, $ \F_{0,t} $ is generated by Brownian
increments on $ [0,t] $, while $ \F_{t,1} $ --- on $ [t,1]
$. Similarly we may define $ \F_{s,t} $ for $ s < t $, and we have
\[
\F_{r,s} \otimes \F_{s,t} = \F_{r,t} \quad \text{whenever } r < s < t
\, . \index{zzz@$ \otimes $!for \sif s}
\]
It means two things: first, independence,\index{independence}
\[
\Pr{ A \cap B } = \Pr{A} \Pr{B} \quad \text{whenever } A \in \F_{r,s},
B \in \F_{s,t} \, ;
\]
and second, $ \F_{r,t} $ is generated by $ \F_{r,s} $ and $ \F_{s,t} $
(that is, $ \F_{r,t} $ is the least sub-\sif\ containing both $
\F_{r,s} $ and $ \F_{s,t} $). Such a two-parameter family $
(\F_{s,t}) $ of sub-\sif s is called a \emph{factorization} (of the
given probability space). Some additional precautions are needed when
dealing with semigroups (like \ref{1e1}), and also, with
discrete time.

Sub-\sif s $ \F_A $ can be defined for some subsets $ A \subset T $
more general than intervals, getting
\[
\F_A \otimes \F_B = \F_C \quad \text{whenever } A \uplus B = C \, .
\]
Models of high symmetry admit arbitrary measurable sets $ A $; models
of low symmetry do not. For some examples (such as \ref{1d1},
\ref{1d2}), a factorization emerges after the limiting
procedure.\footnote{%
 For \ref{1d2}, some factorization is naturally defined for $
 \Om_n $, but is lost in the limiting procedure, and a new
 factorization emerges.}

No factorization at all is given for \ref{1f1}, \ref{1f2}. Still, the
past $ \F_{0,t} = \F_t $ is defined naturally. However, the future is
not defined, since possible continuations depend on the past. Here we
deal with a one-parameter family $ (\F_t) $ of sub-\sif s, satisfying
only a monotonicity condition
\[
\F_s \subset \F_t \quad \text{whenever } s < t \, ;
\]
such $ (\F_t) $ is called a \emph{filtration.}

\section{Abstract Nonsense of the Scaling Limit}
\label{sec:2}
\subsection{More on our limiting procedures}
\label{sec:2.1}

The joint compactification $ K $ of $ \Om_1 \uplus \Om_2 \uplus \dots
$, used in \ref{sec:1.3}, is not quite satisfactory. Return to \ref{1d3}:
\begin{equation}\label{2a1}
f_u (\om) = \frac1{\sqrt n} \sum_{k<un} \tau_k(\om) \quad \text{for }
u \in [0,1] \cap \Q
\end{equation}
($ \Q $ being the set of rational numbers). The limiting model is the
Brownian motion, restricted to $ [0,1] \cap \Q $. What about an
irrational point, $ v \in [0,1] \setminus \Q \, $? The random variable
$ f_v $ may be defined on $ \Om $ as the limit (say, in $ L_2 $) of $
f_u $ for $ u \to v $, $ u \in [0,1] \cap \Q $. On the other hand, $
f_v $ is naturally defined on $ \Om_1 \uplus \Om_2 \uplus \dots $ (by
the same formula \eqref{2a1}). However, $ f_v $ is not a continuous
function on the compact space $ K $.\footnote{%
 There exist $ \om_n \in \Om_n $ such that $ \lim_n f_u (\om_n) $
 exists for all $ u \in [0,1] \cap \Q $, but $ \lim_n f_v (\om_n) $
 does not exist.
 \[
\begin{gathered}\includegraphics[scale=0.9]{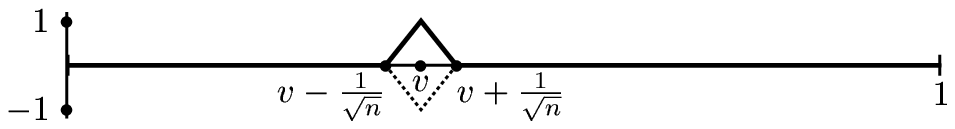}\end{gathered}
 \]}
Thus, the weak convergence $ P_i \to P $ is relevant to $ f_u $ but
not $ f_v $. Something is wrong!

What is wrong is the uniform topology used in
\eqref{1cf}--\eqref{1ci}. A right topology should take measures $ P_i
$ into account. We have two ways, `moderate' and `radical'.

Here is the `moderate' way. We choose some appropriate subsets $ B_m
\subset ( \Om_1 \uplus \Om_2 \uplus \dots ) $, $ B_1 \subset B_2
\subset \dots $, such that
\[
\inf_i P_i ( B_m \cap \Om_i ) \uparrow 1 \quad \text{for } m \to
\infty
\]
and in \eqref{1cg}--\eqref{1ci} replace the assumption ``$ f_n \in C
$, $ f_n \to f $ uniformly \imp $ f \in C $'' with
\begin{equation}\label{2a1a}
f_n \in C, \, f_n \to f \text{ uniformly on each $ B_m $} \imply f \in
C \, .
\end{equation}

\begin{example}\label{2a1b}

Continuing \eqref{2a1} we define $ B_m $ by
\[
B_m \cap \Om_i = \left\{ \om \in \Om_i : \sup_{0\le k<l\le i} \frac{
\left| \frac1{\sqrt i} \sum_{j=k}^l \tau_j (\om) \right| }{
\left( \frac{l-k}{i} \right)^{1/3} } \le m \right\} \, ;
\]
then\footnote{%
 Of course, $ |u-v|^\al $ for any $ \al \in (0,1/2) $ may be used, not
 only $ |u-v|^{1/3} $.}
\[
| f_u (\om) - f_v (\om) | \le m |u-v|^{1/3} \quad \text{for } \om \in
B_m \cap \Om_i
\]
if $ i $ is large enough (namely, $ 2/i < |u-v| $). The set $ C $
(satisfying \eqref{2a1a}) generated by $ f_u $ for all rational $ u $,
also contains $ f_v $ for all irrational $ v $.

\end{example}

Similarly to \ref{sec:1.3}, we may translate \eqref{2a1a} into
the topological
language. For each $ m $, the restriction of $ C $ to $ B_m $
corresponds to a joint compactification $ (K_m,\al_m) $ of $ B_m \cap
\Om_i $. Clearly, $ K_{m_1} \subset K_{m_2} $ for $ m_1 < m_2 $, and $
\al_{m_1} = \al_{m_2} |_{K_{m_1}} $. Thus, we get a \emph{joint $ \si
$-compactification}\index{joint s@joint $ \si $-compactification}
\[
\al : ( \Om_1 \uplus \Om_2 \uplus \dots ) \to K_\infty = K_1 \cup K_2
\cup \dots
\]
We do not need a topology on the union $ K_\infty $ of metrizable
compact spaces $ K_1 \subset K_2 \subset \dots $\footnote{%
 But if you want, $ K_\infty $ may be equipped with the inductive
 limit topology; that is, $ U \subset K_\infty $ is open if and only
 if for every $ m $, $ U \cap K_m $ is open (in $ K_m $). However, the
 topology usually is not metrizable.}
We just define $ C(K_\infty) $ as the set of all functions $ g :
K_\infty \to \R $ such that $ g |_{K_m} $ is continuous (on $ K_m $)
for each $ m $. We have
\[
C = \al^{-1} \( C(K_\infty) \) \, ,
\]
that is, observables $ f \in C $ are functions of the form
\[
f = g \circ \al \, , \quad \text{that is, } f(\om) = g ( \al(\om) ),
\quad g \in C(K_\infty) \, .
\]
If measures $ \al(P_i) $ weakly converge (w.r.t. bounded functions of
$ C(K_\infty) $, recall \eqref{condC}, \eqref{condf}), we
get the limiting model $ (\Om,P) $ by taking $ \Om = K_\infty $ and $
P = \lim_{i\to\infty} \al(P_i) $.

\begin{example}

Continuing \ref{2a1b} we see that the limiting measure $ P $ exists,
and the joint distribution of all $ f_u $ (extended to $ K_\infty $ by
continuity) w.r.t.\ $ P $ is the Wiener measure. The `uniform' metric
on $ K_\infty $,
\[
\dist (x,y) = \sup_{0\le u\le1} | f_u(x) - f_u(y) | \, ,
\]
is continuous on each $ K_m $ (intersected with the support of $ P
$). Therefore, every function continuous in
the `uniform' metric belongs to $ C(K_\infty) $. Our joint $ \si
$-compactification is another form of the usual weak convergence of
random walks to the Brownian motion.

\end{example}

That was the `moderate way'. It requires special subsets $ B_m \subset
( \Om_1 \uplus \Om_2 \uplus \dots ) $, in contrast to the `radical
way'; basically, the latter allows the sequence of sets $ B_m $ to
depend on a sequence of functions $ f_n $, see \eqref{2a1a}. In other
words, instead of uniform (or `locally uniform') convergence, we
introduce a weaker topology by the metric\footnote{%
 Alternatively, we may restrict ourselves to bounded functions $ \Om_1
 \uplus \Om_2 \uplus \dots \to [-1,+1] $ (applying a transformation
 like $ \arctan $) and use, say,
 \[
 \dist (f,g) = \sup_i \int | f(\om)-g(\om) | \, \D P_i (\om) \, .
 \]}
\begin{equation}\label{2a2}
\dist (f,g) = \sup_i \int \frac{ | f(\om)-g(\om) | }{ 1 + |
f(\om)-g(\om) | } \, \D P_i (\om) \, .
\end{equation}
If $ f_n \in C(K) $ and $ \dist (f_n,f) \to 0 $ then $ f_n $ converge
in probability w.r.t.\ $ P $; thus, $ f $ is naturally defined $ P
$-almost everywhere.\footnote{%
 In fact, every (equivalence class of) $ P $-measurable function can
 be obtained in that way provided that, for each $ i $, supports of $
 P_i $ and $ P $ do not intersect. It means that every random variable
 on the limiting probability space is the scaling limit of some
 function on $ \Om_1 \uplus \Om_2 \uplus \dots $ (see also \ref{2c8}).}

Let $ C $ be the closure of $ C(K) $ in the metric \eqref{2a2}. Then
\[
\int \phi (f_1,\dots,f_d) \, \D P_i \tendsi \int \phi (f_1,\dots,f_d)
\, \D P
\]
for every $ d $, every bounded continuous function $ \phi : \R^d \to
\R $, and every $ f_1,\dots,f_d \in C $. The joint
distribution of $ f_1, \dots, f_d $ w.r.t.\ $ P_i $ converges (weakly)
to that w.r.t.\ $ P $. So, the weak convergence $ P_i \to P $ is
relevant for the whole $ C $ (not only $ C(K) $). That is the idea of
the `radical way', presented systematically in \ref{sec:2.2},
\ref{sec:2.3}.

Returning again to \ref{1d3} we see that $ f_v $ (for $ v \in
[0,1] $) is the limit of $ f_u $ (for $ u \in [0,1] \cap \Q $) in the
metric \eqref{2a2}; thus, $ f_v \in C $ for all $ v \in [0,1] $.

However, much more can be said. Not only
\[
\Lim_{i\to\infty} \bigg( \frac1{\sqrt i} \sum_{ai < k < bi} \tau_k
(\om) \bigg) = \int_a^b \D B(t) \, ,
\]
where `$ \Lim $' means the scaling limit (as explained above), but
also
\begin{multline*}
\Lim_{i\to\infty} \bigg( i^{-d/2} \sum_{ai<k_1<\dots<k_d<bi}
\tau_{k_1} (\om)  \dots\tau_{k_d} (\om) \bigg) \\
= \idotsint\limits_{a<t_1<\dots<t_d<b} \D B(t_1) \dots \D B(t_d) =
 \frac1{d!} H_d \( B(b)-B(a), b-a \)
\end{multline*}
where $ H_d $ is the Hermite polynomial (see for instance
\cite[IV.3.8]{RY}). Taking finite linear combinations and their
closure in the metric \eqref{2a2} we get
\begin{multline}\label{2a3}
\Lim_{i\to\infty} \bigg( \sum_{d=0}^\infty i^{-d/2}
\sum_{0<k_1<\dots<k_d<i} \psi_d \( \tfrac{k_1}i, \dots, \tfrac{k_d}i
 \) \tau_{k_1} (\om) \dots \tau_{k_d}(\om) \bigg) \\
= \sum_{d=0}^\infty \;\; \idotsint\limits_{0<t_1<\dots<t_d<1} \psi_d
 (t_1,\dots,t_d) \, \D B(t_1) \dots \D B(t_d)
\end{multline}
provided that functions $ \psi_d $ are Riemann integrable, and vanish
for $ d $ large enough. The right-hand side is well-defined for all $
\psi_d \in L_2 $ such that $ \sum_d \| \psi_d \|_2^2 < \infty $; the
scaling limit may be kept by replacing $ \psi_d \( \tfrac{k_1}i,
\dots, \tfrac{k_d}i \) $ with the mean value of $ \psi_d $ on the $
1/i $-cube centered at $ \( \tfrac{k_1}i, \dots, \tfrac{k_d}i \)
$. Now, $ (0,1) $ may be replaced with the whole $ \R $; $ \psi_d $ is
defined on $ \De_d = \{ (x_1,\dots,x_d) \in \R^d : x_1 < \dots < x_d
\} $. The right-hand side of \eqref{2a3} gives us an isometric linear
correspondence between $ L_2 ( \De_0 \uplus \De_1 \uplus \De_2 \uplus
\dots ) $ and $ L_2 (\Om,\F,P) $, where $ (\Om,\F,P) $ is the
probability space describing the Brownian motion (on the whole $ \R
$).

\subsection{Coarse probability space: definition and simple example}
\label{sec:2.2}

\begin{definition}\label{2b1}
A \emph{coarse probability space}\index{coarse!probability space} $ \(
(\Om[i],\F[i],P[i])_{i=1}^\infty, \A \) $\index{zzOm@$ \Om[i] $}
\index{zzF@$ \F[i] $}\index{zzP@$ P[i] $}
consists of a sequence of
probability spaces $ (\Om[i],\F[i],P[i]) $ and a set $ \A $ of subsets
of the disjoint union $ \Om [\all] = \Om(1) \uplus \Om(2) \uplus
\dots $,\index{zzup@$ \uplus $, disjoint union}
satisfying the following conditions:

\begin{myitemize}
\item[\textup{(a)}]
$ \forall A \in \A \; \forall i \; (A\cap\Om[i]) \in \F[i] $;

\item[\textup{(b)}]
$ \forall A,B \in \A \; \( A \cap B \in \A, \, A \cup B \in \A, \,
\Om[\all] \setminus A \in \A \) $;

\item[\textup{(c)}]
$ \A $ contains every $ A \subset \Om[\all] $ such that $
\forall i \; (A\cap\Om[i]) \in \F[i] $ and $ P[i] \( A \cap \Om[i] \)
\to 0 $ for $ i \to \infty $;

\item[\textup{(d)}]
$ \( \cup_{k=1}^\infty A_k \) \in \A $ for every pairwise disjoint
$ A_1, A_2, \dots \in \A $ such that $ \sum_k \sup_i P[i] \( A_k \cap
\Om[i] \) < \infty $;

\item[\textup{(e)}]
$ \lim_i P[i] \( A \cap \Om[i] \) $ exists for every $ A \in \A $;

\item[\textup{(f)}] there exists a finite or countable subset $ \A_1
\subset \A $ that \emph{generates} $ \A $ in the sense that the least
subset of $ \A $ satisfying \textup{(b)--(d)} and containing $ \A_1 $
is the whole $ \A $.

\end{myitemize}

A set $ \A $ satisfying \textup{(a)--(f)} will be called a
\emph{coarse \sif}\index{coarse!sigma@\sif}\footnote{%
 It is not a \sif, unless $ \A $ contains all sets satisfying
 \ref{2b1}(a).}
\textup{(}on the \emph{coarse sample space}\index{coarse!sample space} $
(\Om[i],\F[i],P[i])_{i=1}^\infty $\textup{)}. Each set $ A $ belonging
to the coarse \sif\ $ \A $ will be called \emph{coarsely
measurable}\index{coarsely measurable!set}
\textup{(}w.r.t.\ $ \A $\textup{)}, or a \emph{coarse
event.}\index{coarse!event}
\end{definition}

\begin{remark}
Condition \ref{2b1}(c) is equivalent to

\begin{myitemize}
\item[(c1)]
$ \forall i \; \F[i] \subset \A $.
That is, if a set $ A \subset \Om[\all] $ is contained in some $
\Om[i] $, and is $ \F[i] $-measurable, then $ A \in \A $.
\end{myitemize}

Also, Condition \ref{2b1}(d) is equivalent to each of the following
conditions (d1)--(d4). There, we assume that $ A \subset
\Om[\all] $, $ \forall i \; \( A \cap \Om[i] \) \in \F[i] $, and
$ \forall k \; A_k \in \A $.

\begin{myitemize}

\item[(d1)]
If $ A_k \uparrow A $ (that is, $ A_1 \subset A_2 \subset \dots $ and
$ A = \cup_k A_k $) and $ \sup_i P[i] \( (A\setminus A_k) \cap \Om[i]
\) \to 0 $ for $ k \to \infty $, then $ A \in \A $.

\item[(d2)]
If $ \sup_i P[i] \( (A\bigtriangleup A_k) \cap \Om[i] \) \to 0 $ for $
k \to \infty $, then $ A \in \A $. (Here $ A\bigtriangleup A_k =
(A\setminus A_k) \cup (A_k \setminus A) $.)

\item[(d3)]
If $ A_k \uparrow A $ and $ \limsup_i P[i] \( (A\setminus A_k) \cap
\Om[i] \) \to 0 $ for $ k \to \infty $, then $ A \in \A $.

\item[(d4)]
If $ \limsup_i P[i] \( (A\bigtriangleup A_k) \cap \Om[i] \) \to 0 $
for $ k \to \infty $, then $ A \in \A $.

\end{myitemize}
So, we have 10 equivalent combinations: (c)\&(d), (c1)\&(d),
(c)\&(d1), (c1)\&(d1), (c)\&(d2), \dots, (c1)\&(d4). (I omit the
proof.)

However, ``$ \sup_i $'' in (d) cannot be replaced with ``$ \limsup_i
$''.
\end{remark}

\begin{lemma}\label{2b3}
Let $ \A_1 $ be a finite or countable set satisfying
\textup{\ref{2b1}(a,e)} and

\textup{(b1)} $ \forall A,B \in \A_1 \; \( A \cap B \in \A_1 \) $.

Then the least set $ \A $ containing $ \A_1 $ and satisfying
\textup{\ref{2b1}(b,c,d)} is a coarse \sif.
\end{lemma}

\begin{proof}
The algebra generated by $ \A_1 $ satisfies (e), since $ P[i] \(
(A \cup B) \cap \Om[i] \) = P[i] ( A \cap \Om[i] ) + P[i] ( B \cap
\Om[i] ) - P[i] \( (A \cap B) \cap \Om[i] \) $. We enlarge the algebra
according to (c), which preserves (e), as well as (a), (b). Finally,
we enlarge it according to (d), which preserves (a), (b), (e); (c) and
(f) hold trivially.
\qqed\end{proof}

In such a case we say that the coarse \sif\ $ \A $ is
\emph{generated}\index{generated!coarse \sif, by sets} by the set $
\A_1 $.

\begin{example}
Let $ \Om[i] = \{ 0, \frac1i, \dots, \frac{i-1}i \} $, and $ P[i] $ be
the uniform distribution on $ \Om[i] $. Every interval $ (s,t) \subset
(0,1) $ gives us a set $ A_{s,t} \subset \Om[\all] $,
\[
A_{s,t} \cap \Om[i] = (s,t) \cap \Om[i] \, .
\qquad
\begin{gathered}\includegraphics{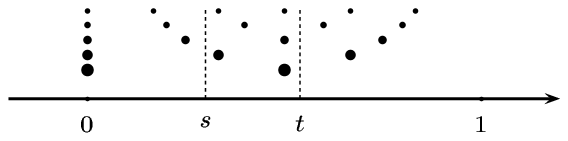}\end{gathered}
\]
We take a dense countable set of pairs $ (s,t) $ (say, rational $ s,t
$) and consider the set $ \A_1 $ of the corresponding $ A_{s,t} $. The
set $ \A_1 $ satisfies the conditions of \ref{2b3}, therefore it
generates a coarse \sif\ $ \A $. In fact, $ \A $ consists of all $ A =
A[1] \uplus A[2] \uplus \dots $ such that sets $ A[i]+(0,1/i) \subset
(0,1) $ converge in probability to some $ A[\infty] \subset (0,1) $;
that is, $ \mes \( A[\infty] \bigtriangleup ( A[i]+(0,1/i) ) \) \to 0
$ for $ i\to\infty $.
\[
\begin{gathered}\includegraphics{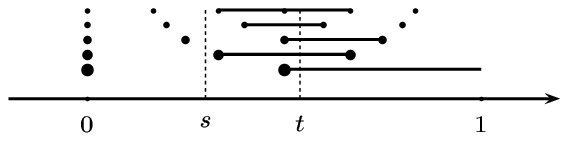}\end{gathered}
\]
If $ A = A_{s,t} $ then, of course, $ A[\infty] = (s,t) $.
\end{example}

\begin{example}
Continuing \ref{1c1}, we take $ \Om[i] = \{-1,+1\}^i $ with
the uniform distribution $ P[i] $. Given $ n $ and $ a =
(a_1,\dots,a_n) \in \{-1,+1\}^n $, we consider $ A_a \subset
\Om[\all] $,
\[
A_a \cap \Om[i] = \{ (\tau_1,\dots,\tau_i) : \tau_1 = a_1, \dots,
\tau_n = a_n \} \quad \text{for } i \ge n \, .
\]
Such sets $ A_a $ (for all $ a $ and $ n $) are a countable collection
$ \A_1 $ satisfying the conditions of \ref{2b3}, therefore it
generates a coarse \sif\ $ \A $. In fact, $ \A $ consists of all $ A =
A[1] \uplus A[2] \uplus \dots $ such that sets $ \be_i^{-1} (A)
\subset (0,1) $ converge in probability to some $ A[\infty] \subset
(0,1) $; here $ \be_i : (0,1) \to \Om[i] $ is such a measure
preserving map:
\[
\be_i (x) = \( (-1)^{c_1}, \dots, (-1)^{c_i} \) \quad \text{when
} x - \Big( \frac{c_1}{2} + \dots + \frac{c_i}{2^i} \Big) \in \Big( 0,
\frac1{2^i} \Big) \, ,
\]
for any $ c_1,\dots,c_i \in \{0,1\} $.
\end{example}

You may guess that some limiting procedure produces a (`true', not
coarse) probability space out of any given coarse probability
space. Indeed, such a procedure, called `refinement', is described in
\ref{sec:2.3}.

\subsection{Good use of joint compactification}
\label{sec:2.3}

Having a coarse probability space $ \(
(\Om[i],\F[i],P[i])_{i=1}^\infty, \A \) $ and its refinement $
(\Om,\F,P) $ (to be defined later), we may hope that the Hilbert space
$ L_2[\infty] = L_2 (\Om,\F,P) $ is in some sense the limit of Hilbert
spaces $ L_2[i] = L_2 \( \Om[i], \F[i], P[i] \) $. That is indeed the
case in the framework of joint compactification, as we'll see. A bad
use of the framework, tried in \ref{sec:1.3}, is a joint
compactification of given probability spaces. A good use, considered
here, is a joint compactification of metric (Hilbert, \dots) spaces
built over the given probability spaces.

\begin{definition}\label{2c1}
A \emph{coarse Polish space}\index{coarse!Polish space} is $ \(
(S[i],\rho[i])_{i=1}^\infty, c \) $, where each $ (S[i],\linebreak[0]
\rho[i]) $ is a Polish space\index{Polish space} (that is, a
complete separable metric space\footnote{%
 Many authors define a Polish space as a metrizable topological space
 admitting a complete separable metric. However, I assume that a
 metric is given.}%
), and $ c \subset S[1] \times S[2] \times \dots $ is a set of
sequences $ x = \( x[1], x[2], \dots \) $ satisfying the following
conditions:

\textup{(a)} if $ x_1, x_2 \in S[1] \times S[2] \times \dots $ are
such that $ \rho[i] \( x_1[i], x_2[i] \) \to 0 $ (for $ i \to \infty
$), then $ ( x_1 \in c ) \equi ( x_2 \in c ) $;

\textup{(b)} if $ x, x_1, x_2, \dots \in S[1] \times S[2] \times \dots
$ are such that $ \sup_i \rho[i] \( x_k[i], x[i] \) \to 0 $
\textup{(}for $ k \to \infty $\textup{),} then $ \( \forall k \;\; x_k
\in c \) \impl \( x \in c \) $;

\textup{(c)} $ \lim_i \rho[i] \( x_1[i], x_2[i] \) $ exists for every
$ x_1, x_2 \in c $;

\textup{(d)} there exists a finite or countable subset $ c_1 \subset c
$ that \emph{generates} $ c $ in the sense that the least subset of $
c $ satisfying \textup{(a), (b)} and containing $ c_1 $ is the whole $
c $.
\end{definition}

\begin{remark}
Condition \ref{2c1}(d) does not change if `satisfying (a), (b)' is
replaced with `satisfying (b)'. That is, \ref{2c1}(d) is just
separability of $ c $ in the metric $ x_1, x_2 \mapsto \sup_i \rho[i]
\( x_1[i], x_2[i] \) $.
\end{remark}

The refinement of a coarse Polish space $ \(
(S[i],\rho[i])_{i=1}^\infty, c \) $ is basically the metric space $ \(
c, \ti\rho \) $, where
\[
\ti\rho ( x_1, x_2 ) = \lim_i \rho[i] \( x_1[i], x_2[i] \) \, .
\]
However, $ \ti\rho $ is a pseudometric (semimetric); it may vanish
for some $ x_1 \ne x_2 $. The equivalence class, denoted by $ x[\infty]
$, of a sequence $ x \in c $ consists of all $ x_1 \in c $ such that $
\rho[i] \( x_1[i], x[i] \) \to 0 $. On the set $ S[\infty] $ of all
equivalence classes we introduce a metric $ \rho[\infty] $,
\[
\rho[\infty] \( x_1[\infty], x_2[\infty] \) = \lim_{i\to\infty}
\rho[i] \( x_1[i], x_2[i] \) \, ;
\]
thus, $ \( S[\infty], \rho[\infty] \) $ is a metric space. We write
\[
\( S[\infty], \rho[\infty] \) = \Lim_{i\to\infty,c} \( S[i], \rho[i]
\)\index{zzLi@$ \Lim $, refinement}
\]
and call $ \( S[\infty], \rho[\infty] \) $ the
\emph{refinement}\index{refinement!of coarse Polish space} of
the coarse Polish space $ \( (S[i],\rho[i])_{i=1}^\infty,
\linebreak[0] c \) $.
Also, for every $ x = (x[1], x[2], \dots) \in c $ we denote its
equivalence class $ x[\infty] \in S[\infty] $ by
\[
x[\infty] = \Lim_{i\to\infty,c} x[i] \, ,
\]
and call it the refinement of $ x $.

\begin{lemma}\label{2c2}
For every coarse Polish space, its refinement $ (S,\rho) $ is a Polish
space.
\end{lemma}

\begin{proof}
Separability follows from \ref{2c1}(d); completeness is to be
proven. Let $
x_1, x_2, \dots $ be a Cauchy sequence in $ (S,\rho) $; we have to
find $ x \in S $ such that $ \rho (x_k,x) \to 0 $. We may assume that
$ \sum_k \rho (x_k,x_{k+1}) < \infty $. Each $ x_k $ is an equivalence
class; using (a) we choose for each $ k=1,2,3,\dots $ a representative
$ s_k \in S[1] \times S[2] \times \dots $ of $ x_k $ such that $
\sup_i \rho[i] \( s_k[i], s_{k+1}[i] \) \le 2 \rho (x_k,x_{k+1})
$. Completeness of $ \( S[i], \rho[i] \) $ ensures existence of $
s_\infty [i] = \lim_k s_k [i] $. Condition (b) ensures $ s_\infty \in
c $. The equivalence class $ x \in S $ of $ s_\infty $ satisfies $
\rho (x_k,x) \le \sup_i \rho[i] \( s_k[i], s_\infty[i] \)
\linebreak[0] \to 0 $ for $ k \to \infty $.
\qqed\end{proof}

Let $ (S[i],\rho[i])_{i=1}^\infty, c \) $ be a coarse Polish space,
and $ (S,\rho) $ its refinement. On the disjoint union $ \( S[1]
\uplus S[2] \uplus \dots \) \uplus S $ we introduce a topology,
namely, the weakest topology making continuous the following functions
$ f_s : \( S[1] \uplus S[2] \uplus \dots \) \uplus S \to [0,\infty) $
for $ s \in c $,
\begin{gather*}
f_s (x) = \rho[i] \( x, s[i] \) \quad \text{for } x \in S[i] \, , \\
f_s (x) = \rho ( x, s[\infty] ) \quad \text{for } x \in S \, ,
\end{gather*}
and an additional function $ f_0 : \( S[1] \uplus S[2] \uplus \dots \)
\uplus S \to [0,\infty) $, $ f_0 (x) = 1/i $ for $ x \in S[i] $, $ f_0
(x) = 0 $ for $ x \in S $. On every $ S[i] $ separately (and also on $
S $), the new topology coincides with the old topology, given by $
\rho[i] $ (or $ \rho $).

We may choose a sequence $ (s_k) $ dense in $ c $; the topology is
generated by functions $ f_{s_k} $ (and $ f_0 $), therefore it is a
metrizable topology. Moreover, the sequence of functions $ \( \frac{
f_{s_k} (\cdot) }{ 1 + f_{s_k} (\cdot) } \)_{k=1}^\infty $ (and $ f_0
$) maps the disjoint union into the metrizable compact space $
[0,1]^\infty $, and is a homeomorphic embedding. Thus, we have a joint
compactification of all $ S[i] $ and $ S $; and so, we treat them as
subsets of a compact metrizable space $ K $;
\[
S[i] \subset K \, , \quad S \subset K \, .
\]

\begin{lemma}
Let $ s_\infty \in S $, $ s_1 \in S[1], s_2 \in S[2], \dots $ Then $
s_i \to s_\infty $ in $ K $ if and only if $ s = (s_1,s_2,\dots) \in c
$ and $ \Lim_{i\to\infty,c} s_i = s_\infty $.
\end{lemma}

\begin{proof}
\emph{The `if' part.} The needed relation, $ f_k (s_i) \to f_k
(s_\infty) $ for $ i \to \infty $, is ensured by \ref{2c1}(c).

\emph{The `only if' part.} We choose $ x \in c $ such that $ x[\infty]
= s_\infty $; then $ \rho[i] \( s_i, x[i] \) \to \rho \( s_\infty,
x[\infty] \) = 0 $, thus $ s \in c $ by \ref{2c1}(a).
\qqed\end{proof}

The assumption `$ s_\infty \in S $' is essential. Other limiting
points (not belonging to $ S $) may exist; corresponding sequences
converge in $ K $ but do not belong to $ c $. And, of course, sets $
S, S[1], S[2], \dots $ are not closed in $ K $, unless they are
compact.

\begin{lemma}\label{2c45}
A set $ c_1 \subset c $ generates $ c $ if and only if the set of
refinements $ \{ x[\infty] : x \in c_1 \} $ is dense in $ S[\infty]
$.
\end{lemma}

\begin{proof}
\emph{The `only if' part} follows from a simple argument: if $ S' $ is
a closed subset of $ S $ then the set $ c' $ of all $ x \in c $ such
that $ x[\infty] \in S' $ satisfies \ref{2c1}(a,b).

\emph{The `if' part.} Let $ \{ x[\infty] : x \in c_1 \} $ be dense in
$ S[\infty] $ and $ s \in c $. We choose $ x_k \in c_1 $ such that $
x_k [\infty] \to s $. Similarly to the proof of \ref{2c2}, we
construct $ y_k \in c_1 $ such that $ \rho[i] \( s_k[i], y_k[i] \) \to
0 $ when $ i \to \infty $ for each $ k $, and $ \sup_i \rho[i] \(
y_k[i], s[i] \) \to 0 $ when $ k \to \infty $. The subset of $ c $
generated by $ c_1 $ contains all $ y_k $ by \ref{2c1}(a). Thus, it
contains $ s $ by \ref{2c1}(b).
\qqed\end{proof}

Given continuous functions $ f[i] : S[i] \to \R $, $ f[\infty] :
S[\infty] \to \R $, we write $ f[\infty] = \Lim_{i\to\infty,c} f[i] $
if $ f[i] (x[i]) \to f[\infty] (x[\infty]) $ whenever $ x[\infty] =
\Lim_{i\to\infty,c} x[i] $. If functions $ f[i] $ are equicontinuous
(say, $ | f[i](x) - f[i](y) | \le \rho[i] (x,y) $ for all $ i $ and $
x,y \in S[i] $), then it is enough to check that $ f[i] (x_k[i]) \to
f[\infty] (x_k[\infty]) $ for some sequence $ (x_k)_{k=1}^\infty $, $
x_k \in c $, such that the sequence $ (x_k[\infty])_{k=1}^\infty $ is
dense in $ S[\infty] $.

Given continuous maps $ f[i] : S[i] \to S[i] $, $ f[\infty] : S \to S
$, we write $ f[\infty] = \Lim_{i\to\infty,c} f[i] $ if $
\Lim_{i\to\infty,c} f[i] (x[i]) = f[\infty] (x[\infty]) $ whenever $
x[\infty] = \Lim_{i\to\infty,c} x[i] $. That is, $ \Lim \( f[i] (x[i])
 \) = \( \Lim f[i] \) \( \Lim x[i] \) $. If maps $ f[i] $ are
equicontinuous then, again, convergence may be checked on $ x_k $ such
that $ x_k[\infty] $ are dense.

\begin{sloppypar}
Given continuous maps $ f[i] : S[\infty] \to S[i] $, we may ask
whether $ \Lim_{i\to\infty,c} f[i](x) = x $ for all $ x \in S[\infty]
$, or not. If maps $ f[i] $ are equicontinuous then, still,
convergence may be checked for a dense subset of $ S[\infty] $.
\end{sloppypar}

If every $ S[i] $ is not only a metric space but also a Hilbert (or
Banach) space, and $ c $ is linear (that is, closed under linear
operations), then the refinement $ S $ is also a Hilbert (or Banach)
space, and linear operations are continuous on $ \( S[1] \cup S[2]
\cup \dots \) \cup S \subset K $ in the sense that
\[
\Lim_{i\to\infty,c} ( a s_1[i] + b s_2[i] ) = a \Lim_{i\to\infty,c}
s_1[i] + b \Lim_{i\to\infty,c} s_2[i]
\]
for all $ s_1, s_2 \in c $.

Consider the case of Hilbert spaces $ S[i] = H[i] $, $ S = H $. Given
linear\footnote{%
 Continuous, of course.}
operators $ R[i] : H[i] \to H[i] $, we may ask about $ \Lim R[i] $. If
it exists, we get
\[
\Lim \( R[i] x[i] \) = \( \Lim R[i] \) \( \Lim x[i] \) \, .
\]
If $ \sup_i \| R[i] \| < \infty $, then $ R[i] $ are equicontinuous, and
convergence may be checked on a sequence $ x_k $ such that vectors $
x_k[\infty] $ span $ H $ (that is, their linear combinations are dense
in $ H $). For example, one-dimensional orthogonal projections; if $
x[\infty] = \Lim x[i] $ then $ \Proj_{x[\infty]} = \Lim \Proj_{x[i]}
$.

Given linear operators $ R[i] : H \to H[i] $, we may ask whether $
\Lim R[i](x) = x $ for all $ x \in H $, or not. If $ \sup_i \| R[i] \|
< \infty $ then convergence may be checked on a sequence that spans $
H $. Such $ R[i] $ always exist; moreover, $ \| R[i] \| \le
1 $ may be ensured. Proof: we take $ x_k $ such that $ x_k[\infty] $
are an orthonormal basis of $ H $. After some correction, $ x_k [i] $
become orthogonal (for each $ i $), and $ \| x_k(i) \| \le 1
$.\footnote{%
 Of course, $ \| x_k[i] \| \to 1 $ for $ i \to \infty $, but in
 general we cannot ensure $ \| x_k[i] \| = 1 $. It may happen that $
 \dim H[i] < \infty $ but $ \dim H = \infty $.}
Now we let $ R[i] x_k[\infty] = x_k [i] $.

We return to coarse \emph{probability} spaces.

Let $ \( (\Om[i],\F[i],P[i])_{i=1}^\infty, \A \) $ be a coarse
probability space. For each $ i $ the pseudometric $ A,B \mapsto P[i] (
A \bigtriangleup B ) $ on $ \F[i] $ gives us the metric space $
\MALG[i] = \MALG \( \Om[i], \F[i], P[i] \) $ of all equivalence
classes of measurable sets. It is not only a metric space but also a
Boolean algebra, and moreover, a separable measure
algebra\index{measure algebra} (as defined in
\cite[17.44]{Ke}). Treating every coarse event $ A \in \A $ as a
sequence of $ A[1] \in \MALG[1], A[2] \in \MALG[2], \dots $ we get a
coarse Polish space $ \( (\MALG[i])_{i=1}^\infty, \A \) $. Its
refinement is a metric space $ \MALG[\infty] $. The set $ \A $ is
closed under Boolean operations (union, intersection,
complement). Therefore $ \MALG[\infty] $ is not only a metric space
but also a Boolean algebra. Using \ref{2c2} it is easy to check
that $ \MALG[\infty] $ is a separable measure algebra. Therefore
\cite[17.44]{Ke} it is (up to isomorphism) of the form
\[
\MALG[\infty] = \MALG (\Om,\F,P)
\]
for some probability space $ (\Om,\F,P) $. In the nonatomic case we
may take $ (\Om,\F,P) = (0,1) $ with Lebesgue measure; in general, we
may take a shorter (maybe, empty) interval plus a finite (maybe,
empty) or countable set of
atoms. Such a probability space $ (\Om.\F,P) $ (unique up to
isomorphism) will be called the
refinement\index{refinement!of coarse probability space}
of the coarse probability space $ \( (\Om[i],\F[i],P[i])_{i=1}^\infty,
\A \) $, and we write
\[
(\Om,\F,P) = \Lim_{i\to\infty,\A} \( \Om[i], \F[i], P[i] \)
\index{zzLi@$ \Lim $, refinement}
\]
(in practice, sometimes I omit ``$ i \to \infty $'' or ``$ \A $'' or
both under the ``$ \Lim $'').

Every sequence $ A = ( A[1], A[2], \dots ) \in \A $ has its refinement
\[
\Lim_{i\to\infty,\A} A[i] = A[\infty] \in \MALG (\Om,\F,P) \, .
\]

\begin{lemma}\label{2c5}
A subset $ \A_1 $ of a coarse \sif\ $ \A $ generates $ \A $ if and
only if the refinement $ \F $ of $ \A $ is generated \textup{($ \bmod
\, 0 $)} by refinements $ A[\infty] $ of all $ A \in \A_1 $.
\end{lemma}

\begin{proof}
We apply \ref{2c45} to the algebra generated by $ \A_1 $.
\qqed\end{proof}

In order to define $ L_2 (\A) $ as a set of functions on $
\Om[\all] $, we start with indicators $ \One_A $ for $ A \in \A
$, form their linear combinations, and take their completion in the
metric
\[
\| f \|_{L_2(\A)} = \sup_i \| f[i] \|_{L_2[i]} \, ,
\]
where $ L_2[i] = L_2 \( \Om[i], \F[i], P[i] \) $; the completion is a
Banach (not Hilbert) space $ L_2 (\A) $.\index{zzl2a@$ L_2 (\A) $}
Each element $ f $
of the completion is evidently identified with a sequence of $ f[i]
\in L_2[i] $, or a function on $ \Om[\all] $. We have a coarse
Polish space $ \( (L_2[i])_{i=1}^\infty, L_2(\A) \) $. It has its
refinement, $ L_2[\infty] $.

\begin{lemma}\label{2c7}
The refinement $ L_2 [\infty] $ of $ \( (L_2[i])_{i=1}^\infty, L_2(\A)
\) $ is \textup{(}canonically isomorphic to\textup{)} $ L_2 (\Om,\F,P)
$, where $
(\Om,\F,P) $ is the refinement of $ \( ( \Om[i], \F[i],\linebreak[0]
P[i] )_{i=1}^\infty, \A \) $.
\end{lemma}

\begin{proof}
We define the canonical map $ L_2 (\A) \to L_2 (\Om,\F,P) $ first on
indicators by $ \One_A \mapsto \One_{A[\infty]} $, and extend it by
linearity and continuity to the whole $ L_2 (\A) $. We note that the
image of $ f \in L_2 (\A) $ in $ L_2 (\Om,\F,P) $ depends only on the
refinement $ f[\infty] \in L_2[\infty] $ of $ f $, and their norms are
equal (both are equal to $ \lim_i \| f[i] \| $). We have a linear
isometric embedding $ L_2 [\infty] \to L_2 (\Om,\F,P) $. Its image is
closed (since $ L_2[\infty] $ is complete by \ref{2c2}), and
contains indicators $ \One_B $ for all $ B \in \MALG (\Om,\F,P) $;
therefore the image is the whole $ L_2 (\Om,\F,P) $.
\qqed\end{proof}

\begin{remark}\label{2c8}
The same holds for $ L_p $ for each $ p \in (0,\infty) $, and for the
space $ L_0 $ of all random variables (equipped with the topology of
convergence in probability). Elements of $ L_0 (\A) $ will be called
coarsely measurable\index{coarsely measurable!function}
(w.r.t.\ $ \A $) functions (on $ \Om[\all]
$), or \emph{coarse random variables;}\index{coarse!random variable}
elements of $ L_2 (\A) $ --- square integrable coarse random
variables.
\end{remark}

Let $ f $ be a coarse random variable. Then (usual) random variables $
f[i] : \Om[i] \to \R $ converge in distribution (for $ i \to \infty $)
to the refinement $ f[\infty] : \Om \to \R $. The distribution of $
f[\infty] $ will be called the 
\emph{limiting distribution}\index{limiting distribution} of $ f $.

It may happen that $ f \in L_2 (\A) $ but $ (\sgn f) \notin L_2 (\A)
$. An example: $ f(\om) = \frac{(-1)^i}{i} $ for all $ \om \in \Om[i]
$. Here, the limiting distribution is an atom at $ 0 $, and the
function `$ \sgn $' is discontinuous at $ 0 $.

\begin{lemma}
\textup{(a)} Let $ f : \Om [\all] \to \R $ be a coarse random variable,
and $ \phi : \R \to \R $ a continuous function. Then $ \phi \circ f :
\Om [\all] \to \R $ is a coarse random variable.

\textup{(b)} The same as \textup{(a)} but $ \phi $ may be
discontinuous at points of a set $ Z \subset \R $, negligible w.r.t.\
the limiting distribution of $ f $.
\end{lemma}

\begin{proof}
If $ f $ is a linear combination of indicators, then $ \phi \circ f $
is another linear combination of the same indicators. A
straightforward approximation gives (a) for uniformly continuous $
\phi $. In general, for every $ \eps $ there exists a compact set $ K
\subset \R \setminus Z $ of probability $ \ge 1-\eps $ w.r.t.\ the
limiting distribution, and also w.r.t.\ the distribution of $ f[i] $
for all $ i $ (since all these distributions are a compact set of
distributions). The restriction of $ f $ to $ K $ is uniformly
continuous. The limit for $ \eps \to 0 $ is uniform in $ i $.
\qqed\end{proof}

For a given Polish space $ S $ we may define a coarse $ S $-valued
random variable\index{coarse!random variable!$ S $-valued}
as a map $ f : \Om[\all] \to S $ such that
(usual) random variables $ f[i] : \Om[i] \to S $ converge in
distribution (for $ i \to \infty $), and $ f^{-1} (B) \in \A $ for
every $ B \subset S $ such that the boundary of $ B $ is negligible
w.r.t.\ the limiting distribution of $ f $.

For $ S = \R $ the new definition conforms with the old one.

A coarse \sif\ generated by a given sequence of sets (coarse events)
was defined after \ref{2b3}. Often it is convenient to generate a
coarse \sif\ by a sequence of functions (coarse random variables). A
function $ f : \Om[\all] \to \R $ is coarsely $ \A $-measurable
if and only if $ \A $ contains sets $ f^{-1} \( (-\infty,x) \) $ for
all $ x \in \R $ except for atoms (if any) of the limiting
distribution of $ f $. A dense countable subset of these $ x $ is
enough. So, a coarse \sif\
generated\index{generated!coarse \sif, by functions}
by a finite or countable set of
functions $ f $ is nothing but the coarse \sif\ generated by a
countable set of sets of the form $ f^{-1} \( (-\infty,x) \) $. More
generally, $ S $-valued (coarse) random variables may be used; they
are reduced to the real-valued case by composing with appropriate
continuous functions $ S \to \R $.

\begin{lemma}\label{2c9}
\begin{sloppypar}
A sequence of functions $ f_k : \Om[\all] \to \R $ generates a
coarse \sif\ if and only if for every $ n $, $ n $-dimensional random
variables $ \( f_1[i], \dots, f_n[i] \) : \Om[i] \to \R^n $ converge
in distribution \textup{(}for $ i \to \infty $\textup{).}
\end{sloppypar}
\end{lemma}

\begin{proof}
\emph{The `only if' part.} Let $ f_1, \dots, f_n $ be coarsely
measurable (w.r.t.\ some coarse \sif), then they have a limiting joint
distribution.

\emph{The `if' part.} For each $ n $ we choose a dense countable set $
Q_n \subset \R $ negligible w.r.t.\ the limiting distribution of $ f_n
$. We apply \ref{2b3} to the set $ \A_1 $ of coarse events of
the form $ \{ f_1(\cdot) \le q_1, \dots, f_n(\cdot) \le q_n \} $ where
$ q_1 \in Q_1, \dots, q_n \in Q_n $, $ n = 1,2,\dots $
\qqed\end{proof}

\begin{remark}
The same holds for an arbitrary Polish space instead of $ \R $.
\end{remark}

\begin{remark}
Comparing \ref{2c9} and \eqref{condf} we see that every joint
compactification of $ \Om_1 \uplus \Om_2 \uplus \dots $ (in the sense
of \ref{sec:1.3}, assuming \eqref{condC}) may be downgraded to a
coarse probability
space. Namely, we take a sequence of functions $ f_k $ that generates
$ C $ and consider the coarse \sif\ $ \A $ generated by $ (f_k)
$. Every $ f \in C $ is a coarse random variable, since $ L_0 (\A) $
is closed under all operations used in \eqref{1cg}, \eqref{1ch}, or
\eqref{1ci}.\footnote{%
 Of course, $ L_0 (\A) $ usually contains no sequence dense in the
 \emph{uniform} topology.}
Therefore $ \A $ does not depend on the choice of $ (f_k) $.
\end{remark}

\section{Scaling Limit and Independence}
\label{sec:3}
\subsection{Product of coarse probability spaces}
\label{sec:3.1}

Having two coarse probability spaces $ \(
(\Om_1[i],\F_1[i],P_1[i])_{i=1}^\infty, \A_1 \) $ and $ \(
(\Om_2[i],\linebreak[0]
\F_2[i],P_2[i])_{i=1}^\infty, \linebreak[0]
\A_2 \) $, we define their 
product\index{product!of coarse probability spaces}
as the coarse probability space $ \(
(\Om[i],\F[i],P[i])_{i=1}^\infty, \A \) $ where for each $ i $,
\[
(\Om[i],\F[i],P[i]) = (\Om_1[i],\F_1[i],P_1[i]) \times
(\Om_2[i],\F_2[i],P_2[i])
\]
is the usual product of probability spaces, and $ \A $ is the smallest
coarse \sif\ that contains $ \{ A_1 \times A_2 : A_1 \in \A_1, A_2 \in
\A_2 \} $, where $ A_1 \times A_2 \subset \Om[\all] $ is defined
by $ \forall i \;\; ( A_1 \times A_2 ) [i] = A_1 [i] \times A_2 [i]
$. Existence of such $ \A $ is ensured by \ref{2b3}. We write $
\A = \A_1 \otimes \A_2 $.\index{zzz@$ \otimes $!for coarse \sif s}

\begin{lemma}\label{3a1}
The refinement of the product of two coarse probability spaces is
(canonically isomorphic to) the product of their refinements.
\end{lemma}

\begin{proof}
\begin{sloppypar}
Denote these refinements by $ (\Om_1,\F_1,P_1) $, $ (\Om_2,\F_2,P_2)
$ and $ (\Om,\F,P) $. Both $ \MALG (\Om_1,\F_1,P_1) $ and $ \MALG
(\Om_2,\F_2,P_2) $ are naturally embedded into $ \MALG (\Om,\F,P) $ as
\emph{independent} subalgebras. They generate $ \MALG (\Om,\F,P) $ due
to \ref{2c5}.
\qqed\end{sloppypar}
\end{proof}

Given an arbitrary coarse \sif\ $ \A $ on the product coarse sample
space $ \( (\Om_1[i],\F_1[i],P_1[i]) \linebreak[0]
\times (\Om_2[i],\F_2[i],P_2[i])
\)_{i=1}^\infty $, we may ask whether $ \A $ is a product, that
is, $ \A = \A_1 \otimes \A_2 $ for some $ \A_1, \A_2 $, or not. No
need to check all $ \A_1, \A_2 $. Rather, we have to check
\[
\A_1 = \{ A_1 : A_1 \times \Om_2 \in \A \} \, , \quad
\A_2 = \{ A_2 : \Om_1 \times A_2 \in \A \} \, ; 
\]
of course, $ A_1 \times \Om_2 \subset \Om[\all] $ is defined by
$ \forall i \;\, (A_1\times\Om_2) [i] = A_1 [i] \times \Om_2 [i] $. If
$ \{ A_1 \times A_2 : A_1 \in \A_1, A_2 \in \A_2 \} $ generates $ \A
$, then $ \A $ is a product; otherwise, it is not.

The refinement $ \F $ of $ \A $ contains two sub-\sif s $ \F_1 = \{
(A_1\times\Om_2) [\infty] : A_1 \in \A_1 \} $, $ \F_2 = \{
(\Om_1\times A_2) [\infty] : A_2 \in \A_2 \} $. They are independent:
\[
P ( A \cap B ) = P(A) P(B) \quad \text{for } A \in \F_1, \, B \in \F_2
\, .
\]

\begin{lemma}
$ \A $ is a product if and only if $ \F_1, \F_2 $ generate $ \F $.
\end{lemma}

\begin{proof}
We apply \ref{2c5} to $ \{ A_1 \times A_2 : A_1 \in \A_1, A_2
\in \A_2 \} $.
\qqed\end{proof}

\begin{remark}
It is well-known that a generating pair of independent sub-\sif s
means that $ (\Om,\F,P) $ is (isomorphic to) the product of two
probability spaces. So, a coarse probability space is a product if and
only if its refinement is a product. (Assuming, of course, that the
coarse \emph{sample} space is a product.)
\end{remark}

Let $ \A = \A_1 \otimes \A_2 $.
Consider Hilbert spaces $ H_1[i] = L_2 (\Om_1[i], \F_1[i], P_1[i]) $, \,
$ H_2[i] = L_2 (\Om_2[i], \F_2[i], P_2[i]) $, \, $ H[i] = L_2 (\Om[i],
\F[i], P[i]) $. For each $ i $, the space $ H[i] $ is (canonically
isomorphic to)
$ H_1[i] \otimes H_2[i] $.\index{zzz@$ \otimes $!for Hilbert spaces}
Indeed, for $ x_1 \in H_1[i]
$, $ x_2 \in H_2[i] $ we define $ x_1 \otimes x_2 \in H[i] $ by $ (x_1
\otimes x_2) (\om_1,\om_2) = x_1(\om_1) x_2(\om_2) $; then $ \ip{ x_1
\otimes x_2 }{ y_1 \otimes y_2 } = \ip{ x_1 }{ y_1 } \ip{ x_2 }{ y_2 }
$, and factorizable vectors (of the form $ x_1 \otimes x_2 $) span the
space $ H[i] $. We know (see \ref{2c7}) that the refinement $
H[\infty] $ of $ \( (H[i])_{i=1}^\infty, L_2(\A) \) $ is $ L_2
(\Om,\F,P) $. Also, $ H_1[\infty] = L_2 (\Om_1,\F_1,P_1) $ and $
H_2[\infty] = L_2 (\Om_2,\F_2,P_2) $. Using \ref{3a1} we get $
H[\infty] = H_1[\infty] \otimes H_2[\infty] $. In that sense,
\[
\Lim \( H_1[i] \otimes H_2[i] \) = \( \Lim H_1[i] \) \otimes \( \Lim
H_2[i] \) \, .
\]
If $ x \in L_2 (\A_1) $, $ y \in L_2 (\A_2) $, we
define $ x \otimes y $ by $ (x \otimes y) [i] = x[i] \otimes y[i] $
for all $ i $. We get $ x \otimes y \in L_2 (\A) $ and $ (x
\otimes y) [\infty] = x[\infty] \otimes y[\infty] $, that is,
\begin{equation}\label{3a2}
\Lim \( x[i] \otimes y[i] \) = \( \Lim x[i] \) \otimes \( \Lim y[i] \)
\, ,
\end{equation}
since it holds for (linear combinations of) indicators of coarse
events. Note also that linear combinations of factorizable vectors are
dense in $ L_2 (\A) $.

Assume that $ R_1 [i] : H_1 [i] \to H_1 [i] $, $ R_2 [i] : H_2 [i] \to
H_2 [i] $ are linear operators, possessing limits $ R_1 [\infty] =
\Lim R_1[i] $, $ R_2 [\infty] = \Lim R_2[i] $. Consider linear
operators $ R_1[i]  \otimes R_2[i] = R[i] : H[i] \to H[i] $. (It means
that $ R[i] x[i] = R_1[i] x_1[i] \otimes R_2[i] x_2[i] $ whenever $
x[i] = x_1[i] \otimes x_2[i] $.) If $
\sup_i \| R_1[i] \| < \infty $, $ \sup_i \| R_2[i] \| < \infty $, then
$ \Lim R[i] = R_1 [\infty] \otimes R_2 [\infty] $, that is,
\begin{equation}
\Lim \( R_1[i] \otimes R_2[i] \) = \( \Lim R_1[i] \) \otimes \( \Lim
R_2[i] \) \, .
\end{equation}
Proof: We have to check that
\[
\Lim \( R_1[i] \otimes R_2[i] \) x[i] = \( \Lim R_1[i] \otimes \Lim
R_2[i] \) \( \Lim x[i] \)
\]
for all $ x \in L_2(\A) $. We may assume that $ x $ is factorizable, $
x = x_1 \otimes x_2 $; then
\begin{multline*}
\Lim \( R_1[i] \otimes R_2[i] \) \( x_1[i] \otimes x_2[i] \) = \\
= \Lim \( R_1[i] x_1[i] \otimes R_2[i] x_2[i] \) = \\
= \( \Lim R_1[i] x_1[i] \) \otimes \( \Lim R_2[i] x_2[i] \) = \\
= \( \Lim R_1[i] \) \( \Lim x_1[i] \) \otimes \( \Lim R_2[i] \) \(
 \Lim x_2[i] \) = \\
= \( \Lim R_1[i] \otimes \Lim R_2[i] \) \( \Lim x_1[i] \otimes \Lim
 x_2[i] \) \, .
\end{multline*}

Especially, let $ R_2[i] $ be the orthogonal projection to the
one-dimensional subspace of constants (basically, the expectation),
and $ R_1[i] $ be the unit (identity) operator. Then $ \( R_1[i]
\otimes R_2[i] \) \( x[i] \) = \cE{ x[i] }{ \F_1[i] } $, since it
holds for factorizable vectors. Further, $ R_2[\infty] = \Lim R_2[i] $
is the expectation on $ (\Om_2,\F_2,P_2) $, since convergence of
vectors implies convergence of one-dimensional projections, and
constant functions on $ \Om_2[\all] $ belong to $ L_2 (\A)
$. So,
\begin{equation}\label{3a6}
\Lim \cE{ x[i] }{ \F_1[i] } = \cE{ \Lim x[i] }{ \F_1 }
\end{equation}
for all $ x \in L_2 (\A) $.

All the same holds for the product of any finite number of spaces (not
just two).

\subsection{Dyadic case}
\label{sec:3.2}

Let $ (\Om[i],\F[i],P[i]) $ be the space of all maps $ \frac1i \Z \to
\{-1,+1\} $ with the usual product measure. That is, we have
independent random signs $ \tau_{k/i} $ for all integers $ k
$;\footnote{%
 Rigorously, I should denote it by $ \tau_k[i] $, but $ \tau_{k/i} $
 is more expressive. Though $ \tau_{2/6} $ is not the same as $
 \tau_{1/3} $, hopefully, it does not harm.}
each random sign takes on two values $ \pm1 $ with probabilities $
50\%,50\% $. The coarse sample space $
(\Om[i],\F[i],P[i])_{i=1}^\infty $ will be called the \emph{dyadic
coarse sample space.}\index{dyadic coarse sample space}\footnote{%
 Sometimes a subsequence is used; say, $ i \in \{ 2,4,8,16,\dots \} $
 only; or equivalently, $ \Om[i] $ is the space of maps $ 2^{-i} \Z \to
 \{-1,+1\} $; see \ref{3b6}, \ref{3b7}.}
Let $ \A $ be a coarse \sif\ on the dyadic
coarse sample space. What about decomposing it, say, into the past and
the future w.r.t.\ a given instant?

Let us define a \emph{coarse instant}\index{coarse!instant}
as a sequence $ t =
\(t[i])_{i=1}^\infty $ such that $ t[i] \in \frac1i\Z $ (that is, $ i
t[i] \in \Z $) for all $ i $, and there exists $ t[\infty] \in \R $
(call it the refinement\index{refinement!of coarse instant}
of the coarse instant) such that $ t[i] \to
t[\infty] $ for $ i \to \infty $. A \emph{coarse time
interval}\index{coarse!time interval}
is a pair $ (s,t) $ of coarse instants $ s,t $ such that $ s \le t $
in the sense that $ s[i] \le t[i] $ for all $ i $.

For every coarse time interval $ (s,t) $ we define the coarse
probability space $ \( ( \Om_{s,t}[i], \F_{s,t}[i], \linebreak[0]
P_{s,t}[i] )_{i=1}^\infty, \A_{s,t} \) $ as follows. First, $
\Om_{s,t}[i] $\index{zzOmst@$ \Om_{s,t}[i] $}
is the space of all maps $ \( \frac1i\Z \cap [ s[i],t[i] ) \) \to \{
-1, +1 \} $.\footnote{%
 It may happen that $ s[i] = t[i] $, then $ \Om_{s,t} [i] $ contains a
 single point.}
Second, $ \F_{s,t} [i] $\index{zzFst@$ \F_{s,t}[i] $}
and $ P_{s,t} [i] $\index{zzPst@$ P_{s,t}[i] $}
are defined
naturally, and we have the canonical measure preserving map $
(\Om[i],\F[i],P[i]) \to ( \Om_{s,t}[i], \F_{s,t}[i], P_{s,t}[i] )
$. Third, each $ A \subset \Om_{s,t}[\all] $ has its inverse
image in $ \Om[\all] $; if the inverse image of $ A $ belongs to
$ \A $ then (and only then) $ A $ belongs to $ \A_{s,t} $, which is
the definition of $ \A_{s,t} $. It is easy to see that $ \A_{s,t} $ is
a coarse \sif.

Given coarse time intervals $ (r,s) $ and $ (s,t) $, we have
\[
\( \Om_{r,t}[i], \F_{r,t}[i], P_{r,t}[i] \) = \( \Om_{r,s}[i],
\F_{r,s}[i], P_{r,s}[i] \) \times \( \Om_{s,t}[i], \F_{s,t}[i],
P_{s,t}[i] \) \, ,
\]
and we may ask whether $ \A_{r,t} $ is a product, that is, $ \A_{r,t}
= \A_{r,s} \otimes \A_{s,t} $, or not.

\begin{definition}\label{3b1}
\begin{sloppypar}
A \emph{dyadic coarse
factorization}\index{dyadic coarse factorization}
is a coarse probability space $
\( (\Om[i],\F[i],P[i])_{i=1}^\infty, \A \) $ such that $
(\Om[i],\F[i],P[i])_{i=1}^\infty $ is the dyadic coarse sample space;
\[
\A_{r,t} = \A_{r,s} \otimes \A_{s,t}
\]
whenever $ r,s,t $ are coarse instants such that $ r[i] \le s[i] \le
t[i] $ for all $ i $; and
\[
\A \text{ is generated by } \bigcup_{(s,t)} \A_{s,t} \, ,
\]
where the union is taken over all coarse time intervals $ (s,t) $.
\end{sloppypar}
\end{definition}

\begin{example}
A single function $ f : \Om[\all] \to \R $, defined by $ f(\om) =
\tau_{0/i}(\om) $ for $ \om \in \Om[i] $, generates a coarse \sif\ $
\A $. However, the coarse probability space $ \(
(\Om[i],\F[i],P[i])_{i=1}^\infty, \A \) $ is \emph{not} a dyadic
coarse factorization. The equality $ \A_{r,t} = \A_{r,s} \otimes
\A_{s,t} $ is violated when $ s[i] $ converges to $ 0 $ from both
sides; say, $ s[i] = (-1)^i/i $. It means that a single point of the
time continuum should not carry a random sign. See also
\ref{3b8}--\ref{3b10}.
\end{example}

Every family $ (\A_{s,t})_{s\le t} $ of coarse \sif s $ \A_{s,t} $ on
coarse sample spaces $ \( \Om_{s,t}[i], \F_{s,t}[i], P_{s,t}[i]
\)_{i=1}^\infty $, indexed by all coarse time intervals $ (s,t) $ and
satisfying $ \A_{r,t} = \A_{r,s} \otimes \A_{s,t} $ whenever $ r \le s
\le t $, corresponds to a dyadic coarse factorization.

\begin{example}\label{3b3}
Given a coarse time interval $ (s,t) $, we consider $ f_{s,t} :
\Om[\all] \to \R $,
\[
f_{s,t} (\om) =
\frac1{\sqrt i} \sum_{k:s[i]\le k/i<t[i]} \tau_{k/i}(\om)
\quad \text{for } \om \in \Om[i] \, .
\]
Only $ s[\infty], t[\infty] $ matter, in the sense that
\begin{equation}\label{3b3a}
\int_{\Om[i]} \frac{ | \ti f[i] - f[i] | }{ 1 + | \ti f[i] - f[i] | }
\, dP[i] \tendsi 0
\end{equation}
if $ f = f_{s,t} $, and
$ \ti f = f_{\ti s,\ti t} $ is such a function built for a different
coarse time interval $ (\ti s, \ti t) $ satisfying $ \ti s[\infty] =
s[\infty] $, $ \ti t[\infty] = t[\infty] $.
Moreover, $ \| \ti f[i] - f[i] \|_{L_2[i]} \to 0 $ for $ i \to \infty
$. We choose a sequence of
coarse time intervals, $ (s_n,t_n)_{n=1}^\infty $, such that the
sequence of their refinements, $ (s_n[\infty],t_n[\infty]) $ is dense
among all (usual, not coarse) intervals.
The sequence $ \(f_{s_n,t_n}\)_{n=1}^\infty $ satisfies the
condition of \ref{2c9} and therefore it generates a coarse \sif\
$ \A $. It is easy to see that $ \A $ does not depend on the choice of
$ (s_n,t_n) $. Clearly, the refinement of $ f_{s,t} $ is the increment $
B(t[\infty]) - B(s[\infty]) $ of the usual Brownian motion $ B(\cdot)
$.

Given three coarse instants $ r \le s \le t $, we have
\[
f_{r,t} = f_{r,s} + f_{s,t} \, .
\]
It shows that $ f_{r,t} $ is
coarsely measurable w.r.t.\ the product of two coarse \sif s $
\A_{r,s} \otimes \A_{s,t} $, which implies $ \A_{r,t} = \A_{r,s} \otimes
\A_{s,t} $. So, we have a dyadic coarse factorization. We may call it the
Brownian coarse factorization.
\end{example}

\begin{example}\label{3b4}
Let $ f_{s,t}(\om) $ be the same as in \ref{3b3} and in addition,
\[
g_{s,t} (\om) = \frac1{\sqrt i} \sum_{k:s[i]\le k/i<t[i]} (-1)^k
\tau_{k/i}(\om) \quad \text{for } \om \in \Om[i] \, .
\]
In the scaling limit we get two independent Brownian motions $ B_1,
B_2 $; the refinement of $ f_{s,t} $ is $ B_1(t[\infty]) -
B_1(s[\infty]) $, the refinement of $ g_{s,t} $ is $ B_2(t[\infty]) -
B_2(s[\infty]) $. By the way, $ (-1)^k $ cannot be replaced with $
(-1)^{k-s[i]} $; it would violate the condition of \ref{2c9}.

We may also consider
\[
f^{(n)}_{s,t} (\om) = \frac1{\sqrt i} \sum_{k:s[i]\le k/i<t[i]} \exp
\bigg( 2\pi \I \frac k n \bigg) \tau_{k/i}(\om) \quad
\text{for } \om \in \Om[i]
\]
for $ n=1,2,3,\dots $ (here $ \I = \sqrt{-1} $, while $ i $ is an
integer). In the scaling limit we get two real-valued
Brownian motions $ B_1, B_2 $ and infinitely many complex-valued
Brownian motion $ B_3, B_4, \dots $ All $ B_n $ are independent.

Another construction of that kind:
\[
f^{(\la)}_{s,t} (\om) = \frac1{\sqrt i} \sum_{k:s[i]\le k/i<t[i]} \exp
\bigg( 2\pi \I \la \frac k {\sqrt i} \bigg)
\tau_{k/i}(\om) \quad \text{for } \om \in \Om[i] \, .
\]
In the scaling limit, each $ \la \in (0,\infty) $ gives a
complex-valued Brownian motion $ B_\la $. Any finite or countable set
of numbers $ \la $ may be used, and leads to independent Brownian
motions. Note that we cannot use more than a countable set of $ \la $,
since separability is stipulated by the definition of a coarse
probability space.
\end{example}

\begin{example}\label{3b5}
For $ n=1,2,\dots $ we introduce
\[
f^{(n)}_{s,t} (\om) = \frac1{\sqrt i} \sum_{k:s[i]\le k/i\le(k+n)/i<t[i]}
\prod_{m=1}^n \tau_{(k+m)/i} (\om) \quad \text{for } \om \in \Om[i] \,
.
\]
In the scaling limit we get independent Brownian motions $ B_n $.

Another construction of that kind:
\[
f^{(\la)}_{s,t} (\om) = \frac1{\sqrt i} \sum_{k:s[i]\le
k/i\le(k+\la\sqrt i)/i<t[i]}
\prod_{m=1}^{\entier(\la\sqrt i)} \tau_{(k+m)/i} (\om) \quad \text{for
} \om \in \Om[i] \, ;
\]
any finite or countable set of numbers $ \la \in (0,\infty) $ may be
used, and leads to independent Brownian motions $ B_\la $.

Note that we cannot take the product over $ m=1,\dots,\entier(\la i)
$; that would destroy factorizability.
\end{example}

\begin{example}\label{3b6}
Here we restrict ourselves to $ i \in \{2,4,8,16,\dots\} $, thus
violating a little of our framework. We let for $ \om \in \Om[i] $, $
i = 2^n $,
\[
g_{s,t} (\om) = \sum_{k:s[i]\le k/i<(k+n-1)/i<t[i]}
\frac{1+\tau_{k/i}(\om)}2 \prod_{m=1}^{n-1}
\frac{1-\tau_{(k+m)/i}(\om)}2 \, .
\]
That is, $ g_{s,t} : \Om[\all] \to \{ 0,1,2,\dots \} $ counts
combinations `$ +-\dotsc- $' of one plus sign and $ (n-1) $ minus signs
in succession. In the scaling limit we get the Poisson process.
\end{example}

\begin{example}\label{3b7}
Let $ f_{s,t} $ be as in \ref{3b3} (Brownian), while $ g_{s,t}
$ is as in \ref{3b6} (Poisson). Taken together, they generate
a coarse \sif. The corresponding scaling limit consists of two
\emph{independent} processes, Brownian and Poisson.
\end{example}

Let $ \( (\Om[i],\F[i],P[i])_{i=1}^\infty, \A \) $ be a dyadic coarse
factorization. Being a coarse probability space, it has a refinement $
(\Om,\F,P) $. For every coarse time interval $ (s,t) $ we have a
coarse sub-\sif\ $ \A_{s,t} \subset \A $ and its refinement, a
sub-\sif\ 
$ \F_{s,t}[\infty] \subset \F $.\index{zzFst@$ \F_{s,t}[\infty] $}
By \ref{3a1},
\[
\F_{r,t}[\infty] = \F_{r,s}[\infty] \otimes \F_{s,t}[\infty] \quad
\text{whenever } r \le s \le t \, .
\]

\begin{lemma}\label{3b8}
If $ s[\infty] = t[\infty] $ then $ \F_{s,t}[\infty] $ is degenerate
(that is, contains sets of probability $ 0 $ or $ 1 $ only).
\end{lemma}

\begin{proof}
Consider the coarse instant $ r $,
\[
r[i] = \begin{cases}
 s[i] &\text{for $ i $ even},\\
 t[i] &\text{for $ i $ odd}.
\end{cases}
\]
For every $ A \in \A_{s,r} $,
\[
P (A[\infty]) = \lim_{i\to\infty} P[i] \( A[i] \) = \lim_{i\to\infty} P[2i] \(
A[2i] \) \in \{ 0,1 \} \, ,
\]
since $ \A_{s,r} [2i] $ is degenerate. So, $ \F_{s,r}[\infty] $ is
degenerate. Similarly, $ \F_{r,t}[\infty] $ is degenerate. However, $
\F_{s,t}[\infty] = \F_{s,r}[\infty] \otimes \F_{r,t}[\infty] $.
\qqed\end{proof}

\begin{lemma}
$ \F_{s,t}[\infty] $ depends only on $ s[\infty], t[\infty] $.
\end{lemma}

\begin{proof}
Let $ (u,v) $ be another coarse time interval such that $ u[\infty] =
s[\infty] $ and $ v[\infty] = t[\infty] $; we have to prove that $
\F_{s,t}[\infty] = \F_{u,v}[\infty] $. Assume that $ s[\infty] <
t[\infty] $ (otherwise both $ \F_{s,t}[\infty] $ and $
\F_{u,v}[\infty] $ are degenerate). Assume also that $ s[i] \le v[i] $
and $ u[i] \le t[i] $ for all $ i $ (otherwise we correct them on a
finite set of indices $ i $). 

Further, we may assume that $ s \le u \le v \le t $; otherwise we
turn to $ s \wedge u \le s \vee u \le t \wedge v \le t \vee v $, where
$ ( s \wedge u ) [i] = s[i] \wedge u[i] = \min \( s[i], u[i] \) $, etc.
Both $ \F_{s,t}[\infty] $ and $ \F_{u,v}[\infty] $ are sandwiched
between $ \F_{s\wedge u,t\vee v}[\infty] $ and $ \F_{s\vee u,t\wedge
v}[\infty] $.

Finally, $ \F_{s,t}[\infty] = \F_{s,u}[\infty] \otimes
\F_{u,v}[\infty] \otimes \F_{v,t}[\infty] = \F_{u,v}[\infty] $, since
$ \F_{s,u}[\infty] $ and $ \F_{v,t}[\infty] $ are degenerate by \ref{3b8}.
\qqed\end{proof}

So, a sub-\sif\ $ \F_{s,t} \subset \F $ is well-defined for every
interval $ (s,t) \subset \R $ (rather than a coarse time interval),
and
\[
\F_{r,t} = \F_{r,s} \otimes \F_{s,t} \quad \text{whenever } -\infty <
r \le s \le t < +\infty \, .
\]

\begin{lemma}\label{3b10}
The union of sub-\sif s $ \F_{s+\eps,t-\eps} $ over $ \eps > 0 $
generates $ \F_{s,t} $.
\end{lemma}

\begin{proof}
Consider $ \F_{\eps,1} $. We have to prove that $ \cE{ x }{
\F_{\eps,1} } $ converges to $ x $ (in $ L_2 (\Om) $, for $ \eps\to0+
$) for every $ x \in L_2 (\F_{0,1}) $, or for $ x[\infty] $ where $ x
\in L_2 (\A_{0,1}) $. Assume the contrary. Then
\[
\| \cE{ x[\infty] }{ \F_{\eps,1} } \| < c < \| x[\infty] \|
\]
for all $ \eps $ small enough, and some constant $ c $. We know that
\[
\cE{ x[\infty] }{ \F_{\eps,1} } = \Lim \cE{ x[i] }{ \F_{\eps,1} [i] }
\]
for each $ \eps $.\footnote{%
 Or rather, an appropriate coarse instant is meant in $ \F_{\eps,1}[i]
 $.}
Therefore
\[
\| \cE{ x[i] }{ \F_{\eps,1}[i] } \| \tendsi \| \cE{ x[\infty] }{
 \F_{\eps,1} } \| < c \, .
\]
We choose a sequence $ \eps[i] \tendsi 0 $ such that $ \| \cE{ x[i] }{
\F_{\eps[i],1} [i] } \| < c $ for all $ i $ large enough. However, $
\Lim \cE{ x[i] }{ \F_{\eps[i],1} [i] } = \cE{ x[\infty] }{
\F_{\eps[\infty],1} } = \cE{ x[\infty] }{ \F_{0,1} } = x[\infty] $; a
contradiction.
\qqed\end{proof}

\subsection{Scaling limit of Fourier-Walsh coefficients}
\label{sec:3.3}

We still consider a dyadic coarse factorization. The Hilbert space $
L_2[i] = L_2 \( \Om[i], \F[i], P[i] \) $ consists of all functions of
random signs $ \tau_m $, $ m \in \frac1i\Z $. The well-known
Fourier-Walsh orthonormal basis of $ L_2[i] $ consists of products
\[
\tau_M = \prod_{m\in M} \tau_m \, , \quad M \in \cC[i] \, , \quad
\cC[i] = \{ M \subset \tfrac1i\Z : M \text{ is finite} \} \, .
\index{zzCi@$ \cC[i] $, the set of finite sets}
\]
Every $ f \in L_2[i] $ is of the form
\[
f = \sum_M \hat f_M \tau_M = \hat f_\emptyset + \sum_{m\in\frac1i\Z}
\hat f_{\{m\}} \tau_m + \sum_{m_1,m_2\in\frac1i\Z,m_1<m_2} \hat
f_{\{m_1,m_2\}} \tau_{m_1} \tau_{m_2} + \dots \, ;
\]
coefficients $ \hat f_M $ are called Fourier-Walsh coefficients of $ f
$. We define the
\emph{spectral measure}\index{spectral measure (discrete case)}
$ \mu_f $ on the countable set $ \cC[i] $ by
\[
\mu_f (\cM) = \sum_{M\in\cM} | \hat f_M |^2 \quad \text{for } \cM
\subset \cC[i] \, ;
\index{zzmuf@$ \mu_f $, spectral measure (discrete)}
\]
it is a finite positive measure,
\[
\mu_f ( \cC[i] ) = \| f \|^2 \, ; \quad \mu_f ( \{ \emptyset \} ) =
(\Ex f)^2 \, ; \quad \mu_f ( \cC[i] \setminus \{\emptyset\} ) =
\Var(f) \, .
\]

Let $ (s,t) $ be a coarse time interval. We have
\begin{gather*}
\cE{ \tau_M }{ \F_{s,t}[i] } = \begin{cases}
 \tau_M &\text{if $ M \subset [ s[i], t[i] ) $},\\
 0 &\text{otherwise;}
\end{cases} \\
\| \cE{ f }{ \F_{s,t}[i] } \|^2 = \mu_f \( \{ M \in \cC[i] : M \subset
 [ s[i], t[i] ) \} \) \, .
\end{gather*}
We apply it to $ f = x[i] $ for an arbitrary $ x \in L_2 (\A) $ and
arbitrary $ i $; $ \mu_f $ becomes $ \mu_{x[i]} $ or $ \mu_x[i] $; by
\eqref{3a6},
\begin{multline*}
\mu_x[i] \( \{ M \in \cC[i] : M \subset [ s[i], t[i] ) \} \) = \| \cE{
 x[i] }{ \F_{s,t}[i] } \|^2 \\
\tendsi \| \cE{ x[\infty] }{ \F_{s,t}[\infty] } \|^2 \, .
\end{multline*}
For every $ \eps > 0 $ we can choose $ s,t $ so that $
\| x[\infty] \|^2 - \| \cE{ x[\infty] }{ \F_{s,t}[\infty] } \|^2 \le
\eps $, and moreover,
\begin{equation}\label{3c1}
\mu_x[i] \( \{ M \in \cC[i] : M \subset [ s[i], t[i] ) \} \) \le \eps
\quad \text{for all } i \, .
\end{equation}
We consider each $ \mu_x[i] $ as a measure on the space $ \cC[\infty]
$\index{zzCinf@$ \cC[\infty] $, space of compact sets}
of all compact subsets of $ \R $, equipped with the Hausdorff
metric; the metric is
\begin{equation}\label{3c2}
\dist ( M_1, M_2 ) = \sup_{x\in\R} \Big| \min_{y\in M_1} |x-y| -
\min_{y\in M_2} |x-y| \Big|
\end{equation}
for nonempty $ M_1, M_2 $, and $ \dist ( \emptyset, M ) = 1 $ for $ M
\ne \emptyset $. Clearly, $ \cC[i] \subset \cC[\infty] $ for each $ i
$; thus, a measure on $ \cC[i] $ is also a measure on $ \cC[\infty]
$.\footnote{%
 One may turn $ (\cC[i])_{i=1}^\infty $ into a coarse Polish space,
 and identify its refinement with $ \cC[\infty] $. It leads to a joint
 compactification of all $ \cC[i] $ and $ \cC[\infty] $, which is a
 suitable framework for weak convergence of measures on $ \cC[i] $ to a
 measure on $ \cC[\infty] $. However, it is simpler to use natural
 embeddings, $ \cC[i] \subset \cC[\infty] $.}
The set $ \{ M \in \cC[\infty] : M \subset [u,v] \} $ is well-known
to be compact, for every $ [u,v] \subset \R $. Thus, \eqref{3c1} shows
that the sequence of measures $ \mu_x[i] $ on $ \cC[\infty] $ is
tight.

Let $ (s_1,t_1) $ and $ (s_2,t_2) $ be two coarse time intervals, $
s_1 \le t_1 \le s_2 \le t_2 $. Sub-\sif s $ \F_{s_1,t_1}[i] $ and $
\F_{s_2,t_2}[i] $ are independent; they generate a sub-\sif\ that may
be denoted by
\[
\F_{(s_1,t_1)\cup(s_2,t_2)} [i] = \F_{s_1,t_1} [i] \otimes
\F_{s_2,t_2} [i] \, .
\]
We have
\begin{gather*}
\cE{ \tau_M }{ \F_{(s_1,t_1)\cup(s_2,t_2)}[i] } = \begin{cases}
 \tau_M &\text{if $ M \subset [s_1[i],t_1[i]) \cup [s_2[i],t_2[i]) $},\\
 0 &\text{otherwise;}
\end{cases} \\
\| \cE{ f }{ \F_{(s_1,t_1)\cup(s_2,t_2)}[i] } \|^2 = \mu_f \( \{ M \in
 \cC[i] : M \subset [s_1[i],t_1[i]) \cup [s_2[i],t_2[i]) \} \) \, ; \\
\mu_x[i] \( \{ M \in \cC[i] : M \subset [s_1[i],t_1[i]) \cup
 [s_2[i],t_2[i]) \} \) = \\
= \| \cE{ x[i] }{
\F_{(s_1,t_1)\cup(s_2,t_2)}[i] } \|^2 \tendsi \| \cE{ x[\infty] }{
\F_{(s_1,t_1)\cup(s_2,t_2)} [\infty] } \|^2 \, ,
\end{gather*}
where $ \F_{(s_1,t_1)\cup(s_2,t_2)}[\infty] = \F_{s_1,t_1} [\infty]
\otimes \F_{s_2,t_2} [\infty] = \F_{s_1[\infty],t_1[\infty]} \otimes
\F_{s_2[\infty],t_2[\infty]} $. A generalization of \eqref{3a6} to the
product of more than two spaces was used here.

The same holds for more than two coarse time intervals:
\begin{multline}\label{3c3}
\mu_x[i] \( \{ M \in \cC[i] : M \subset [s_1[i],t_1[i]) \cup \dotsc
 \cup [s_n[i],t_n[i]) \} \) \\
\tendsi \| \cE{ x[\infty] }{
 \F_{(s_1,t_1)\cup\dotsc\cup(s_n,t_n)} [\infty] } \|^2 \, .
\end{multline}
We have convergence of spectral measures on a special class of subsets
of $ \cC[\infty] $. Note that the intersection of two such subsets is
again such a subset. Therefore, the convergence holds on the algebra
of subsets generated by the class. A generic element of the algebra
is the union of a finite number of `cells' of the form
\begin{equation}\label{3c4}
\{ M \in \cC[\infty]: M \subset \cup_{k=1}^n [ s_k, t_k ) \text{
and } M \cap [ s_k, t_k ) \ne \emptyset \text{ for } k=1,\dots,n
\} \, ;
\end{equation}
here $ [ s_k, t_k ) \subset \R $ are usual (rather than coarse) time
intervals. (Endpoints may be neglected, as we will see soon.)
The diameter of the cell \eqref{3c4} (w.r.t.\ the metric \eqref{3c2})
does not exceed $ \max_k ( t_k - s_k ) $. Thus, we get weak
convergence of measures, which proves the following result.

\begin{theorem}\label{3c5}
For every dyadic coarse factorization $ \(
(\Om[i],\F[i],P[i])_{i=1}^\infty, \A \) $ and every $ x \in L_2 (\A)
$, the sequence $ \( \mu_x [i] \)_{i=1}^\infty $ of spectral measures
converges weakly to a (finite, positive) measure $ \mu_x [\infty] $ on
the Polish space $ \cC[\infty] $.
\end{theorem}

Convergence of measures $ \mu_x[i] $ on a `cell' of the form
\eqref{3c3} (or \eqref{3c4}) does not ensure that the limit is $
\mu_x[\infty] $ on the `cell'.\footnote{%
 Think for example about an atom at the point $ \frac1n $ of $ \R $,
 and `cells' of the form $ (x,y] $.}
Rather, the limit lies between $ \mu_x[\infty] $-measures of the
interior and the closure of the cell,
\begin{multline}\label{3c6}
\mu_x[\infty] \( \{ M \in \cC[\infty] : M \subset (s_1,t_1) \cup
 \dotsc \cup (s_n,t_n) \} \) \\
\le \| \cE{ x[\infty] }{
 \F_{(s_1,t_1)\cup\dotsc\cup(s_n,t_n)} } \|^2 \\
\le \mu_x[\infty] \( \{ M \in \cC[\infty] : M \subset [s_1,t_1] \cup
 \dotsc \cup [s_n,t_n] \} \) \, .
\end{multline}

\begin{lemma}\label{3c7}
For every $ t \in \R $,
\[
\mu_x[\infty] \( \{ M \in \cC[\infty] : M \ni t \} \) = 0 \, .
\]
\end{lemma}

\begin{proof}
Lemma \ref{3b10} gives us
\[
\| \cE{ x[\infty] }{ \F_{(-\infty,-\eps)\cup(\eps,+\infty)} }
\|^2 \tendseps \| x[\infty] \|^2 \, ;
\]
therefore
\[
\mu_x[\infty] \( \{ M \in \cC[\infty] : M \subset (-\infty,\eps] \cup
[\eps,+\infty) \} \) \tendseps \mu_x[\infty] \( \cC[\infty] \) \, .
\]
\qqed\end{proof}

Applying Fubini's theorem we see that $ \mu_x[\infty] $ is
concentrated on (the set of all) compact sets $ M $ of Lebesgue
measure $ 0 $ (therefore, nowhere dense).

Due to \ref{3c7} we see that the boundary of a `cell' is
negligible (of measure $ 0 $); inequalities \eqref{3c6} are, in fact,
equalities. So,
\begin{equation}\label{22}
\mu_x[\infty] \( \{ M \in \cC[\infty] : M \subset E \} \) =
\| \cE{ x[\infty] }{ \F_E } \|^2 \, ,
\end{equation}
where $ E \subset \R $ is an arbitrary elementary
set,\index{elementary set} that is, a
finite union of intervals (treated modulo finite sets), $ E =
(s_1,t_1) \cup \dotsc \cup (s_n,t_n) $, and $ \F_E = \F_{s_1,t_1}
\otimes \dots \otimes \F_{s_n,t_n} $.\index{zzFE@$ \F_E $, sub-\sif}

For a finite $ i $, the Fourier-Walsh basis decomposes $ L_2[i] $ into
one-dimensional subspaces indexed by $ M \in \cC[i] $, and each subset
$ \cM \subset \cC[i] $ leads to a subspace $ H_\cM $ of $ L_2[i] $
spanned by $ \tau_M $, $ M \in \cM $. In particular, for a subset of
the form $ \cM_E = \{ M \in \cC[i] : M \subset E \} $ we have $
H_{\cM_E} = L_2 (\Om[i],\F_E[i],P[i]) $.

Similarly, for the limiting object, the subspace $ H_{\cM_E} = L_2
(\Om,\F_E,P) $ of $ L_2[\infty] $ corresponds
to the set $ \cM_E = \{ M \in \cC[\infty] : M \subset E \} $. In
\ref{sec:3.4} a subspace $ H_\cM \subset L_2[\infty] $ will be
defined for every Borel set $ \cM \subset \cC[\infty] $.

\subsection{The limiting object}
\label{sec:3.4}

\begin{definition}\label{3d1}
A \emph{continuous
factorization}\index{continuous factorization}\index{factorization}
\textup{(}of probability
spaces, over $ \R $\textup{)} consists of a probability space $
(\Om,\F,P) $ and a two-parameter family $ (\F_{s,t})_{s\le t} $ of
sub-\sif s $ \F_{s,t} \subset \F $ such that\footnote{%
 Here $ r,s,t $ are real numbers; coarse instants are not used in
 \ref{sec:3.4}, \ref{sec:3.5}.}
\[
\F_{r,t} = \F_{r,s} \otimes \F_{s,t} \quad  \text{whenever } r \le s
\le t \leqno\text{\textup{(a)}}
\]
\textup{(}that is, $ \F_{r,s} $ and $ \F_{s,t} $ are independent, and
together generate $ \F_{r,t} $\textup{),}
\[
\bigcup_{\eps>0} \F_{s+\eps,t-\eps} \text{ generates $ \F_{s,t} $
whenever $ s < t $,} \leqno\text{\textup{(b)}}
\]
and
\[
\bigcup_{n=1}^\infty \F_{-n,n}\text{ generates } \F \,
. \leqno\text{\textup{(c)}}
\]
\end{definition}

The refinement of any dyadic coarse factorization is a continuous
factorization (as was shown in \ref{sec:3.2}).

\begin{definition}\label{3.17}
Let $ \( (\Om,\F,P), (\F_{s,t})_{s\le t} \) $ be a continuous
factorization, and $ x \in L_2 (\Om,\F,P) $. The \emph{spectral
measure}\index{spectral measure (continuous case)}
$ \mu_x $\index{zzmux@$ \mu_x $, spectral measure (continuous)}
of $ x $ is the \textup{(}finite, positive\textup{)} measure
on the space $ \cC = \cC[\infty] $ of compact subsets of $ \R $ such
that
\[
\mu_x \( \{ M \in \cC : M \subset E \} \) = \| \cE{ x }{ \F_E } \|^2
\]
for all elementary sets $ E \subset \R $.
\end{definition}

Uniqueness of $ \mu_x $ is checked easily. Existence of $ \mu_x $ is
proven in \ref{sec:3.3} by discrete approximation, assuming that
the continuous factorization is the refinement of a dyadic coarse
factorization. Another proof, without approximation, will be given
by \ref{3.22}.

The spectral measure is concentrated on (the set of all) nowhere dense
compact sets, and
\begin{equation}\label{225}
\mu_x \( \{ M \in \cC : M \ni t \} \) = 0 \quad \text{for each
} t \in \R \, ,
\end{equation}
which follows from \ref{3.19} for $ s=t $, since $ \F_{t,t} =
\F_{t,t} \otimes \F_{t,t} $ is degenerate.

\begin{example}\label{3d2}

The refinement of the Brownian coarse factorization (see
\ref{3b3}) is the Brownian continuous factorization,
\[
\F_{s,t} \text{ is generated by } \{ B(v) - B(u) : s \le u \le v \le t
\} \, ,
\]
where $ B(\cdot) $ is the usual Brownian motion. Every $ x \in L_2 $
admits It\^o's decomposition into multiple stochastic integrals,
\begin{multline*}
x = \hat x ( \emptyset ) + \int \hat x ( \{ t_1 \} ) \, \D B(t_1) +
 \iint\limits_{t_1<t_2} \hat x ( \{t_1,t_2\} ) \, \D B(t_1) \D B(t_2)
 + \dots \\
= \sum_{n=0}^\infty \, \idotsint\limits_{t_1<\dots<t_n} \hat x ( \{
 t_1,\dots,t_n \} ) \, \D B(t_1) \dots \D B(t_n) \, ,
\end{multline*}
where $ \hat x \in L_2 (\Cfinite) $,
$ \Cfinite $\index{zzCf@$ \Cfinite $, space of finite sets}
being the space of
all finite subsets of $ \R $, equipped with the natural (Lebesgue)
measure, making the transform $ x \leftrightarrow \hat x $ unitary,
according to the formula
\begin{multline*}
\Ex |x|^2 = | \hat x ( \emptyset ) |^2 + \int | \hat x ( \{ t_1 \} )
|^2 \, \D t_1 + 
 \iint\limits_{t_1<t_2} | \hat x ( \{t_1,t_2\} ) |^2 \, \D t_1 \D t_2
 + \dots \\
= \sum_{n=0}^\infty \, \idotsint\limits_{t_1<\dots<t_n} | \hat x ( \{
 t_1,\dots,t_n \} ) |^2 \, \D t_1 \dots \D t_n \, .
\end{multline*}
The spectral measure $ \mu_x $ of $ x $ is
\[
\mu_x (A) = \sum_{n=0}^\infty \;
\idotsint\limits_{t_1<\dots<t_n,\{t_1,\dots,t_n\}\in A} | \hat x ( \{
 t_1,\dots,t_n \} ) |^2 \, \D t_1 \dots \D t_n \, . 
\]
This is an important property of the Brownian continuous
factorization: the spectral measure (of any random variable) is
concentrated on the subset $ \Cfinite \subset \cC $, and absolutely
continuous w.r.t.\ the Lebesgue measure on $ \Cfinite $.

In particular, for $ x = \exp \( \I \sqrt\la B(t) \) $ the measure $
\mu_x $ is just the distribution of the Poisson process of rate $ \la
$ on $ (0,t) $. Indeed,
\[
\exp \( \I \sqrt\la B(t) \) = \E^{-\la t/2} \sum_{n=0}^\infty \la^{n/2}
\idotsint\limits_{0<t_1<\dots<t_n<t} \, \D B(t_1) \dots \D B(t_n) \, .
\]

\end{example}

\begin{example}\label{3d3}

\begin{sloppypar}
Recall the process $ Y_\eps $ of \ref{1b3};
\[
Y_\eps (t) = \exp \( \I B(\ln t) - \I B(\ln\eps) \) \, .
\]
We define $ \F_{s,t} $ as the \sif\ generated by `multiplicative
increments' $ Y_\eps(v) / Y_\eps(u) $ for all $ (u,v) \subset
(s,t) $, that is, by (usual) Brownian increments on $ (\ln s, \ln t)
$. The spectral measure $ \mu_{Y_\eps(t)} $ is the distribution of a
non-homogeneous Poisson process on $ (\eps,t) $, the image of the
usual Poisson process (of rate $ 1 $) on $ (\ln\eps, \ln t) $ under
the time change $ u \mapsto \E^u $. The rate of the non-homogeneous
Poisson process is $ \la(s) = 1/s $.
\end{sloppypar}

The limiting process $ Y $ was discussed in \ref{1b3}. It may be
treated as the refinement of $ Y_\eps $ for $ \eps \to 0 $ (I leave the
details to the reader). The spectral measure $ \mu_{Y(t)} $ should be
the distribution of a non-homogeneous Poisson process on $ (0,t) $, at
the rate $ \la(s) = 1/s $. Random points accumulate to $ 0 $; we add $
0 $ to the random set, making it compact. However, the equality $ \mu
( \{ M : M \ni 0 \} ) = 1 $ does not conform to \ref{3c7}! It
happens because the limiting object is not a \emph{continuous}
factorization. Denote by $ \F_{0+,1} $ the \sif\ generated by $
\cup_{\eps>0} \F_{\eps,1} $. Every $ Y(1) / Y(t) $ for $ t > 0 $ is $
\F_{0+,1} $-measurable, but $ Y(1) $ is not. The global phase is
missing. Of course, for every $ t>0 $, there exists an independent
complement of $ \F_{0+,t} $ in $ \F_{-\infty,t} $ (for example, the
\sif\ generated by $ Y(t) $). However, we cannot choose a single
complement (to be denoted by $ \F_{-\infty,0+} $) for all $ t>0 $,
since the tail \sif\ $ \cap_{t>0} \F_{-\infty,t} $ is degenerate.

\end{example}

\begin{lemma}\label{3.19}

For every continuous factorization $ \( (\Om,\F,P), (\F_{s,t})_{s\le
t} \) $ and every $ s \le t $,
\[
\F_{s,t} = \bigcap_{\eps>0} \F_{s-\eps,t+\eps} \, .
\]

\end{lemma}

\begin{proof}
\begin{sloppypar}
The \sif\ $ \cap_{\eps>0} \F_{0,\eps} $ is degenerate by Kolmogorov's
zero-one law applied to $ \F_{1,\infty}, \F_{1/2,1}, \F_{1/3,1/2},
\dots \, $ Further, $ \F_{-\infty,\eps} = \F_{-\infty,0} \otimes
\F_{0,\eps} \tendseps \F_{-\infty,0} $. Though the equality $ \lim (
\A \vee \B_n ) = \A \vee ( \lim \B_n ) $ does not hold in general, it
does hold for independent $ \A $ and $ \B_1 $ ($ \B_1 \supset \B_2
\supset \dots $), which is a rather trivial part of Weizs\"acker's
criteria \cite{We}. The rest of the proof is left to the reader.
\qqed\end{sloppypar}
\end{proof}

The theory of direct integrals of Hilbert spaces may be used on the
way to Theorem \ref{3.25}. In fact, I did so in
\cite[Th.~2.3]{Ts98}. Here, however, I choose a self-contained
presentation.
First, a general result of measure theory, useful for proving the
existence of $ \mu_x $ (without dyadic approximation).

\begin{lemma}\label{3.20}
Let $ X $ be a compact topological space, $ \A $ an algebra
of subsets of $ X $, and $ \mu : \A \to [0,\infty) $ an additive
function satisfying the following regularity condition:

For every $ A \in \A $ and $ \eps > 0 $ there exists
$ B \in \A $ such that $ \overline B \subset A $ \textup{(here $
\overline B $ is the closure of $ B $)} and $ \mu(B) \ge \mu(A) - \eps
$.

Then $ \mu $ has a unique extension to a measure on the \sif\
generated by $ \A $.
\end{lemma}

\begin{proof}
Due to a well-known theorem, it is enough to prove that $ \mu $ is $
\sigma $-additive on $ \A $. Let $ A_1 \supset A_2 \supset \dots $, $
A_1, A_2, \dots \in \A $, $ \cap A_k = \emptyset $; we have to prove
that $ \mu (A_k) \to 0 $. Given $ \eps > 0 $, we can choose $ B_k \in
\A $ such that $ \overline B_k \subset A_k $ and $ \mu (B_k) \ge \mu
(A_k) -
2^{-k} \eps $. Due to compactness, the relation $ \cap \overline B_k
\subset \cap A_k = \emptyset $ implies $ \overline B_1 \cap \dots \cap
\overline B_n = \emptyset $ for some $ n $. Thus, $ \mu (A_n) = \mu (
A_1 \cap \dots \cap A_n ) \le \mu ( B_1 \cap \dots \cap B_n ) + \mu (
A_1 \setminus B_1 ) + \dots + \mu ( A_n \setminus B_n ) < \eps $.
\qqed\end{proof}

\begin{remark}\label{3.21}
All $ A \in \A $ such that $ A $ and $ X \setminus A $ both satisfy
the regularity condition, are a subalgebra of $ \A $. (The proof is
left to the reader.) Therefore it is enough to check the condition for
$ A $ and $ X \setminus A $ where $ A $ runs over a set that generates
the algebra $ \A $.
\end{remark}

\begin{lemma}\label{3.22}
The spectral measure $ \mu_x $ exists for every $ x \in L_2 (\Om,\F,P)
$ and every continuous factorization $ (\F_{s,t})_{s\le t} $.
\end{lemma}

\begin{proof}
First, compactness. We have $ \| \cE{ x }{ \F_{-m,m} } \|^2 \to
\|x\|^2 $ for $ m \to \infty $ by \ref{3d1}(c); thus we may restrict
ourselves to $ x $ measurable w.r.t.\ $ \F_{-m,m} $ for some $ m
$. The corresponding part $ \cC_m = \{ M \in \cC : M \subset [-m,m] \}
$ of $ \cC $ is compact.

Second, additivity on an algebra. We have an algebra $ \A $ of subsets
of $ \cC_m $, generated by `cells' of the form \eqref{3c4}. Such a
cell leads to a subspace of $ L_2 (\Om,\F_{-m,m},P) $ spanned by
products $ f_1 \dots f_n $ where each $ f_k $ is measurable w.r.t.\ $
\F_{s_k,t_k} $, square integrable, and $ \Ex f_k = 0 $. A partition of
the interval $ [-m,m] $ into $ n $ subintervals leads to a partition
of $ \cC_m $ into $ 2^n $ parts, and a decomposition of $ L_2
(\Om,\F_{-m,m},P) $ into $ 2^n $ orthogonal subspaces. Thus, $ x $
decomposes into $ 2^n $ orthogonal vectors; their squared norms give
us $ \mu_x $ on a finite subalgebra (of cardinality $ 2^{2^n} $) of $
\A $. We see that $ \mu_x $ is additive on such subalgebras. Their
union (over all partitions of $ [-m,m] $) is the whole $ \A $, and any
two of them are contained in some third; therefore, $ \mu_x $ is
additive on $ \A $.

Third, regularity (required by \ref{3.20}). Due to
\ref{3.21}, regularity may be checked only for sets $ A_E = \{ M \in
\cC_m : M \subset E \} $ and $ \cC_m \setminus A_E $. It follows
easily from \ref{3d1}(b) and \ref{3.19}.
\qqed\end{proof}

\begin{remark}
In the proof of \ref{3.22}, an \emph{orthogonal
decomposition}\index{orthogonal decomposition}
of the Hilbert space $ H = L_2 (\Om,\F,P) $ over the algebra $ \A $ is
constructed; that is, a family $ (H_A)_{A\in\A} $ of (closed linear)
subspaces $ H_A \subset H $ such that $ H_{A\cup B} = H_A \oplus H_B
$\index{zzz@$ \oplus $, orthogonal sum of Hilbert spaces}
(it means that $ H_A $ and $ H_B $ are orthogonal, and their sum is $
H_{A\cup B} $) whenever $ A \cap B = \emptyset $, and $ H_\cC = H
$. The decomposition satisfies
\[
H_{\cM_E} = L_2 (\Om,\F_E,P) \, ,
\]
where $ \cM_E = \{ M \in \cC : M \subset E \} $, and is uniquely
determined by this property.
\end{remark}

The following general result will help us construct $ H_\cM $ for
all Borel sets $ \cM \subset \cC $.

\begin{lemma}\label{3.24}
Let $ X $ be a set, $ \A $ an algebra of subsets of $ X $, $ H $ a
Hilbert space, and $ (H_A)_{A\in\A} $ an orthogonal decomposition of $
H $ over $ \A $. Assume that for every $ x \in H $ the additive
function\footnote{%
 Here $ \Proj_{H_A} $ is the orthogonal projection $ H \to H_A $.}
$ A \mapsto \| \Proj_{H_A} x \|^2 $ on $ \A $ can be extended to a
measure on the \sif\ $ \si(\A) $ generated by $ \A $. Then the
orthogonal decomposition can be extended to an orthogonal
decomposition $ (H_B)_{B\in\si(\A)} $,
$ \si $-additive\index{si@$ \si $-additive orthogonal decomposition}
in the sense that\footnote{%
 That is, $ H_{B_1\cup B_2\cup\dots} $ is the closure of the algebraic
 sum of $ H_{B_k} $.}
$ H_{B_1\cup B_2\cup\dots} = H_{B_1} \oplus H_{B_2} \oplus \dots
$ whenever $ B_1, B_2, \dots \in \si(\A) $ are pairwise disjoint.
\end{lemma}

\begin{proof}
The extension of the additive function $ \mu_x : \A \to [0,\infty) $,
$ \mu_x (A) = \| \Proj_{H_A} x \|^2 $, to a measure on $ \si(\A) $ is
unique; denote it by $ \mu_x $ again. Consider the set of all $ B \in
\si(\A) $ such that there exists a subspace $ H_B \subset H $
satisfying $ \| \Proj_{H_B} x \|^2 = \mu_x (B) $
for all $ x \in H $. The set contains $ \A $, and is a monotone
class (that is, closed under the limit of monotone sequences), which
is easy to check. Therefore the set is the whole $ \si(\A) $.
\qqed\end{proof}

Combining \ref{3.22} and \ref{3.24} we conclude.

\begin{theorem}\label{3.25}
\begin{sloppypar}
For every continuous factorization $ \( (\Om,\F,P), (\F_{s,t})_{s\le t}
\) $ there exists one and only one $ \si $-additive orthogonal
decomposition $ (H_\cM) $\index{zzHM@$ H_\cM $, subspace}
of the Hilbert space $ L_2 (\Om,\F,P) $ over
the Borel \sif\ of the space $ \cC $ \textup{(}of compact subsets of $
\R $\textup{)} such that $ H_{\cM_E} = L_2(\Om,\F_E,P) $ for every
elementary set $ E \subset \R $ \textup{(}that is, a finite union of
intervals\textup{)}; here $ \cM_E = \{ M \in \cC : M \subset E \}
$. The orthogonal decomposition is related to spectral measures by
\begin{equation}\label{326}
\| \Proj_{H_\cM} f \|^2 = \mu_f (\cM)
\end{equation}
for all $ f \in L_2 (\Om,\F,P) $ and all Borel sets $ \cM \subset \cC
$.
\end{sloppypar}
\end{theorem}

\subsection{Time shift; noise}
\label{sec:3.5}

Let $ \( (\Om[i],\F[i],P[i])_{i=1}^\infty, \A \) $ be a dyadic coarse
factorization. For each $ i $ the lattice $ \frac1i \Z $ acts on $
\Om[i] $ by measure preserving transformations $ \al_t : \Om[i] \to
\Om[i] $ (time shift),
\[
\al_t (\om) (s) = \om (s-t) \quad \text{for all } s \in \frac1i \Z \,
.
\]
For each coarse instant $ t = (t[i])_{i=1}^\infty $ we have a map $
\al_t : \Om [\all] \to \Om [\all] $,
\[
\al_t (\om) [i] (s) = \om [i] (s-t[i]) \quad \text{for all } s \in
\frac1i \Z \, .
\]
Such $ \al_t $ is an automorphism of the dyadic coarse sample space,
but the coarse \sif\ $ \A $ need not be invariant under $ \al_t $. We
consider such a condition:
\begin{equation}\label{3e1}
\text{$ \A $ is invariant under $ \al_t $ for every coarse instant $ t
$.}
\end{equation}
Dyadic coarse factorizations of \ref{3b3}, \ref{3b5},
\ref{3b6}, \ref{3b7} satisfy \eqref{3e1}, but that of
\ref{3b4} does not.

If \eqref{3e1} is satisfied, then the refinement $ \al_t[\infty] =
\Lim_{i\to\infty,\A} \al_t[i] $ is an automorphism of the refinement $
(\Om,\F,P) $ of the dyadic coarse factorization. Existence of the
limit for \emph{every} converging sequence $ t = (t[i]) $ implies that
$ \al_t [\infty] $ depends on $ t[\infty] $ only (see \ref{3e4}
below), and we get a one-parameter group $ (\al_t)_{t\in\R} $ of
automorphisms (that is, invertible measure preserving transformations
$ \bmod \, 0 $) of $ (\Om,\F,P) $. The group is continuous in the
sense that $ \Pr{ A \bigtriangleup \al_t(A) } { \tendst 0 } $ for all
$ A \in \F $, which is ensured by \eqref{3e1} (see \ref{3e4}
again).

\begin{definition}

A \emph{noise}\index{noise} $ \( (\Om,\F,P), (\F_{s,t})_{s\le t},
(\al_t)_{t\in\R}
\) $ consists of a continuous factorization $ \( (\Om,\F,P),
(\F_{s,t})_{s\le t} \) $ and a one-parameter group of automorphisms $
\al_t $ of $ (\Om,\F,P) $ such that
\begin{gather*}
\al_t^{-1} (\F_{r,s}) = \F_{r-t,s-t} \quad \text{for all } r,s,t \in
 \R, \, r \le s \, , \\
P \( A \bigtriangleup \al_t^{-1} (A) \) \tendst 0 \quad \text{for all
} A \in \F \, .
\end{gather*}

\end{definition}

Unfortunately, the latter assumption (continuity of the group action)
is missing in my former publications, which opens the door for
pathologies.\footnote{%
 Most results of these former publications do not depend on the (missing)
 continuity condition. But anyway, a discontinuous group action is a
 pathology, no doubt. (In particular, it cannot be Borel measurable.)
 The proof of Lemma 2.9 of \cite{Ts98}, based on Weyl's relation,
 depends on the continuity condition.}

\begin{remark}

Continuity of the factorization follows from other assumptions, see
\cite[Lemma 2.1]{Ts98}. For arbitrary factorizations, continuity is
restrictive (recall \ref{3d3}); waiving it, we get discontinuity
points $ t \in \R $ which are a finite or countable set. For a noise,
however, the set is invariant under time shifts, and therefore, empty.

\end{remark}

\begin{lemma}\label{3e4}
For every dyadic coarse factorization satisfying \eqref{3e1}, its
refinement is a noise.
\end{lemma}

\begin{proof}
Our first argument parallels the proof of \ref{3b8}. Namely, let
$ s,t $ be two coarse instants such that $ s[\infty] = t[\infty] $. We
introduce a coarse event $ r $:
\[
r[i] = \begin{cases}
 s[i] &\text{for $ i $ even},\\
 t[i] &\text{for $ i $ odd}.
\end{cases}
\]
We have
\[
\Lim \al_s [i] = \Lim \al_s [2i] = \Lim \al_r [2i] = \Lim \al_r [i] \,
.
\]
Similarly, $ \Lim \al_t [i] = \Lim \al_r [i] $. Thus, $ \Lim \al_s [i]
= \Lim \al_t [i] $, and we may define a one-parameter group of
automorphisms $ (\al_t)_{t\in\R} $ on $ (\Om,\F,P) $ by $
\al_{t[\infty]} = \Lim \al_t[i] $.

Our second argument resembles the proof of \ref{3b10}. Namely,
assume existence of $ A_\infty \in \F $, $ \eps > 0 $ and $ t_n \to 0
$ such that $ P \( A_\infty \bigtriangleup \al^{-1}_{t_n} (A_\infty)
\) \ge \eps $ for all $ n $. We choose a coarse event $ A \in \A $
such that $ A[\infty] = A_\infty $, and coarse instants $ s_n $ such
that $ s_n[\infty] = t_n $ for all $ n $. Taking into account that $
P[i] \( A[i] \bigtriangleup \al^{-1}_{s_n} [i] A[i] \) \to P \(
A_\infty \bigtriangleup \al^{-1}_{t_n} (A_\infty) \) \ge \eps $ and $
s_n[i] \to t_n $ when $ i \to \infty $, we choose integers $ i_1 < i_2
< \dots $ such that $ P[i] \( A[i] \bigtriangleup \al^{-1}_{s_n} [i]
A[i] \) \ge \eps/2 $ and $ | s_n[i] | \le |t_n|+1/n $ whenever $ i \ge
i_n $. We define a coarse instant $ r $ by $ r[i] = s_n[i] $ whenever
$ i_n \le i < i_{n+1} $. Clearly, $ r[\infty] = 0 $; therefore $ \Lim
\al^{-1}_r [i] A[i] = \al_0^{-1} A[\infty] = A[\infty] $, and $ P[i]
\( A[i] \bigtriangleup \al^{-1}_r [i] A[i] \) \to 0 $, which is
impossible: these probabilities exceed $ \eps/2 $. The contradiction
proves continuity of the group $ (\al_t)_{t\in\R} $.
\qqed\end{proof}

\begin{question}
Is every noise the refinement of some dyadic coarse factorization
satisfying \eqref{3e1}? I do not know; I guess that the answer is
negative. It would be interesting to find some special features of
such refinements among all noises. It is also unclear what happens to
the class of such refinements, if subsequences are permitted (like in
\ref{3b6}).
\end{question}

\section{Example: The Noise Made by a Poisson Snake}
\label{sec:4}
This section is based on a paper by J.~Warren entitled ``The noise
made by a Poisson snake'' \cite{Wa}.

\subsection{Three discrete semigroups: algebraic definition}
\label{sec:4.1}

A discrete semigroup (with unit; non-commutative, in general) may be
defined by generators and relations.

Two generators $ f_+, f_- $ with two relations $ f_+ f_- = 1 $, $ f_-
f_+ = 1 $ generate a semigroup $ G_1^\discrete $ that is in fact a
group, just the cyclic group $ \Z $. Indeed, every word reduces to
some $ f_+^k $ or $ f_-^k $ (or $ 1 $).

Two generators $ f_+, f_- $ with a single relation $ f_+ f_- = 1 $
generate a semigroup $ G_2^\discrete $. Every word reduces to some $
f_-^k f_+^l $. The composition is
\begin{equation}\label{4a1}
( f_-^{k_1} f_+^{l_1} ) ( f_-^{k_2} f_+^{l_2} ) = f_-^k f_+^l \, ,
\quad 
\begin{aligned}
 k &= k_1 + \max (0,k_2-l_1) \, , \\
 l &= l_2 + \max (0,l_1-k_2) \, .
\end{aligned}
\end{equation}
The canonical homomorphism $ G_2^\discrete \to G_1^\discrete $ maps $
f_+ $ to $ f_+ $, $ f_- $ to $ f_- $, and $ f_-^k f_+^l $ into $
f_-^{k-l} $ (if $ k>l $), $ f_+^{l-k} $ (if $ k<l $), or $ 1 $ (if
$ k=l $). Accordingly, the composition law \eqref{4a1} satisfies
\[
l-k = (l_1-k_1) + (l_2-k_2) \, .
\]
There is a more convenient pair of parameters, $ a=l-k $, $ b=k $;
that is,\footnote{%
 Parameters $ a,b $ of \eqref{4a2} and $ a,b,c $ of \eqref{4a4} are
 suggested by S.~Watanabe.}
\begin{equation}\label{4a2}
\begin{gathered}
f_{a,b} = f_-^b f_+^{a+b} \quad \text{for } a,b \in \Z, \, b \ge 0 ,
 \, a+b \ge 0 \, ; \\
f_{a_1,b_1} f_{a_2,b_2} = f_{a,b} \, , \quad
\begin{aligned}
 a &= a_1 + a_2 \, , \\
 b &= \max ( b_1, b_2-a_1 ) \, .
\end{aligned}
\end{gathered}
\end{equation}
The canonical homomorphism $ G_2^\discrete \to G_1^\discrete $ maps $
f_{a,b} $ to $ f_a $, where $ f_a \in G_1^\discrete $ is $ f_+^a $ for
$ a>0 $, $ f_-^{|a|} $ for $ a<0 $, and $ 1 $ for $ a=0 $.

Three generators $ f_-, f_+, f_* $ with three relations
\begin{equation}\label{4a3}
f_+ f_- = 1 \, , \quad f_* f_- = 1 \, , \quad f_* f_+ = f_* f_*
\end{equation}
generate a semigroup $ G_3^\discrete $. Every word reduces to some $
f_-^k f_+^l f_*^m $. The following homomorphism $ G_3^\discrete \to
G_2^\discrete $ will be called canonical: $ f_- \mapsto f_- $, $ f_+
\mapsto f_+ $, $ f_* \mapsto f_+ $. We have $ f_-^k f_+^l f_*^m
\mapsto f_-^k f_+^{l+m} $, which suggests such a triple of parameters
for $ G_3^\discrete $: $ a = l+m-k $, $ b=k $, $ c=m $; that is,
\begin{equation}\label{4a4}
\begin{gathered}
f_{a,b,c} = f_-^b f_+^{a+b-c} f_*^c \quad \text{for } a,b,c \in \Z, \;
 b \ge 0 , \; 0 \le c \le a+b \, ; \\
f_{a_1,b_1,c_1} f_{a_2,b_2,c_2} = f_{a,b,c} \, , \quad
\begin{aligned}
 a &= a_1 + a_2 \, , \\
 b &= \max ( b_1, b_2-a_1 ) \, ,
\end{aligned}
\quad
 c = \begin{cases}
  a_2 + c_1 &\text{if $ c_1 > b_2 $},\\
  c_2 &\text{otherwise}.
 \end{cases}
\end{gathered}
\end{equation}
The canonical homomorphism $ G_3^\discrete \to G_2^\discrete $ is just
$ f_{a,b,c} \mapsto f_{a,b} $.

Note that $ G_1^\discrete $ is commutative, but $ G_2^\discrete $ and
$ G_3^\discrete $ are not.

\subsection{The three discrete semigroups: representation}
\label{sec:4.2}

By a representation of a semigroup $ G $ on a set $ S $ we mean a map
$ G \times S \ni (g,s) \mapsto g(s) \in S $ such that
\[
(g_1 g_2) (s) = g_2 \( g_1(s) \) \quad \text{and} \quad 1(s) = s
\]
for all $ g_1, g_2 \in G $, $ s \in S $. The representation is called
faithful, if
\[
g_1 \ne g_2 \imply \exists s \in S \;\; \( g_1(s) \ne g_2(s) \) \, .
\]
Every $ G $ has a faithful representation on itself, $ S=G $, namely,
the regular representation, $ g(g_0) = g_0 g $. Fortunately, $
G_2^\discrete $ and $ G_3^\discrete $ have more economical faithful
representations on the set $ \Z_+ = \{ 0,1,2,\dots \} $. Namely, for
$ G_2^\discrete $,
\begin{equation}\label{4b1}
\begin{gathered}\includegraphics{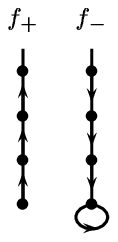}\end{gathered}
\qquad
\begin{gathered}
f_+ (x) = x+1 \, , \quad f_- (x) = \max ( 0, x-1 ) \, , \\
f_{a,b} (x) = a + \max (x,b) \, ,
\end{gathered}
\quad
\begin{gathered}\includegraphics{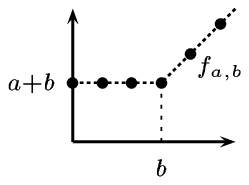}\end{gathered}
\end{equation}
$ x \in \Z_+ $. For $ G_3^\discrete $,
\begin{equation}\label{4b2}
\begin{gathered}
 \includegraphics{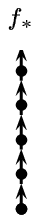}\hspace{-2mm}
 \includegraphics{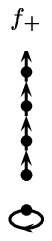}\hspace{-2mm}
 \includegraphics{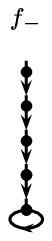}
\end{gathered}
\quad
\begin{gathered}
f_* (x) = x+1 \, , \quad f_- (x) = \max ( 0, x-1 ) \, , \\
f_+ (x) = \begin{cases}
 x+1 &\text{ for $ x>0 $},\\
 0 &\text{for $ x=0 $};
\end{cases} \\
f_{a,b,c} (x) = \begin{cases}
 c &\text{for $ 0 \le x \le b $}, \\
 x+a &\text{for $ x > b $}.
\end{cases}
\end{gathered}
\hspace{-1mm}
\begin{gathered}\includegraphics{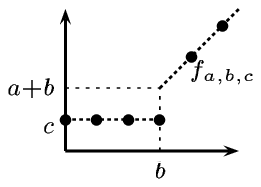}\end{gathered}
\end{equation}

\subsection{Random walks and stochastic flows in discrete semigroups}
\label{sec:4.3}

\begin{example}\label{4ce1}

The standard random walk on $ \Z $ may be described by $ G_1^\discrete
$\nobreakdash-\hspace{0pt}valued random variables
\begin{equation}\label{4c1}
\begin{gathered}
\xi_{s,t} = \xi_{s,s+1} \xi_{s+1,s+2} \dots \xi_{t-1,t} \quad
\text{for } s,t \in \Z, \, s \le t \, ; \\
\xi_{t,t+1} \text{ are independent random variables } (t\in\Z) \, ; \\
\Pr{ \xi_{t,t+1} = f_- } = \frac12 = \Pr{ \xi_{t,t+1} = f_+ } \quad
\text{for each } t \in \Z \, .
\end{gathered}
\end{equation}
Note that $ \xi_{r,s} \xi_{s,t} = \xi_{r,t} $ whenever $ r \le s \le t
$. Everyone knows that
\begin{equation}\label{4c1a}
\Pr{ \xi_{0,t} = f_a } = \frac1{2^t} \binom{ t }{ \frac{t+a}2 }
\end{equation}
for $ a = -t, -t+2, -t+4, \dots, t $.

In fact, `the standard random walk' is the random process $ t \mapsto
\xi_{0,t} $. Taking into account that $ G_1^\discrete $ is a group, $
\xi_{s,t} $ may be thought of as an increment, $ \xi_{s,t} =
\xi_{0,s}^{-1} \xi_{s,t} $.

\end{example}

\begin{example}\label{4ce2}

Formulas \eqref{4c1} work equally well on $ G_2^\discrete $. Still, $
\xi_{r,s} \xi_{s,t} = \xi_{r,t} $. However, $ G_2^\discrete $ is not a
group, and $ \xi_{s,t} $ is not an increment; moreover, it is not a
function of $ \xi_{0,s} $ and $ \xi_{0,t} $. Indeed, knowing $ a_1,b_1
$ and $ a_1+a_2 $, $ \max ( b_1, b_2-a_1 ) $ (recall \eqref{4a2}) we
can find $ a_2 $ but not $ b_2 $. Thus, the two-parameter family $
(\xi_{s,t})_{s\le t} $ of random variables is more than just a random
walk. Let us call such a family an \emph{abstract stochastic
flow.}\index{abstract stochastic flow (discrete)}
Why `abstract'? Since $ G_2^\discrete $ is an abstract semigroup
rather than a semigroup of transformations (of some set). So, we have
the standard abstract flow in $ G_2^\discrete $. In order to
get a (usual, not abstract) stochastic
flow,\index{stochastic flow (discrete)}
we have to choose a
representation of $ G_2^\discrete $. Of course, the regular
representation could be used, but the representation \eqref{4b1} is
more useful. Introducing integer-valued random variables $ a(s,t),
b(s,t) $ by
\[
\xi_{s,t} = f_{a(s,t),b(s,t)}
\]
we express the stochastic flow as
\[
\xi_{s,t} (x) = a(s,t) + \max ( x, b(s,t) ) \, .
\]
Fixing $ s $ and $ x $ we get a random process called a single-point
motion of the flow. Namely, it is a reflecting random
walk. Especially, for $ s=0 $ and $ x=0 $, the process
\[
t \mapsto \xi_{0,t} (0) = a(0,t) + b(0,t)
\]
is a reflecting random walk. It is easy to see that two processes
\begin{align*}
t &\mapsto \xi_{0,t} (0) = a(0,t) + b(0,t) \, , \\
t &\mapsto \Big| a(0,t) + \frac12 \Big| - \frac12
\qquad\begin{gathered}\includegraphics{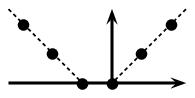}\end{gathered}
\end{align*}
are identically distributed. Also,
\begin{equation}\label{4c1c}
\begin{aligned}
b(0,t) &= - \min_{s=0,1,\dots,t} a(0,s) \, , \\
a(0,t) + b(0,t) &= \max_{s=0,1,\dots,t} a(s,t) \, , \\
\end{aligned} \qquad
\begin{gathered}\includegraphics{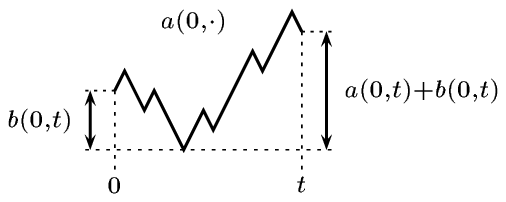}\end{gathered}
\end{equation}
and $ a(\cdot,\cdot) $ is the standard random walk on $ G_1^\discrete
= \Z $. That is, the canonical homomorphism $ G_2^\discrete \to
G_1^\discrete $ transforms the standard flow on $ G_2^\discrete
$ into the standard flow (or random walk) on $ G_1^\discrete $. Using
the reflection principle, one gets
\begin{equation}\label{4c2}
\Pr{ \xi_{0,t} = f_{a,b} } = \frac{ a+2b+1 }{ 2^t } \frac{ t! }{ \Big(
\frac{t+a}2 + b + 1 \Big)! \Big( \frac{t-a}2 - b \Big)! } \, .
\end{equation}
Note that $ a,b $ occur only in the combination $ a+2b $.

\end{example}

\begin{example}\label{4ce3}

On $ G_3^\discrete $, we have no `standard' random walk or flow;
rather, we introduce a one-parameter family of abstract stochastic
flows,
\begin{equation}\label{4c3}
\begin{gathered}
\xi_{s,t} = \xi_{s,s+1} \xi_{s+1,s+2} \dots \xi_{t-1,t} \quad
\text{for } s,t \in \Z, \, s \le t \, ; \\
\xi_{t,t+1} \text{ are independent random variables } (t\in\Z) \, ; \\
\Pr{ \xi_{t,t+1} = f_- } = \frac12, \quad \Pr{ \xi_{t,t+1} = f_+ } =
\frac{1-p}2, \quad \Pr{ \xi_{t,t+1} = f_* } = \frac{p}2 \, ;
\end{gathered}
\end{equation}
$ p \in (0,1) $ is the parameter. The canonical homomorphism $
G_3^\discrete \to G_2^\discrete $ glues together $ f_+ $ and $ f_* $,
thus eliminating the parameter $ p $ and giving the standard abstract
flow on $ G_2^\discrete $. Defining $ a(\cdot,\cdot), b(\cdot,\cdot),
c(\cdot,\cdot) $ by
\[
\xi_{s,t} = f_{a(s,t),b(s,t),c(s,t)}
\]
we see that the joint distribution of $ a(\cdot,\cdot) $ and $
b(\cdot,\cdot) $ is the same as before.

Representation \eqref{4b2} of $ G_3^\discrete $ turns the abstract
flow into a stochastic flow on $ \Z_+ $. Its single-point
motion is a sticky random walk,\index{sticky random walk}
\[
t \mapsto \xi_{0,t} (0) = c(0,t) \, .
\]

In order to find the conditional distribution of $ c(\cdot,\cdot) $
given $ a(\cdot,\cdot) $ and $ b(\cdot,\cdot) $ we observe that
\begin{gather}
a(0,t) - c(0,t) = \min \( a(0,t), \min \{ x : \xi_{\si(x),\si(x)+1} =
 f_* \} \) \label{4c4c} \\
\quad \text{where } \si(x) = \max \{ s = 0,\dots,t : a(0,s) = x \} \,
 , \quad -b(0,t) \le x < a(0,t) . \notag \\
\begin{gathered}\includegraphics{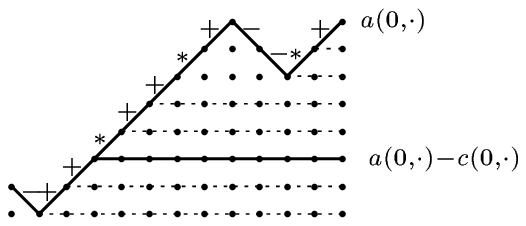}\end{gathered}
\notag
\end{gather}
Therefore the conditional distribution of $ c(0,t) $ is basically the
truncated geometric distribution. More exactly, it is the
(conditional) distribution of
\begin{equation}\label{4c4d}
\max \( 0, a(0,t) + b(0,t) - G + 1 \) \, , \qquad G \sim \Geom (p) \,
;
\end{equation}
here $ G $ is a random variable, independent of $ a(\cdot,\cdot),
b(\cdot,\cdot) $, such that $ \Pr{ G = g } = p(1-p)^{g-1} $ for $ g =
1,2,\dots\, $ This is the discrete counterpart of a well-known result
of J.~Warren \cite{Wa97}. So,
\begin{equation}\label{4c4}
\Pr{ \xi_{0,t} = f_{a,b,c} } = \frac{ a+2b+1 }{ 2^t } \frac{ t! }{ \Big(
\frac{t+a}2 + b + 1 \Big)! \Big( \frac{t-a}2 - b \Big)! } \cdot p
(1-p)^{a+b-c}
\end{equation}
for $ c > 0 $; for $ c = 0 $ the factor $ p(1-p)^{a+b-c} $ turns into
$ (1-p)^{a+b} $, rather than $ p(1-p)^{a+b} $, because of truncation.

\end{example}

\subsection{Three continuous semigroups}
\label{sec:4.4}

The continuous counterpart of the discrete semigroup $ G_1^\discrete =
\Z $ is the semigroup $ G_1 = \R = \{ f_a : a \in \R \} $, $ f_{a_1}
f_{a_2} = f_{a_1+a_2} $.

The continuous counterpart of the discrete semigroup $ G_2^\discrete =
\{ f_{a,b} : a,b \in \Z, \, b \ge 0, \, a+b \ge 0 \} $ is the
semigroup
\begin{equation}\label{4d1}
\begin{gathered}
G_2 = \{ f_{a,b} : a,b \in \R, \, b \ge 0, \, a+b \ge 0 \} \, , \\
f_{a_1,b_1} f_{a_2,b_2} = f_{a,b} \, , \quad
\begin{aligned}
 a &= a_1 + a_2 \, , \\
 b &= \max ( b_1, b_2-a_1 )
\end{aligned}
\end{gathered}
\end{equation}
(recall \eqref{4a2}). The canonical homomorphism $ G_2 \to G_1 $ maps
$ f_{a,b} $ to $ f_a $.

The continuous counterpart of the discrete semigroup $ G_3^\discrete =
\{ f_{a,b,c} : a,b,c \in \Z, \, b \ge 0, \, 0 \le c \le a+b \} $ is
the semigroup
\begin{equation}\label{4d2}
\begin{gathered}
G_3 = \{ f_{a,b,c} : a,b,c \in \R, \, b \ge 0, \, 0 \le c \le a+b \}
 \, , \\
f_{a_1,b_1,c_1} f_{a_2,b_2,c_2} = f_{a,b,c} \, , \quad
\begin{aligned}
 a &= a_1 + a_2 \, , \\
 b &= \max ( b_1, b_2-a_1 ) \, ,
\end{aligned}
\quad
 c = \begin{cases}
  a_2 + c_1 &\text{if $ c_1 > b_2 $},\\
  c_2 &\text{otherwise}
 \end{cases}
\end{gathered}
\end{equation}
(recall \eqref{4a4}). The canonical homomorphism $ G_3 \to G_2 $ maps
$ f_{a,b,c} $ to $ f_{a,b} $.

Note that $ G_1 $ is commutative but $ G_2, G_3 $ are not. Also, $ G_1
$ and $ G_2 $ are topological semigroups, but $ G_3 $ is not (since
the composition is discontinuous at $ c_1 = b_2 $).

There are two one-parameter semigroups in $ G_2 $, $ \{ f_{a,0} : a
\in [0,\infty) \} $ and $ \{ f_{-b,b} : b \in [0,\infty) \} $. They
generate $ G_2 $ according to the relation $ f_{b,0} f_{-b,b} = 1 $;
namely, $ f_{a,b} = f_{-b,b} f_{a+b,0} $.

There are three one-parameter semigroups in $ G_3 $, $ \{ f_{a,0,0} :
a \in [0,\infty) \} $, $ \{ f_{-b,b,0} : b \in [0,\infty) \} $ and $
\{ f_{c,0,c} : c \in [0,\infty) \} $. They generate $ G_3 $ according
to relations $ f_{b,0,0} f_{-b,b,0} = 1 $, $ f_{b,0,b} f_{-b,b,0} = 1
$, and $ f_{c,0,c} f_{a,0,0} = f_{c,0,c} f_{a,0,a} $ for $ c > 0 $;
namely, $ f_{a,b,c} = f_{-b,b,0} f_{a+b-c,0,0} f_{c,0,c} $. 

Here is a faithful representation of $ G_2 $ on $ [0,\infty) $ (recall
\eqref{4b1}):
\begin{equation}
\begin{gathered}
f_{a,b} (x) = a + \max (x,b) \, ,
\end{gathered}
\qquad
\begin{gathered}\includegraphics{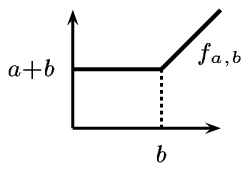}\end{gathered}
\end{equation}
$ x \in [0,\infty) $.

Here is a faithful representation of $ G_3 $ on $ [0,\infty) $ (recall
\eqref{4b2}):
\begin{equation}
\begin{gathered}
f_{a,b,c} (x) = \begin{cases}
 c &\text{for $ 0 \le x \le b $}, \\
 x+a &\text{for $ x > b $}.
\end{cases}
\end{gathered}
\qquad
\begin{gathered}\includegraphics{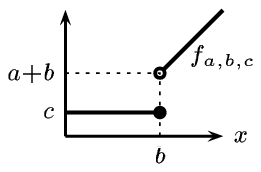}\end{gathered}
\end{equation}
All functions are increasing, but $ f_{a,b} $ are continuous, while $
f_{a,b,c} $ are not.

\subsection{Convolution semigroups in these continuous semigroups}
\label{sec:4.5}

\begin{example}\label{4ee1}

Everyone knows that the binomial distribution \eqref{4c1a} is
asymptotically normal. That is, the distribution of $ \sqrt\eps
a(0,t/\eps) $ converges weakly (for $ \eps\to0 $) to the normal
distribution $ \mu_t^{(1)} = \N(0,t) $. These form a convolution
semigroup, $ \mu_s^{(1)} * \mu_t^{(1)} = \mu_{s+t}^{(1)} $.

Note however, that $ a(s,t) $ and $ \xi_{s,t} $ are defined (see
\eqref{4c1}) only for integers $ s,t $. We may extend them, in one way
or another, to real $ s,t $. Or alternatively, we may use coarse
instants $ t = \( t[i] \)_{i=1}^\infty $, $ t[i] \in \frac1i \Z $, $
t[i] \to t[\infty] $, introduced in \ref{sec:3.2}. For every
coarse instant $ t
$, the distribution of $ i^{-1/2} a(0,it[i]) $ converges weakly (for $
i \to \infty $) to $ \mu_{t[\infty]}^{(1)} = \N(0,t[\infty]) $.

\end{example}

\begin{example}\label{4ee2}

The two-dimensional distribution \eqref{4c2} on $ G_2^\discrete $ has
its asymptotics. Namely, the joint distribution of $ i^{-1/2}
a(0,it[i]) $ and $ i^{-1/2} b(0,it[i]) $ converges weakly (for $ 
i \to \infty $) to the measure $ \mu_{t[\infty]}^{(2)} $ with density
(on the relevant domain $ b>0 $, $ a+b>0 $; $ t $ means $ t[\infty]
$):
\begin{equation}\label{4e1}
\frac{ \mu_t^{(2)} (\D a \D b) }{ \D a \D b } = \frac{ 2 (a+2b) }{
\sqrt{2\pi} \, t^{3/2} } \exp \bigg( - \frac{ (a+2b)^2 }{ 2t } \bigg)
\, .
\end{equation}
Treating $ \mu_t^{(2)} $ (for $ t \in [0,\infty) $) as a measure on $ G_2 $,
we get a convolution semigroup: $ \mu_s^{(2)} * \mu_t^{(2)} =
\mu_{s+t}^{(2)} $. Of course, the convolution is taken according to
the composition \eqref{4d1}.

\end{example}

\begin{example}\label{4ee3}

What about the three-dimensional distribution \eqref{4c4} on $
G_3^\discrete \, $? It has a parameter $ p $. In order to get a
non-degenerate asymptotics, we let $ p $ depend on $ i $, namely,
\[
p = \frac1{\sqrt i} \to 0 \, .
\]
Then the distribution of $ i^{-1/2} G $, where $ G \sim \Geom (p) $
(recall \eqref{4c4d}), converges weakly to the exponential
distribution $ \Exp(1) $, and the joint distribution of $ i^{-1/2}
a(0,it[i]) $, $ i^{-1/2} b(0,it[i]) $ and $ i^{-1/2} c(0,it[i]) $
converges weakly to a measure $ \mu_{t[\infty]}^{(3)} $. The measure
has an absolutely continuous part and a singular part (at $ c=0 $),
and may be described (somewhat indirectly) as the joint distribution
of three random variables $ a $, $ b $ and $ (a+b-\eta)^+ $, where the
pair $ (a,b) $ is distributed $ \mu_t^{(2)} $ (see \eqref{4e1}), $
\eta $ is independent of $ (a,b) $, and $ \eta \sim \Exp(1)
$. Treating $ \mu_t^{(3)} $ (for $ t \in [0,\infty) $) as a measure on
$ G_3 $, we get a convolution semigroup: $ \mu_s^{(3)} * \mu_t^{(3)} =
\mu_{s+t}^{(3)} $, the convolution being taken according to the
composition \eqref{4d2}. No need to check the relation `by hand'; it
follows from its discrete counterpart. The latter follows from the
construction of \ref{sec:4.3} (since random variables $
\xi_{0,1}, \xi_{1,2}, \dots, \xi_{s+t-1,s+t} $ are independent). It
may seem that the limiting procedure does not work, since $ G_3 $ is
not a topological semigroup; the composition \eqref{4d2} is
discontinuous at $ c_1 = b_2 $. However, that is not an obstacle,
since the equality $ c_1 = b_2 $ is of zero probability, as far as
triples $ (a_1,b_1,c_1) $ and $ (a_2,b_2,c_2) $ are independent and
distributed $ \mu_s^{(3)} $, $ \mu_t^{(3)} $, respectively ($ s,t > 0
$). The atom of $ c_1 $ at $ 0 $ does not matter, since $ b_2 $ is
nonatomic. The composition is continuous almost everywhere!

\end{example}

\subsection{Getting dyadic}
\label{sec:4.6}

Our flows in $ G_1^\discrete $ and $ G_2^\discrete $ are dyadic (two
equiprobable possibilities in each step), which cannot be said about $
G_3^\discrete $; here, in each step, we have three possibilities $ f_-,
f_+, f_* $ of probabilities $ 1/2, (1-p)/2, p/2 $. Can a
dyadic model produce
the same asymptotic behavior? Yes, it can, at the expense of using $ i
\in \{ 1, 4, 16, 64, \dots \} $ only (recall \ref{3b6}); and,
of course, the dyadic model is more complicated.\footnote{%
 Maybe, a still more complicated construction can use all $ i $; I do
 not know.}
Instead of the trap at $ 0 $, we design a trap near $ 0 $ as follows:
\begin{gather*}
\begin{gathered}\includegraphics{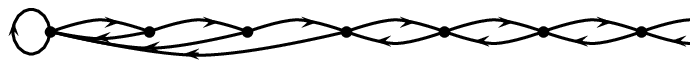}\end{gathered} \\
g_+ = f_* = f_{1,0,1} \, ; \quad g_- = f_-^m f_+^{m-1} = f_{-1,m,0} \,
 ; \\
\Pr{ \xi_{t,t+1} = g_- } = \frac12 = \Pr{ \xi_{t,t+1} = g_+ } \, .
\end{gather*}
The old (small) parameter $ p $ disappears, and a new (large)
parameter $ m $ appears. We'll see that the two models are
asymptotically equivalent, when $ p = 2^{-m} $.

As before, we may denote
\[
\xi_{s,t} = f_{a(s,t),b(s,t),c(s,t)} \, .
\]
Note, however, that only $ a(s,t) $ is the same as before; $ b(s,t) $,
$ c(s,t) $ and $ \xi_{s,t} $ are modified. Formula \eqref{4c1c} for $
b(0,t) $ fails, but still,
\begin{equation}\label{4fb}
b(0,t) = - \min_{s=0,1,\dots,t} a(0,s) + O(m) \, ,
\end{equation}
which is asymptotically the same. Formula \eqref{4c4c} for $ c(0,t) $
also fails. Instead,
\begin{gather}
\begin{gathered}
\includegraphics{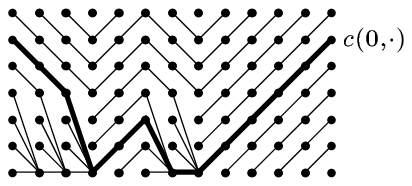}\qquad
\includegraphics{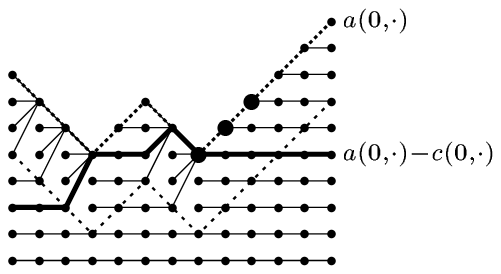}
\end{gathered} \notag \\
a(0,t) - c(0,t) = \min \{ x : \si(x+m-1) - \si(x) = m-1 \} \, ,
 \label{4fd}
\end{gather}
if such $ x $ exists in the set $ \Z \cap [ \min_{[0,t]} a(0,\cdot),
a(0,t)-m+1 ] $; otherwise, $ c(0,t) = O(m) $. (Here $ \si $ is the
same as in \eqref{4c4c}.)

The conditional distribution of $ c(0,t) $, given the path $
a(0,\cdot) $, is not at all geometric (unlike \eqref{4c4d}), since now
$ c(0,t) $ is uniquely determined by $ a(0,\cdot) $. However,
according to \eqref{4fd}, $ a(0,t)-c(0,t) $ is determined by small
increments of the process $ \si(\cdot) $. On the other hand, the
large-scale structure of the path $ a(0,\cdot) $ is correlated mostly
with large increments of $ \si(\cdot) $; small increments are
numerous, but contribute little to the sum. Using this argument, one
can show that $ c(0,t) $ is asymptotically independent of $ a(0,t) $
(and $ b(0,t) $, due to \eqref{4fb}).

The unconditional distribution of $ c(0,t) $ can be found from
\eqref{4fd}, taking into account that increments $ \si(x+1) - \si(x) $
are independent, and each increment is equal to $ 1 $ with probability
$ 1/2 $. We have Bernoulli trials, and we wait for the first block of
$ m-1 $ `successes'. For large $ m $, the waiting time is
approximately exponential, with the mean $ 2^m $.\footnote{%
 Such a block appears, in the mean, after $ 2^{m-1} $ shorter blocks,
 of mean length $ \approx 2 $ each.}
Thus, $ 2^{-m} \( a(0,t) - c(0,t) - \min_{[0,t]} a(0,\cdot) \) $ is
asymptotically $ \Exp(1) $, truncated (at $ c=0 $) as in \ref{sec:4.5}.

Taking the limit $ i = 2^{2m} \to \infty $, we get for $ i^{-1/2}
a(0,it[i]) $, $ i^{-1/2} b(0,it[i]) $, $ i^{-1/2} c(0,it[i]) $ the
limiting distribution $ \mu_{t[\infty]}^{(3)} $, the same as in
\ref{sec:4.5}.

\subsection{Scaling limit}
\label{sec:4.7}

For any coarse instants $ s,t $ such that $ s \le t $, the
distribution $ \mu^{(n)}_{s,t}[i]
$ of $ i^{-1/2} \xi^{(n)}_{is[i],it[i]} $ converges weakly (for $ i
\to \infty $) to the measure $ \mu^{(n)}_{s,t}[\infty] =
\mu^{(n)}_{t[\infty]-s[\infty]} $ on $ G_n $, for our three models, $
n=1,2,3 $. Of course, multiplication of $ \xi $ by $ i^{-1/2} $ is
understood as multiplication of $ a(\cdot,\cdot) $, $ b(\cdot,\cdot)
$, $ c(\cdot,\cdot) $ by $ i^{-1/2} $, which is a homomorphic
embedding of $ G_n^\discrete $ into $ G_n $.

Let $ r, s, t $ be coarse instants, $ r \le s \le t $. Due to
independence, the joint distribution $ \mu^{(n)}_{r,s}[i] \otimes
\mu^{(n)}_{s,t}[i] $ of random variables $ i^{-1/2}
\xi^{(n)}_{ir[i],is[i]} $ and $ i^{-1/2} \xi^{(n)}_{is[i],it[i]} $
converges weakly to $ \mu^{(n)}_{r,s}[\infty] \otimes
\mu^{(n)}_{s,t}[\infty] $. However, we need the joint distribution of
three random variables,
\[
i^{-1/2} \xi^{(n)}_{ir[i],is[i]} \, , \quad i^{-1/2}
\xi^{(n)}_{is[i],it[i]} \, , \quad i^{-1/2} \xi^{(n)}_{ir[i],it[i]} \,
,
\]
the third being the product of the first and the second in the
semigroup $ G_n $. For $ n=1,2 $ weak convergence for the triple
follows immediately from weak convergence for the pair, since the
composition is continuous. For $ n=3 $, discontinuity of the
composition in $ G_3 $ does not invalidate the argument, since the
composition is continuous almost everywhere w.r.t.\ the relevant
measure (recall \ref{sec:4.5}).

Similarly, for every $ k $ and all coarse instants $ t_1 \le \dots
\le t_k $, the joint distribution of $ k(k-1)/2 $ random variables $ 
i^{-1/2} \xi^{(n)}_{it_l[i],it_m[i]} $, $ 1 \le l < m \le k $,
converges weakly (for $ i \to \infty $). We choose a sequence $
(t_k)_{k=1}^\infty $ of coarse instants such that the sequence of
numbers $ (t_k[\infty])_{k=1}^\infty $ is dense in $ \R $, and use
\ref{2c9}, getting a coarse probability space.

The H\"older condition, the same as in \ref{2a1b}, holds for
all three
models. I mean H\"older continuity of $ a(\cdot,\cdot) $, $
b(\cdot,\cdot) $, $ c(\cdot,\cdot) $. Indeed, $ a(\cdot,\cdot) $ is
the same as in \ref{2a1b}; $ b(\cdot,\cdot) $ is related to $
a(\cdot,\cdot) $ via \eqref{4c1c} or \eqref{4fb}, and $ c(\cdot,\cdot)
$ satisfies (on any interval)
\[
\max_{|s-t|\le x} | c(0,s) - c(0,t) | \le \max_{|s-t|\le x} | a(0,s) -
a(0,t) | \, ,
\]
though, for the model of \ref{sec:4.6}, $ O(m) $ must be added.

Thus, a joint $ \si $-compactification is constructed for all three
models (the third model --- in two versions, \ref{4ce3} and
\ref{sec:4.6}).

\subsection{Noises}
\label{sec:4.8}

\begin{example}\label{4h1}

The standard flow in $ G_1^\discrete $, rescaled by $ i^{-1/2} $,
gives us a coarse probability space, identical to that of
\ref{3b3}. It is a dyadic coarse factorization. Its refinement is the
Brownian continuous factorization. Equipped with the natural time
shift, it is a noise.

\end{example}

\begin{example}\label{4h2}

The standard flow in $ G_2^\discrete $, rescaled by $ i^{-1/2} $,
gives us another coarse probability space. It is also a dyadic coarse
factorization (the proof is similar to the previous case). Its
`two-dimensional nature' is a delusion; the dyadic coarse
factorization is identical to that of \ref{4h1}. The second
dimension $ b(\cdot,\cdot) $ reduces to the first dimension, $
a(\cdot,\cdot) $, by \eqref{4c1c}.

\end{example}

\begin{example}\label{4h3}

The flow in $ G_3 $, introduced in \ref{4ce3}, rescaled by $
i^{-1/2} $ with $ p = i^{-1/2} $ (recall \ref{4ee3}), gives us
a coarse probability space. It is not a dyadic coarse factorization,
since it is not dyadic. However, it satisfies a natural generalization
of \ref{3b1} to the non-dyadic case (the proof is as
before). Its refinement is a continuous factorization, and (with
natural time shift), a noise; it may be called the noise of
stickiness.\index{noise of stickiness}

Once again, the second dimension, $ b(\cdot,\cdot) $, reduces to the
first dimension, $ a(\cdot,\cdot) $. Indeed, the joint distribution of
$ a(\cdot,\cdot) $ and $ b(\cdot,\cdot) $ is the same as in
\ref{4h2}. What about the third dimension, $ c(\cdot,\cdot) \, $?

The conditional distribution of $ c(s,t) $, given $ a(s,t) $ and $
b(s,t) $, is basically truncated exponential. Namely, it is the
distribution of $ \( a(s,t) + b(s,t) - \eta \)^+ $ where $ \eta \sim
\Exp(1) $; see \ref{4ee3}. Moreover, for any $ r < s < t $, the
conditional distribution of $ c(r,t) $ given $ a(r,s), b(r,s) $ and $
a(s,t), b(s,t) $, is still the distribution of $ \( a(r,t) + b(r,t) -
\eta \)^+ $. In other words, $ c(r,t) $ is conditionally independent
of $ a(r,s), b(r,s), a(s,t), b(s,t) $, given $ a(r,t), b(r,t) $. That
is a property of the composition \eqref{4d2}; if $ c_1 \sim \( a_1 +
b_1 - \eta_1 \)^+ $ and $ c_2 \sim \( a_2 + b_2 - \eta_2 \)^+ $ then $
c \sim \( a + b - \eta \)^+ $.
\[
\begin{gathered}\includegraphics{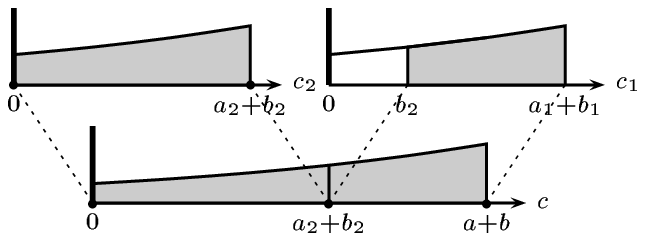}\end{gathered}
\]
It follows by induction that the conditional distribution of $
c(t_1,t_n) $, given all $ a(t_i,t_j) $ and $ b(t_i,t_j) $, is given by
the same formula $ \( a(t_1,t_n) + b(t_1,t_n) - \eta \)^+ $, $ \eta
\sim \Exp(1) $, for every $ n $ and $ t_1 < \dots < t_n $. Therefore,
the same holds for the conditional distribution of $ c(s,t) $ given
all $ a(u,v) $ and $ b(u,v) $ for $ u,v $ such that $ s \le u \le v
\le t $ (a well-known result of J.~Warren \cite{Wa97}). We see that $
c(\cdot,\cdot) $ is not a function of $ a(\cdot,\cdot) $ (and $
b(\cdot,\cdot) $).

\end{example}

\begin{example}

Another flow in $ G_3^\discrete $, introduced in \ref{sec:4.6},
being rescaled by $ i^{-1/2} $ with $ i = 2^{2m} $, gives us a dyadic
coarse factorization. Its refinement is the same continuous
factorization (and noise) as in \ref{4h3}.

\end{example}

\subsection{The Poisson snake}
\label{sec:4.9}

Formula \eqref{4c4c} suggests a description of the sticky flow in $
G_3^\discrete $ by a combination of a simple random walk $
a(\cdot,\cdot) $ and a random subset of the set of its `chords'. A
chord may be defined as an interval $ [s,t] $, $ s,t \in \Z $, $ s < t
$, such that $ a(s,t) = 0 $ and $ a(s,u) > 0 $ for all $ u \in (s,t)
\cap \Z $. Or equivalently, a chord is a horizontal straight segment
on the plane that connects points $ \( s, a(0,s) \) $ and $ \( t,
a(0,t) \) $ and goes below the graph of $ a(0,\cdot) $. The random
subset of chords is very simple: every chord belongs to the subset
with probability $ p $, independently of others. Note that $ p =
i^{-1/2} $ is equal to the vertical pitch (after rescaling $
a(\cdot,\cdot) $ by $ i^{-1/2} $). The scaling limit suggests itself:
a Poisson random subset of the set of all chords of the Brownian
sample path.

\begin{definition}

A \emph{finite chord} of a continuous function $ f : \R \to \R $ is a
set of the form $ [s,t] \times \{x\} \subset \R^2 $ where $ s < t $, $
x = f(s) $ and $ t = \inf \{ u \in (s,\infty) : f(u) \le x \} $. An
\emph{infinite chord} of $ f $ is a set of the form $ [s,\infty)
\times \{x\} \subset \R^2 $ where $ x = f(s) $ and $ f(t) > x $ for
all $ t \in (s,\infty) $. A \emph{chord} of $ f $ is either a finite
chord of $ f $, or an infinite chord of $ f $.

\end{definition}

\[
\begin{gathered}\includegraphics{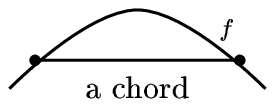}\end{gathered} \qquad
\begin{gathered}\includegraphics{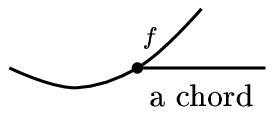}\end{gathered}
\]

If $ f $ decreases, it has no chords. Otherwise it has a continuum of
chords. The set of chords is, naturally, a standard Borel
space,\footnote{%
 For a definition, see \cite[Sect.~12.B]{Ke}.}
due to the one-one correspondence between a chord and its initial
point $ (s,x) \in \R^2 $.

\begin{lemma}\label{4.12}

For every continuous function $ f : \R \to \R $ there exists one and
only one $ \si $\nobreakdash-finite positive Borel measure\footnote{%
 For a definition, see \cite[Sect.~17.A]{Ke}.}
on the space of all
chords of $ f $, such that the set of chords that intersect a vertical
segment $ \{t\} \times [x,y] $ is of measure $ y - x $, whenever $
t,x,y $ are such that $ \inf_{s\in(-\infty,t)} f(s) \le x < y \le f(t)
$.

\end{lemma}

\[
\begin{gathered}\includegraphics{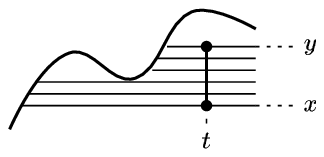}\end{gathered}
\]
The proof is left to the reader. Hint: for every $ \eps > 0 $, the set
of chords longer than $ \eps $ is elementary; on this set, the measure
is locally finite.

The map $ [s,t] \times \{x\} \mapsto s $ (also $ [s,\infty) \times
\{x\} \mapsto s $, of course) sends the measure on the set of chords
(described in \ref{4.12})
into a measure on $ \R $. If $ f $ is of locally finite variation,
then the measure on $ \R $ is just $ (df)^+ $, the positive part of
the Lebesgue-Stieltjes measure. However, we need the opposite case: $
f $ is of infinite variation on every interval, and the measure is
also infinite on every interval. Nevertheless, it is $ \si
$-finite (but not locally finite). We denote it $ (df)^+ $ anyway.

The measure $ (df)^+ $ is concentrated on the set of points of `local
minimum from the right'. If $ f $ is a Brownian sample path then such
points are a set of Lebesgue measure $ 0 $.

So, the set of all chords is a measure space; it carries a natural $
\si $-finite (sometimes, finite) measure. The latter is the intensity
measure of a unique Poisson random measure.\footnote{%
 See for instance \cite[XII.1.18]{RY}.}
This way, (the distribution of) a random set of chords is
well-defined.

Or equivalently, we may consider a Poisson random subset of $ \R $,
whose intensity measure is $ (df)^+ $.

However, it is not so easy to substitute a Brownian sample path $
B(\cdot) $ for $ f(\cdot) $. In order to get a (Poisson) random
variable, we may ask how many random points belong to a given Borel
set $ A \subset \R $ such that $ (dB)^+ (A) < \infty $. Note that for
any \emph{interval} $ A $, $ (dB)^+ (A) = \infty $ a.s. We cannot
choose an appropriate $ A $ without knowing the path $ B(\cdot)
$. The set of all countable dense subsets of $ \R $ does not carry a
natural (non-pathological) Borel structure.

In this aspect, chords are better than points. Chords are
parameterized by three (or two) numbers, and thus, carry a natural
Borel structure, irrespective of $ B(\cdot) $. The random countable
set of chords is not dense; rather, it accumulates toward short
chords.

A point $ (t,x) $ belongs to a random chord of $ B(\cdot) $ if and
only if
\begin{gather*}
x \in \si_t^{-1} (\Pi) \, , \quad \text{that is,} \quad \si_t (x) \in
 \Pi \, , \\
\qquad \text{where } \si_t (x) = \sup \{ s \in (-\infty,t] : B(s) \le
 x \} \text{ for } x \in (-\infty,B(t))
\end{gather*}
(recall \eqref{4c4c}), and $ \Pi $ is the Poisson random subset of $
\R $, whose intensity measure is $ (dB)^+ $. Do not confuse the
inverse image $ \si_t^{-1} (\Pi) $ with the image $ B(\Pi) $. True, $
B(\si_t(x)) = x $, but $ \si_t (B(s)) \ne s $. Sets $ \Pi $ and $
B(\Pi) $ are dense, but the set $ \si_t^{-1} (\Pi) $ is locally
finite. Moreover, $ \si_t^{-1} (\Pi) $ is a Poisson random subset of $
(-\infty,B(t)] $, its intensity being just $ 1 $.

The random countable dense set $ \Pi $ itself is bad; we have no
measurable functions of it. However, the pair $ \( B(\cdot), \Pi \) $
of the Brownian path and the set is good; we have measurable functions
of the pair. In particular, we may use measurable functions of the
locally finite set $ \si_t^{-1} (\Pi) $. Especially,
\[
a(0,t) - c(0,t) = \min \( a(0,t), \min \{ x : \si_t (x) \in \Pi \cap
(0,\infty) \} \) \, .
\]

\enlargethispage*{3mm}

\begin{lemma}

The \sif\ $ \F_{s,t} $ of the noise of stickiness \textup{(}see
\ref{4h3}\textup{)} is generated by Brownian
increments $ B(u) - B(s) $ for $ u \in (s,t) $ and random sets $
\si_u^{-1} \( \Pi \cap (s,t) \) $ for $ u \in (s,t) $ \textup{(}treated as
random variables whose values are finite subsets of $ \R $\textup{).}

\end{lemma}

The proof is left to the reader.

\section{Stability}
\label{sec:5}
\subsection{Discrete case}
\label{sec:5.1}

Fourier-Walsh coefficients, introduced in \ref{sec:3.3} for an
arbitrary dyadic coarse factorization,
\[
f = \sum_{M\in\cC[i]} \hat f_M \tau_M = \hat f_\emptyset +
\sum_{m\in\frac1i\Z} \hat f_{\{m\}} \tau_m +
\sum_{m_1,m_2\in\frac1i\Z,m_1<m_2} \hat f_{\{m_1,m_2\}} \tau_{m_1}
\tau_{m_2} + \dots
\]
help us to examine the stability of a function $ f $, as explained
below. Imagine another array of random signs $
(\tau'_m)_{m\in\frac1i\Z} $ (also independent equiprobable $ \pm1 $)
correlated with the array $ (\tau_m)_{m\in\frac1i\Z} $,
\[
\Ex \tau_m \tau'_m = \rho \quad \text{for each } m \in \frac1i\Z \, ;
\]
$ \rho \in [-1,+1] $ is a parameter. Other correlations vanish. That
is, the joint distribution of all $ \tau_m $ and $ \tau'_m $ is the
product (over $ m \in \frac1i\Z $) of (copies of) such a four-atom
distribution:
\[
\begin{array}{cc@{\hspace{3pt}}@{}lcccr@{}}
 & & & \multicolumn{3}{c}{ \tau_m } & \\
 & & & -1 & & +1 & \\
 \cline{3-7}
 & -1 & \vline & \rule[-2.5mm]{0mm}{7.5mm} \frac{1+\rho}4
  & \vline & \frac{1-\rho}4 & \vline \\
 \cline{3-7}
 \raisebox{5mm}[0pt]{$ \tau'_m $} & +1
  & \vline & \rule[-2.5mm]{0mm}{7.5mm} \frac{1-\rho}4
  & \vline & \frac{1+\rho}4 & \vline \\
 \cline{3-7}
\end{array}
\]
Denoting by $ \ti\Om[i] $ the product of these four-point probability
spaces, we have a natural measure preserving map $ \al : \ti\Om[i] \to
\Om[i] $; as before, $ \Om[i] $ is the product of two-point
probability spaces. In addition, we have another measure preserving
map $ \al' : \ti\Om[i] \to \Om[i] $,
\[
\tau_m \circ \al = \tau_m \, , \quad \tau_m \circ \al' = \tau'_m \, ;
\]
we use the same `$ \tau_m $' for denoting a coordinate function on $
\Om[i] $ and $ \ti\Om[i] $.

For products
\[
\tau_M = \prod_{m\in M} \tau_m \, , \quad M \in \cC[i] \, , \quad
\cC[i] = \{ M \subset \tfrac1i\Z : |M| < \infty \}
\]
we have
\[
\Ex \tau_M \tau'_M = \rho^{|M|} \, , \quad \tau_M \circ \al = \tau_M
\, , \quad \tau_M \circ \al' = \tau'_M \, ,
\]
where $ |M| $\index{zzz@$ "|M"| $, number of elements}
is the number of elements of $ M $. Therefore
\begin{gather*}
\Ex (f\circ\al)(g\circ\al') = \sum_M \rho^{|M|} \hat f_M \hat g_M =
 \ip{g}{ \rho^{\bN[i]} f } \, , \\
\rho^{\bN[i]} : L_2[i] \to L_2[i] \, , \quad \rho^{\bN[i]} \tau_M =
\rho^{|M|} \tau_M \, , \quad \rho^{\bN[i]} f = \sum_M \rho^{|M|} \hat
f_M \tau_M \, .
\end{gather*}
The Hermite operator $ \rho^{\bN[i]} $ is a function of a self-adjoint
operator $ \bN[i] $ defined by $ \bN[i] \tau_M = |M| \tau_M $ for $ M
\in \cC[i] $.

Every bounded function $ \phi : \cC[i] \to \R $ acts on $ L_2[i] $ by
the operator $ f \mapsto \sum_{M\in\cC[i]} \phi(M) \hat f_M
\tau_M $. A commutative operator algebra is isomorphic to the algebra
of functions. The operator $ \rho^{\bN[i]} $ corresponds to the
function $ M \mapsto \rho^{|M|} $. (In some sense, the unbounded
operator $ \bN $ corresponds to the unbounded function $ M \mapsto |M|
$.)

A function $ \phi : \cC[i] \to \{0,1\} $, the indicator of a subset of
$ \cC[i] $, corresponds to a projection operator. Say, for the
(indicator of) the set $ \{ \emptyset \} $, the
operator projects to the one-dimensional space of constants (the
expectation). For the set $ \{ M : M \subset (0,\infty) \}
$, the operator is the conditional expectation, $ \cE{
\cdot }{ \F_{0,\infty}[i] } $.

The function $ M \mapsto |M| $ is the sum (over $ m \in \frac1i\Z $)
of localized functions $ M \mapsto | M \cap \{m\} | $. The latter is
the indicator of the set $ \{ M : M \ni m \} $, corresponding to the
projection operator $ 1 - \cE{ \cdot }{ \F_{\frac1i\Z \setminus \{m\}
} } $. Thus,
\[
\bN f = \sum_m \( f - \cE{ f }{ \F_{\frac1i\Z \setminus \{m\} } } \, .
\]

The operator $ \rho^{\bN[i]} $ may be interpreted as the conditional
expectation w.r.t.\ the sub-\sif\ $ \al^{-1} (\F) $ generated by $
\tau_m \circ \al $, $ m \in \frac1i\Z $:
\[
\cE{ f \circ \al' }{ \al^{-1}(\F) } = ( \rho^{\bN[i]} f ) \circ \al
\quad \text{for } f \in L_2[i] \, .
\]
We may imagine that our data $ \tau_m $ are an unreliable copy of the
true data $ \tau'_m $; each sign $ \tau_m $ is either correct (with
probability $ (1+\rho)/2 $) or inverted (with probability $ (1-\rho)/2
$). If $ \rho $ is close to $ 1 $, our knowledge of $ \tau'_M $ is
satisfactory for moderate $ |M| $ (when $ \rho^{|M|} \approx 1 $) but
very bad for large $ |M| $ (when $ \rho^{|M|} \approx 0 $). The
position of a given function $ f $ between the two extremes is
indicated by the number $ \| f - \rho^\bN f \| $.

\begin{example}\label{5ae1}

In the Brownian coarse factorization (recall \ref{3b3}),
\[
\sup_i \| f[i] - \rho^{\bN[i]} f[i] \| \to 0 \quad \text{for } \rho
\to 1
\]
for all $ f \in L_2 (\A) $. This follows easily from convergence of
operators (recall \ref{sec:2.3} and \ref{3d2}):
\begin{gather*}
\Lim_{i\to\infty} \rho^{\bN[i]} = \rho^{\bN[\infty]} \, , \\
\rho^{\bN[\infty]} f = \sum_{n=0}^\infty \rho^n
\idotsint\limits_{t_1<\dots<t_n} \hat f ( \{ t_1,\dots,t_n \} ) \,
\D B(t_1)\dots \D B(t_n) \, .
\end{gather*}
Convergence of operators follows from \eqref{2a3}. The same holds for
\ref{3b4}.

\end{example}

\begin{example}\label{5ae2}

A very different situation appears in \ref{3b5}. The second
Brownian motion $ B_2 $ (or rather, its discrete approximation) is
not linear but quadratic in random signs $ \tau_m $, $ m \in \frac1i\Z
$. It is two times less stable:
\[
\bN[i] f_{s,t}^{(2)} [i] = 2 f_{s,t}^{(2)} [i] \, ; \qquad
\Lim_{i\to\infty} \rho^{\bN[i]} = \rho^{2\bN[\infty]} \, ,
\]
if $ \bN[\infty] $ is defined in the same way as in \ref{5ae1}. For $
B_3 $ it is $ \rho^{3\bN[\infty]} $, and so on. Still, $ \sup_i \|
f[i] - \rho^{\bN[i]} f[i] \| \to 0 $ for $ \rho \to 1 $. For $ B_\la
$, however, the change is dramatic. Namely,
\[
\bN[i] f_{s,t}^{(\la)} [i] = \entier (\la\sqrt i) f_{s,t}^{(\la)} [i]
 \, ; \qquad
\Lim_{i\to\infty} \rho^{\bN[i]} = 0^{\bN[\infty]}
\]
for all $ \rho \in (-1,+1) $; here $ 0^{\bN[\infty]} = \lim_{\rho\to0}
\rho^{\bN[\infty]} $ is the orthogonal projection to the
one-dimensional subspace of constants (just the expectation). The same
holds for \ref{3b6}.

\end{example}

Notions of stability and sensitivity are introduced in
\cite[Sects.~1.1, 1.4]{BKS} for a sequence of two-valued functions of
$ 1,2,3,\dots $ two-valued variables. For arbitrary (not just
two-valued) functions, a number of equivalent definitions can be found
in \cite[Sect.~1]{ST}. They may be adapted to our framework as
follows. We consider a function $ f : \Om [\all] \to \R $ such that $
0 < \liminf_i \| f[i] \| \le \limsup_i \| f[i] \| < \infty $. We say
that $ f $ is
stable,\index{stable (discrete case)}
if $ \sup_i \| f[i] - \rho^{\bN[i]} f[i] \| \to
0 $ when $ \rho \to 1 $. We say that $ f $ is
sensitive,\index{sensitive (discrete case)}
if $ \|
\rho^{\bN[i]} f[i] - 0^{\bN[i]} f[i] \| \to 0 $ when $ i \to \infty $,
for some (therefore, every) $ \rho \in (0,1) $. These definitions
conform to \cite{ST} when $ f[i] $ depends only on $ i $ signs $
\tau_{1/i}, \dots, \tau_{i/i} $. In terms of the two $ \rho
$-correlated arrays $ (\tau_m) $, $ (\tau'_m) $, stability means that
$ \Ex \( (f[i]\circ\al')(f[i]\circ\al) \) \to \| f[i] \|^2 $ for $
\rho \to 1 $, uniformly in $ i $. Or, equivalently, $ \Ex \( \cVar{
f[i]\circ\al' }{ \al^{-1}(\F) } \) \to 0 $ when $ \rho \to 1 $,
uniformly in $ i $. Sensitivity means that $ \Ex \( (f[i]\circ\al')
(f[i]\circ\al) \) \to \( \Ex f[i] \)^2 $ when $ n \to \infty $, for
some (therefore, every) $ \rho \in (0,1) $. Or, equivalently, $ \Ex
\big| \cE{ f[i]\circ\al' }{ \al^{-1}(\F) } - \Ex f[i] \big|^2 \to 0 $
when $ n \to \infty $, for some (therefore, every) $ \rho \in (0,1) $.

In particular, those definitions can be applied to any $ f \in L_2(\A)
$ such that $ \| f[\infty] \| \ne 0 $.

Example \ref{5ae1} shows that everything is stable in the Brownian
coarse factorization. In contrast, everything is sensitive in the
coarse factorization generated by $ B_\la $ in \ref{5ae2}. In
\ref{sec:5.3} we will find a reason to rename this `stability'
and `sensitivity' as `micro-stability' and `micro-sensitivity'.

A sufficient condition for sensitivity is found by Benjamini, Kalai
and Schramm in terms of the influence of a (two-valued)
variable on a function, see \cite[Sect.~1.2]{BKS}. In our framework,
the influence\index{influence of variable} of the variable $ \tau_m
$ on a function $ f[i] : \Om[i]
\to \R $ may be defined as the expectation of the square root of the
conditional variance,
\[
\Ex \sqrt{ \cVar{ f[i] }{ \F_{\frac1i\Z\setminus\{m\}} } } \, ;
\]
here $ \F_{\frac1i\Z\setminus\{m\}} $ is the sub-\sif\ of $ \F[i] $
generated by all random signs except for $ \tau_m $. The root of the
conditional variance is simply one half of the difference between two
values of the function $ f[i] $, one value for $ \tau_m = +1 $, the
other for $ \tau_m = -1 $. Thus, our formula gives two times less than
\cite[(1.3)]{BKS}, but the coefficient does not matter. Similarly, for
any set $ M \subset \frac1i \Z $, the influence of $ M $ (that is, of all
variables $ \tau_m $, $ m \in M $) on $ f[i] $ may be defined as
\[
\Ex \sqrt{ \cVar{ f[i] }{ \F_{\frac1i\Z\setminus M} } } \, .
\]
By the way, for a \emph{linear} function, the \emph{squared} influence
is additive (in $ M $); indeed, if $ f[i] = \sum_m c_m \tau_m $, then
$ \cVar{ f[i] }{ \F_{\frac1i\Z\setminus M} } = \Ex \( \sum_{m\in M}
c_m \tau_m \)^2 = \sum_{m\in M} c_m^2 $. The sum of squared influences
appears in the following remarkable result (adapted to our framework).

\begin{theorem}[Benjamini, Kalai, Schramm]\label{thBKS}
Let a function $ f : \Om[\all] \to \{ 0,1 \} $ be such that each $
f[i] $ depends on $ i $ variables $ \tau_{1/i}, \dots, \tau_{i/i} $
only. If
\[
\sum_{k=1}^i \Big( \Ex \sqrt{ \cVar{ f[i] }{
\F_{\frac1i\Z\setminus\{k/i\}} } } \, \Big)^2 \tendsi 0 \, ,
\]
then $ f $ is sensitive.
\end{theorem}

See \cite[Th.~1.3]{BKS}. We will return to the point in \ref{sec:6.4}.

\subsection{Continuous case}
\label{sec:5.2}

We start with the \emph{Brownian} continuous factorization $ \(
(\Om,\F,P), (\F_{s,t})_{s\le t} \) $. Using the Wiener-It\^o
decomposition of $ L_2(\Om,\F,P) $,
\[
f = \sum_{n=0}^\infty \, \underbrace{ \idotsint\limits_{t_1<\dots<t_n}
\hat f ( \{ t_1, \dots, t_n \} ) \, \D B(t_1) \dots \D B(t_n)
}_{\text{belongs to $ n $-th Wiener chaos} } \, , \quad \hat f \in L_2
(\Cfinite) \, ,
\]
we can define a self-adjoint operator $ \bN : L_2 \to L_2 $ such that
for each $ n $, $\bN f = n f $ for all $ f $ of $ n $-th Wiener
chaos. Accordingly, $ \rho^\bN f = \rho^n f $ for these $ f
$. Informally, $ \bN ( \D B(t_1) \dots \D B(t_n) ) = n \D B(t_1) \dots
\D B(t_n) $.

Every bounded Borel function $ \phi $ on $ \Cfinite $
acts on $ L_2 (\Om,\F,P) $ by the operator $ R_\phi $,
\begin{equation}\label{5b1}
R_\phi f = \sum_{n=0}^\infty \, \idotsint\limits_{t_1<\dots<t_n} \phi
\( \{ t_1,\dots,t_n \} \) \hat f \( \{ t_1,\dots,t_n \} \) \, \D
B(t_1)\dots \D B(t_n)  \, .
\end{equation}
The operator $ \rho^\bN $ corresponds to the function $ M \mapsto
\rho^{|M|} $. (In some sense, the unbounded operator $ \bN $
corresponds to the unbounded function $ M \mapsto |M| $.) The
decomposition $
|M| = | M \cap (-\infty,t) | + | M \cap (t,\infty) $ (it holds for $
\mu_f $-almost all $ M $) leads to the operator decomposition $ \bN =
\bN_{-\infty,t} + \bN_{t,\infty} $.
Informally, $ \bN_{-\infty,t} \( \D B(t_1) \dots \D B(t_n) \) = k \D
B(t_1) \dots \D B(t_n) $ and $ \bN_{t,\infty} \( \D B(t_1) \dots \D
B(t_n) \) \linebreak[0]
= (n-k) \D B(t_1) \dots \D B(t_n) $ whenever $ t_1 < \dots < 
t_k < t < t_{k+1} < \dots < t_n $. Accordingly, $ \rho^{\bN} =
\rho^{\bN_{-\infty,t}} \otimes \rho^{\bN_{t,\infty}} $.

A function $ \phi : \Cfinite \to \{0,1\} $, the indicator of a Borel
subset $ \cM $ of $ \Cfinite $, corresponds to the orthogonal
projection operator onto the corresponding (recall Theorem \ref{3.25})
subspace $ H_\cM $. Say, for the (indicator of the) set $ \{ \emptyset
\} $, the operator projects onto the one-dimensional space of constants
(the expectation). For the set $ \{ M : M \subset (0,\infty) \} $ the
operator is the conditional expectation, $ \cE{ \cdot }{ \F_{0,\infty}
} $.

The function
\[
\phi_{s,t} (M) = \begin{cases}
 1 &\text{if $ M \cap (s,t) \ne \emptyset $},\\
 0 &\text{if $ M \cap (s,t) = \emptyset $}
\end{cases}
\]
acts by the operator $ \One - \cE{ \cdot }{
\F_{(-\infty,s)\cup(t,\infty)} } $.

For a finite set $ L = \{ s_1,\dots,s_n \} \subset \R $, $ s_1 < \dots
< s_n $, the function $ \phi_L (M) = \phi_{s_1,s_2} (M) + \dots +
\phi_{s_{n-1},s_n} (M) $ counts intervals $ (s_j,s_{j+1}) $ that
intersect $ M $. Clearly, $ \phi_L (M) \le |M|$, and
\[
\phi_{L_n} (M) \uparrow |M| \quad \text{for $ \mu_f $-almost all $ M $}
\]
if $ L_1 \subset L_2 \subset \dots $ are chosen so that their union is
dense in $ \R $. Accordingly,
\begin{gather}
\bN_{L_n} \uparrow \bN \, , \notag \\
\label{5b2}
\bN_{ \{s_1,\dots,s_n\} } = \sum_{j=1}^{n-1} \( \One - \cE{ \cdot }{
 \F_{(-\infty,s_j)\cup(s_{j+1},\infty)} } \) \, .
\end{gather}
The operator $ \bN $ is thus expressed in terms of the factorization
only, irrespective of the Wiener-It\^o decomposition, which gives us a
bridge to \emph{arbitrary} continuous factorizations. Operators $
R_\phi $ described in the next lemma generalize \eqref{5b1}.

\begin{lemma}\label{5.3}
For every continuous factorization $ \( (\Om,\F,P), (\F_{s,t})_{s\le
t} \) $ there exists one and only one map $ \phi \mapsto
R_\phi $\index{zzrp@$ R_\phi $, operator on $ L_2 $}
from the set of all bounded Borel functions $ \phi : \cC \to \R $ to
the set of (bounded linear) operators on $ L_2 (\Om,\F,P) $ such that

\textup{(a)} the map is a homomorphism of algebras; that is, $
R_{a\phi} = a R_\phi $, $ R_{\phi+\psi} = R_\phi + R_\psi $, $ R_{\phi
\psi} = R_\phi R_\psi $;

\textup{(b)} $ \| R_\phi \| \le \sup_{M\in\cC} |\phi(M)| $;

\textup{(c)} $ R_{\One_\cM} = \Proj_{H_\cM} $ for every Borel set $ \cM
\subset \cC $; here $ \One_\cM $ is the indicator of $ \cM $, and $
(H_\cM) $ is the orthogonal decomposition provided by Theorem
\ref{3.25}.

The map also satisfies the condition

\textup{(d)} let $ \phi, \phi_1, \phi_2, \dots : \cC \to [0,1] $ be
Borel functions such that $ \phi_k \to \phi $ pointwise \textup{(}that
is, $ \phi_k (M) \tendsk \phi(M) $ for each $ M \in \cC $\textup{)};
then $ R_{\phi_k} \to R_\phi $ strongly \textup{(}that is, $ \|
R_{\phi_k} x - R_\phi x \| \tendsk 0 $ for every $ x \in L_2
(\Om,\F,P) $\textup{).}
\end{lemma}

\begin{proof}
Uniqueness and existence are easy: Condition (c) and linearity
determine the map on the algebra of Borel functions $ \phi : \cC \to
\R $ having finite sets of values; it remains to extend the map by
continuity.

For proving Condition (d) we note the equality
\[
\ip{ R_\phi x }{ x } = \int \phi \, d\mu_x \, ,
\]
where $ \mu_x $ is the spectral measure of $ x $; it holds for $ \phi
$ having finite sets of values, and therefore, for all $ \phi $. The
bounded convergence theorem gives us not only $ \ip{ R_{\phi_k} x }{ x
} \to \ip{ R_\phi x }{ x } $, but also $ \ip{ R_{(\phi_k-\phi)^2} x }{
x } \to 0 $. However, $ \| R_{\phi_k} x - R_\phi x \|^2 = \ip{
R_{\phi_k-\phi} x }{ R_{\phi_k-\phi} x } = \ip{ R_{(\phi_k-\phi)^2} x
}{ x } $.
\qqed\end{proof}

\begin{lemma}\label{5b3}

For every continuous factorization $ \( (\Om,\F,P), (\F_{s,t})_{s\le
t} \) $, all finite sets $ L_1 \subset L_2 \subset \dots $ whose
union is dense in $ \R $, and every $ \la \in [0,\infty) $, the limit
\[
U_\la = \lim_n \exp ( -\la \bN_{L_n} ) \, ,
\]
where $ \bN_L $ is defined by \eqref{5b2}, exists in the strong
operator topology, and does not depend on the choice of $ L_1, L_2,
\dots $ Also,
\[
U_\la U_\mu = U_{\la+\mu} \quad \text{for all } \la,\mu \in [0,\infty)
\, .
\]

\end{lemma}

\begin{proof}

We have $ \phi_L = \sum \phi_{s_k,s_{k+1}} $ and $ R_{\phi_{s,t}} =
\One - \cE{\cdot}{\F_{(-\infty,s)\cup(t,\infty)}} $; thus $ R_{\phi_L}
= \bN_L $. It follows that $ R_{\exp(-\la\phi_L)} = \exp (-\la \bN_L)
$. However, $ \exp(-\la\phi_{L_n}) \to \phi_\la $, where $ \phi_\la
(M) = \exp (-\la|M|) $ (and $ e^{-\infty} = 0 $, of course). By
\ref{5.3}(d), $ \exp (-\la \bN_{L_n}) \to R_{\phi_\la} = U_\la $. The
semigroup relation $ U_\la U_\mu = U_{\la+\mu} $ for operators follows
from the corresponding relation $ \phi_\la \phi_\mu = \phi_{\la+\mu} $
for functions.

\qqed\end{proof}

In the Brownian factorization we know that $ U_\la = \exp (-\la \bN) $, $
\bN = \lim_n \bN_{L_n} $. In general, however, the semigroup $
(U_\la)_{\la\ge0} $ is discontinuous at $ \la=0 $ (and $ \bN $ is
ill-defined).

\begin{definition}
Let $ \( (\Om,\F,P), (\F_{s,t})_{s\le t} \) $ be a continuous
factorization, and $ f \in L_2 (\Om,\F,P) $.

\textup{(a)} $ f $ is called
\emph{stable,}\index{stable (continuous case)}
if $ \| f - U_\la f \| \to
0 $ for $ \la \to 0 $, or equivalently, if $ \mu_f $ is concentrated
on $ \Cfinite = \{ M \in \cC : |M| < \infty \} $.

\textup{(b)} $ f $ is called
\emph{sensitive,}\index{sensitive (continuous case)}
if $ U_\la f = 0 $ for
all $ \la > 0 $, or equivalently, if $ \mu_f $ is concentrated on $
\cC \setminus \Cfinite = \{ M \in \cC : |M| = \infty \} $.
\end{definition}

Of course, $ U_0 f = f $ anyway.
For proving equivalence, apply \ref{5.3}(d) to $ U_\la =
R_{\phi_\la} $, $ \phi_\la (M) = \E^{-\la |M|} $.

The space $ L_2(\Om,\F,P) $ decomposes into the direct sum of two
subspaces, stable and sensitive, according to the
decomposition of $ \cC $ into the union of two disjoint subsets, $
\Cfinite $ and $ \cC \setminus \Cfinite $.

A continuous factorization is called
\emph{classical}\index{classical factorization}
(or \emph{stable}), if the stable subspace is the whole $ L_2
(\Om,\F,P) $.

A noise is called
classical,\index{classical noise}
if its continuous factorization is classical.

In order to understand probabilistic meaning of $ U_\la $, consider
first $ \rho^{\bN_L} $, $ L = \{ s_1,\dots,s_n \} $, $ s_1 < \dots <
s_n $. We have
\[
\Om = \Om_{-\infty,s_1} \times \Om_{s_1,s_2} \times \dots \times
\Om_{s_{n-1},s_n} \times \Om_{s_n,\infty}
\]
or rather, $ (\Om,\F,P) =
(\Om_{-\infty,s_1},\F_{-\infty,s_1},P_{-\infty,s_1}) \times \dots $,
but let me use the shorter notation. Each $ \om \in \Om $ may be
thought of as a sequence $ ( \om_{-\infty,s_1}, \om_{s_1,s_2},
\linebreak[0]
\dots \om_{s_{n-1},s_n}, \om_{s_n,\infty} ) $ of local portions of
data. Imagine another portion of data $ \om'_{s_1,s_2} \in
\Om_{s_1,s_2} $, either equal
to $ \om_{s_1,s_2} $ (with probability $ \rho $), or independent of it
(with probability $ 1-\rho $). The joint distribution of $
\om_{s_1,s_2} $ and $ \om'_{s_1,s_2} $ is a convex combination of two
probability measures on $ \ti \Om_{s_1,s_2} = \Om_{s_1,s_2} \times
\Om_{s_1,s_2} $. One measure is concentrated on the diagonal and is
the image of $ P_{s_1,s_2} $ under the map $ \Om_{s_1,s_2} \ni
\om_{s_1,s_2} \mapsto ( \om_{s_1,s_2}, \om_{s_1,s_2} ) \in \ti
\Om_{s_1,s_2} $; this measure occurs with the coefficient $ \rho
$. The other measure is the product measure $ P_{s_1,s_2} \otimes
P_{s_1,s_2} $; it occurs with the coefficient $ 1-\rho $.

Similarly we introduce $ \ti \Om_{s_2,s_3} , \dots, \ti
\Om_{s_{n-1},s_n} $ and construct $ \ti \Om = \Om_{-\infty,s_1} \times
\ti \Om_{s_1,s_2} \times \dots \times \ti \Om_{s_{n-1},s_n} \times
\Om_{s_n,\infty} $ (the factors being equipped with corresponding
measures). It is the same idea as in \ref{sec:5.1}. Again, we
have two measure preserving maps $ \al, \al' : \ti \Om \to \Om $. It
appears that
\[
\cE{ f \circ \al' }{ \al^{-1} (\F) } = ( \rho^{\bN_L} f ) \circ \al
\quad \text{for } f \in L_2 (\Om,\F,P) \, .
\]
This is the probabilistic interpretation of $ \rho^{\bN_L} $; each
portion of data is either correct (with probability $ \rho $), or wrong
(with probability $ 1-\rho $).\footnote{%
 This time, $ \rho \in [0,1] $ rather than $ [-1,1] $. The relation to
 the approach of \ref{sec:5.1} is expressed by the equality
 \begin{multline*}
 \frac{1+\rho}2 \begin{pmatrix} 1/2 & 0 \\ 0 & 1/2 \end{pmatrix} +
  \frac{1-\rho}2 \begin{pmatrix} 0 & 1/2 \\ 1/2 & 0 \end{pmatrix} =
 \begin{pmatrix} (1+\rho)/4 & (1-\rho)/4 \\
   (1-\rho)/4 & (1+\rho)/4 \end{pmatrix}  \\
 = \rho \begin{pmatrix} 1/2 & 0 \\ 0 & 1/2 \end{pmatrix} +
  (1-\rho) \begin{pmatrix} 1/4 & 1/4 \\ 1/4 & 1/4 \end{pmatrix} \, .
 \end{multline*}
 }
However, the portions are not small yet. The limit $ n \to \infty $
makes them infinitesimal, and turns $ \rho^{\bN_L} $ into $ U_\la $,
where $ \rho $ and $ \la $ are related by $ \rho = \E^{-\la} $.

The interpretation above motivates the terms `stable' and `sensitive'.

Constant functions on $ \Om $ are stable; sensitive functions are of
zero mean. This is a terminological deviation from the discrete case;
according to \ref{sec:5.1}, constant functions are both stable
and sensitive.

Two limiting cases of $ U_\la $ are projections. Namely, $ U_\infty =
\lim_{\la\to\infty} U_\la $ is the expectation, and $ U_{0+} =
\lim_{\la\to0+} U_\la $ is the projection onto the stable
subspace. Restricting the `perturbation of local data' to a given interval $
(s,t) $ we get operators $ U_\la^{(s,t)} $. These correspond to
functions $ \cC \ni M \mapsto \exp ( -\la | M \cap (s,t) | ) $ and
satisfy
\begin{equation}
\begin{gathered}
U_\la^{(s,t)} U_\mu^{(s,t)} = U_{\la+\mu}^{(s,t)} \, ; \quad
 U_\la^{(r,s)} U_\la^{(s,t)} = U_\la^{(r,t)} \, ; \\
U_\infty^{(s,t)} = \cE{ \cdot }{ \F_{-\infty,s} \otimes \F_{t,\infty}
 } \, ; \\
U_{0+}^{(s,t)} = \cE{ \cdot }{ \F_{-\infty,s} \otimes \F_{s,t}^\stable
 \otimes \F_{t,\infty} } \, .
\end{gathered}
\end{equation}
Note that \eqref{5b2} may be written as
\begin{equation}\label{5b7}
\bN_{ \{ s_1,\dots,s_n \} } = \( \One - U_\infty^{(s_1,s_2)} \) +
\dots + \( \One - U_\infty^{(s_{n-1},s_n)} \) \, .
\end{equation}

\begin{lemma}\label{5.5}
Let $ \( (\Om,\F,P), (\F_{s,t})_{s\le t} \) $ be a continuous
factorization, $ f \in L_2 (\Om,\F,P) $, and $ g = \eta \circ f $
where $ \eta : \R \to \R $ satisfies $ | \eta(x) - \eta(y) | \le |x-y|
$ for all $ x,y \in \R $. Then
\[
\mu_g ( \cC \setminus \cM_E ) \le \mu_f ( \cC \setminus \cM_E )
\]
for all elementary sets $ E \subset \R $; here $ \cM_E = \{ M \in \cC
: M \subset E \} $.
\end{lemma}

\begin{proof}
We have (up to isomorphism) $ \Om = \Om_E \times \Om_{\R\setminus E} $
(the product of probability spaces is meant). We introduce $ \ti\Om =
\Om \times \Om = (\Om_E \times \Om_E) \times (\Om_{\R\setminus E}
\times \Om_{\R\setminus E}) $ and equip the second factor $
\Om_{\R\setminus E} \times \Om_{\R\setminus E} $ with the product
measure, while the first factor $ \Om_E \times \Om_E $ is equipped
with the measure concentrated on the diagonal, such that (equipping $
\ti\Om $ with the product of these two measures), the measure
preserving `coordinate' maps $ \al, \al' : \ti\Om \to \Om $ satisfy
\begin{gather*}
f \circ \al = f \circ \al' \quad \text{for all $ \F_E $-measurable $ f
 $}, \\
f \circ \al \text{ and } g \circ \al' \text{ are independent, for all
 $ \F_{\R\setminus E} $-measurable $ f,g $}.
\end{gather*}
Then
\[
\cE{ f \circ \al' }{ \al^{-1} (\F) } = \cE{ f }{ \F_E } \circ \al
\quad \text{for all } f \in L_2(\Om,\F,P) \, .
\]
Therefore (recall Theorem \ref{3.25}),
\begin{gather*}
\Ex \( (f\circ\al') (g\circ\al) \) = \Ex \( g \cE{ f }{ \F_E } \) \, ;
 \\
\Ex \( (f\circ\al') (f\circ\al) \) = \ip{ \Proj_{H_{\cM_E}} f }{ f } =
 \mu_f (\cM_E) \, ; \\
\frac12 \Ex ( f\circ\al' - f\circ\al )^2 = \mu_f (\cC) - \mu_f (\cM_E)
= \mu_f ( \cC \setminus \cM_E ) \, .
\end{gather*}
The same holds for $ g $. It remains to note that $ | g\circ\al' -
g\circ\al | = | \eta\circ f\circ\al' - \eta\circ f\circ\al | \le |
f\circ\al' - f\circ\al | $ everywhere on $ \ti\Om $.
\qqed\end{proof}

We introduce a special set $ S $ of Borel functions $ \phi : \cC \to
[0,1] $ in three steps. First, we take all functions of the form $
\One_{\cM_E} $,
\[
\One_{\cM_E} (M) = \begin{cases}
 1 &\text{if $ M \subset E $},\\
 0 &\text{otherwise},
\end{cases}
\]
where $ E \subset \R $ runs over all elementary sets. Second, we
consider all (finite) convex combinations of these $ \One_{\cM_E} $. Third,
we consider the least set $ S $ containing these convex combinations
and closed under pointwise convergence (that is, if $ \phi_k \in S $
and $ \phi_k (M) \to \phi (M) $ for each $ M \in \cC $ then $ \phi \in
S $).

The set $ S $ is convex (since the third step preserves convexity). It
is also closed under multiplication: $ \phi\psi \in S $ for all $
\phi, \psi \in S $. Indeed, multiplicativity holds in the first step,
and is preserved in the second and third steps.

\begin{lemma}\label{5.6}
Let $ \( (\Om,\F,P), (\F_{s,t})_{s\le t} \) $ be a continuous
factorization, $ f \in L_2 (\Om,\F,P) $, and $ g = \eta \circ f $
where $ \eta : \R \to \R $ satisfies $ | \eta(x) - \eta(y) | \le |x-y|
$ for all $ x,y \in \R $. Then
\[
\int (1-\phi) \, d\mu_g \le \int (1-\phi) \, d\mu_f
\]
for all $ \phi \in S $.
\end{lemma}

\begin{proof}
In the first step, for $ \phi = \One_{\cM_E} $, the inequality is
stated by \ref{5.5}. The second step evidently preserves the
inequality. And the third step preserves it due to the bounded
convergence theorem.
\qqed\end{proof}

\begin{lemma}\label{5.7}
Let a Borel set $ \cM \subset \cC $ be such that its indicator
function $ \One_\cM $ belongs to the set $ S $. Then for every
continuous factorization $ \( (\Om,\F,P), (\F_{s,t})_{s\le t} \) $,
the subspace $ H_\cM = \{ f : \mu_f (\cC \setminus \cM) = 0 \} $ of $
L_2(\Om,\F,P) $ is of the form
\[
H_\cM = L_2(\Om,\F_\cM,P)
\]
where $ \F_\cM $ is a sub-\sif\ of $ \F $.
\end{lemma}

\begin{proof}
The subspace satisfies
\[
f \in H_\cM \quad \text{implies} \quad |f| \in H_\cM
\]
(here $ |f|(M) = |f(M)| $ for $ M \in \cC $). Indeed,
\[
\int (1-\One_\cM) \, d\mu_{|f|} \le \int (1-\One_\cM) \, d\mu_{f}
\]
by \ref{5.6}; that is, $ \mu_{|f|} (\cC\setminus\cM) \le
\mu_{f} (\cC\setminus\cM) $. A subspace satisfying such a condition is
necessarily of the form $ L_2(\Om,\F_\cM,P) $.
\qqed\end{proof}

Recall the decomposition of $ L_2 (\Om,\F,P) $ into the sum of two
orthogonal subspaces, stable and sensitive, according to the
decomposition of $ \cC $ into the union of two disjoint subsets, $
\Cfinite $ and $ \cC \setminus \Cfinite $.

\begin{theorem}\label{5.9}
For every continuous factorization $ \( (\Om,\F,P), (\F_{s,t})_{s\le
t} \) $ there exists a sub-\sif\
$ \F_\stable $\index{zzFs@$ \F_\stable $, stable \sif}
of $ \F $ such that for all $ f \in L_2(\Om,\F,P) $
\begin{gather*}
f \text{ is stable if and only if } f \text{ is $ \F_\stable
 $-measurable}; \\
f \text{ is sensitive if and only if } \cE{ f }{ \F_\stable } = 0 \, .
\end{gather*}
\end{theorem}

\begin{proof}
The second statement (about sensitive functions) follows from the
first (about stable functions). By \ref{5.7} it is enough to
prove that the indicator of $ \Cfinite $ belongs to $ S $.

\begin{sloppypar}
For every $ \la \in (0,\infty) $ the function $ \phi_\la : \cC \to
[0,1] $ defined by $ \phi_\la (M) = \exp (-\la |M|) $ belongs to $ S $
due to the limiting procedure $ \phi_\la = \lim \exp (-\la \phi_{L_n})
$ used in the proof of \ref{5b3}. For each $ n $ the function $
\exp (-\la \phi_{L_n}) = \prod \exp ( -\la \phi_{s_k,s_{k+1}} ) $
belongs to $ S $, since each $ \exp ( -\la \phi_{s,t} ) $ is a convex
combination of two indicators, of $ \cM_{(-\infty,s)\cup(t,\infty)} $
and of the whole $ \cM $.
\end{sloppypar}

It remains to note that $ \phi_\la $ converges for $ \la \to 0 $ to
the indicator of $ \Cfinite $.
\qqed\end{proof}

So, a continuous factorization (or a noise) is classical if and only
if $ \F_\stable = \F $.

\subsection{Back to discrete: two kinds of stability}
\label{sec:5.3}

The operator equality $ \Lim \rho^{\bN[i]} = \rho^{\bN[\infty]} $
holds for some dyadic coarse factorizations (recall
\ref{5ae1}) but fails for some others (recall
\ref{5ae2}). Nothing like that happens for spectral measures; $ \mu_f
[i] \to \mu_f [\infty] $ always (see Theorem \ref{3c5} and
\ref{sec:3.4}). However, the operator $ \rho^{\bN[i]} $
corresponds to the function $ \cC[i] \ni M \mapsto \rho^{|M|} $
treated as an element of $ L_\infty ( \mu_f[i] ) $, and the operator $
\rho^{\bN[\infty]} $ corresponds to the function $ \cC[\infty] \ni M
\mapsto \rho^{|M|} $ treated as an element of $ L_\infty (
\mu_f[\infty] ) $. How is it possible?  Where is the origin of the
clash between discrete and continuous?

The origin is discontinuity of functions $ M \mapsto \rho^{|M|} $ and
$ M \mapsto |M| $ w.r.t.\ the Hausdorff topology on $ \cC $.

\begin{example}

Return to the equality $ \bN[i] f^{(2)}_{s,t} [i] = 2 f^{(2)}_{s,t}
[i] $ for $ f^{(2)}_{s,t} [i] = i^{-1/2} \sum \tau_m \tau_{m+(1/i)} $
(see \ref{5ae2} and \ref{3b5}). The spectral measure of $
f^{(2)}_{s,t} [i] $ is concentrated on two-point sets $ M \subset
\frac1i\Z $, namely, on pairs of two adjacent points $ \{ m, m+(1/i)
\} $. However, $ f^{(2)}_{s,t} [\infty] $ is just a Brownian
increment; its spectral measure is concentrated on single-point
sets. Now we see what happens; two close points merge in the limit!
Multiplicity of spectral points eludes the continuous model.

The effect becomes dramatic for $ f^{(\la)}_{s,t} [i] $; everything is
stable in the continuous model ($ i = \infty $), while everything
is sensitive (for $ i \to \infty $) in the discrete model. A finite
spectral set on the continuum hides the infinite multiplicity of each
point.

\end{example}

Conformity between discrete and continuous can be restored by
modifying the idea of stability introduced in
\ref{sec:5.1}. Instead of inverting each $ \tau_m $ (with
probability $ (1-\rho)/2 $)
independently of others, we may invert blocks $ \tau_{s[i]},
\tau_{s[i]+(1/i)}, \dots, \tau_{t[i]} $ where coarse instants $ s,t $
satisfy $ t[\infty] - s[\infty] = \eps $. Each block is inverted with
probability $ (1-\rho)/2 $, independently of other blocks. Ultimately
we let $ \eps \to 0 $, but the order of limits is crucial: $
\lim_{\eps\to0} \lim_{i\to\infty} (\dots) $. This way, we can define
(in discrete time setup)
\emph{block stability}\index{block stable}
and \emph{block sensitivity},\index{block sensitive}
equivalent to stability and sensitivity (resp.) of the
refinement. In contrast, the approach of \ref{sec:5.1} leads to
what may be called
\emph{micro-stability}\index{micro-stable}
and \emph{micro-sensitivity}\index{micro-sensitive}
(for discrete time only).

The function $ \cC \ni M \mapsto \rho^{|M|} $ is not continuous, but
it is upper semicontinuous. Therefore, every micro-stable function is
block stable, and every block sensitive function is micro-sensitive.

\begin{example}

The function $ g_{s,t} $ of \ref{3b6} is micro-sensitive but
block stable. The same holds for all coarse random variables in that
dyadic coarse factorization. It holds also for the second construction
of \ref{3b5} (I mean $ f^{(\la)}_{s,t} $).

\end{example}

\section{Generalizing Wiener Chaos}
\label{sec:6}
\subsection{First chaos, decomposable processes, stability}
\label{sec:6.1}

We consider an arbitrary continuous factorization. As was shown in
Theorem \ref{3.25} and \ref{5.3}, Borel functions $ \phi : \cC
\to \R $ act on $ L_2 (\Om,\F,P) $ by linear operators $ R_\phi $, and
(indicators of) Borel subsets $ \cM \subset \cC $ act by orthogonal
projections to subspaces $ H_\cM $.

In particular, for the Brownian factorization, only $ \Cfinite $ is
relevant. The set $ \{ M \in \Cfinite : |M| = n \} $ corresponds to
the subspace called $ n $-th Wiener chaos.

In general, we may define $ n $-th chaos\index{chaos, $n$-th} as the
subspace of $ L_2(\Om,\F,P) $ that corresponds to $ \{ M \in \cC : |M|
= n \} $. These subspaces are orthogonal, and span the stable subspace
--- not the whole $ L_2(\Om,\F,P) $, unless the noise is classical.

For each $ t \in \R $ the set $ \cM_t = \{ M : M \ni t \} $ is
negligible in the sense that $ H_{\cM_t} = \{ 0 \} $ (recall
\ref{3c7} and \eqref{225}). Neglecting $ \cM_t $ we may treat $ \cC $
as the product,\footnote{%
 Sorry, the formula `$ \cC = \cC_{-\infty,t} \times \cC_{t,\infty} $'
 may be confusing since, on the other hand, $ \cC_{-\infty,t} \subset
 \cC $ and $ \cC_{t,\infty} \subset \cC $. The same can be said about
 the next formula, $ H = H_{-\infty,t} \otimes H_{t,\infty} $.}
\begin{equation}\label{495}
\cC = \cC_{-\infty,t} \times \cC_{t,\infty} \, ,
\end{equation}
where $ \cC_{a,b} $ is the space of all compact subsets of $ (a,b) $;
namely, we treat a set $ M \in \cC $ as the pair of sets $ M \cap
(-\infty,t) $ and $ M \cap (t,\infty) $, assuming $ t \notin M $.

On the other hand, the Hilbert space $ H = H_\cC = L_2 (\Om,\F,P) $
may be treated as the tensor product,
\[
H = H_{-\infty,t} \otimes H_{t,\infty} \, ,
\]
of two Hilbert spaces $ H_{-\infty,t} = H_{\cC_{-\infty,t}} = L_2
(\Om,\F_{-\infty,t},P) $ and $ H_{t,\infty} = H_{\cC_{t,\infty}} = L_2
(\Om,\F_{t,\infty},P) $. Namely, $ f \otimes g $ is just the usual
product $ fg $ of random variables $ f \in L_2 (\Om,\F_{-\infty,t},P)
$ and $ g \in L_2 (\Om,\F_{t,\infty},P) $; note that $ f $ and $ g $
are necessarily independent, therefore $ \Ex |fg|^2 = \( \Ex |f|^2 \)
\( \Ex |g|^2 \) $.

Subspaces $ H_\cM \subset H_{-\infty,t} $ for Borel subsets $ \cM
\subset \cC_{-\infty,t} $ are a $ \si $-additive orthogonal
decomposition of $ H_{-\infty,t} $. The same holds for $ (t,\infty)
$.

\begin{lemma}\label{6.1}
$ H_{\cM_1 \times \cM_2} = H_{\cM_1} \otimes H_{\cM_2} $ for all
Borel sets $ \cM_1 \subset \cC_{-\infty,t} $ and $ \cM_2  \subset
\cC_{t,\infty} $.
\end{lemma}

\begin{proof}
The equality holds for the special case $ \cM_1 = \{ M : M \subset E_1
\} $, $ \cM_2 = \{ M : M \subset E_2 \} $ where $ E_1 \subset
(-\infty,t) $ and $ E_2 \subset (t,\infty) $ are elementary sets;
indeed, $ L_2 (\Om,\F_{E_1},P) \otimes L_2 (\Om,\F_{E_2},P) = L_2
(\Om,\F_{E_1\cup E_2},P) $ since $ \F_{E_1\cup E_2} = \F_{E_1} \otimes
\F_{E_2} $. The general case follows by the monotone class theorem.
\qqed\end{proof}

\begin{theorem}\label{6a1}
The sub-\sif\ generated by the first chaos is equal to $
\F_\stable $.
\end{theorem}

\begin{proof}
The \sif\ is evidently included in $ \F_\stable $.
Given a finite set $ L = \{ s_1,\dots,s_n \} \subset \R $, $ s_1 <
\dots < s_n $, we consider the set $ \cM_L $ of all $ M \in \cC $ such
that $ M \subset (s_1,s_n) $ and each $ [s_k,s_{k+1}] $ contains at
most one point of $ M $. The set $ \cM_L $ being the product (over $ k
$), \ref{6.1} shows that $ H_{\cM_L} $ is the tensor product
(over $ k $) of subspaces of $ L_2 (\Om,\F_{s_k,s_{k+1}},P) $; each
factor is the first chaos on $ (s_k,s_{k+1}) $ plus
constants. Therefore each function of $ H_{\cM_L} $ is measurable
w.r.t.\ the \sif\ generated by the first chaos. We choose $ L_1
\subset L_2 \subset \dots $ whose union is dense in $ \R $; then $
\cM_{L_n} \uparrow \Cfinite $, and corresponding subspaces span the
stable subspace.
\qqed\end{proof}

A random variable $ X \in L_2(\Om,\F,P) $ belongs to the first chaos
if and only if
\[
X = \cE X { \F_{-\infty,t} } + \cE X { \F_{t,\infty} } \quad \text{for
all } t \in \R \, .
\]
For such $ X $, letting $ X_{s,t} = \cE{ X }{ \F_{s,t} } $ we get a
\emph{decomposable process},\index{decomposable process}
that is, a family $ (X_{s,t})_{s\le t} $
of random variables such that $ X_{s,t} $ is $ \F_{s,t} $-measurable
and $ X_{r,s} + X_{s,t} = X_{r,t} $ whenever $ r \le s \le t $. This way
we get decomposable processes satisfying $ \Ex |X_{s,t}|^2 < \infty $
and $ \Ex X_{s,t} = 0 $. Waiving these additional conditions we get a
larger set of processes, but the sub-\sif\ generated by these
processes is still $ \F_\stable $. We may also consider complex-valued
\emph{multiplicative decomposable processes;} it means that $ X_{s,t} :
\Om \to \C $ is $ \F_{s,t} $-measurable and $ X_{r,s} X_{s,t} =
X_{r,t} $. The generated sub-\sif\ is $ \F_\stable $, again. The same
holds under the restriction $ | X_{s,t} | = 1 $ a.s. See
\cite[Th.~1.7]{TV}.

Dealing with a noise (rather than factorization) we may restrict
ourselves to stationary Brownian and Poisson decomposable
processes. `Stationary' means $ X_{r,s} \circ \al_t = X_{r-t,s-t} $.
`Brownian' means $ X_{s,t} \sim \N(0,t-s) $. `Poisson' means $ X_{s,t}
\sim \Poisson (\la(t-s)) $ for some $ \la \in (0,\infty) $. The
generated sub-\sif\ is still $ \F_\stable $. See \cite[Lemma
2.9]{Ts98}. (It was written for the Brownian component, but works also
for the Poisson component.)

For a finite set $ L = \{ s_1,\dots,s_n \} \subset \R $, $ s_1 < \dots
< s_n $, we introduce an operator $ Q_L $ on the space $ L_2^0 = \{ X
\in L_2 (\Om,\F,P) : \Ex X = 0 \} $ by
\[
Q_L = \cE{ \cdot }{ \F_{-\infty,s_1} } + \cE{ \cdot }{ \F_{s_1,s_2} }
+ \dots + \cE{ \cdot }{ \F_{s_{n-1},s_n} } + \cE{ \cdot }{
\F_{s_n,\infty} } \, .
\]

\begin{theorem}\label{6a2}

If finite sets $ L_1 \subset L_2 \subset \dots $ are such that their
union is dense in $ \R $, then operators $ Q_{L_n} $ converge in the
strong operator topology to the orthogonal projection from $ L_2^0 $
onto the first chaos.

\end{theorem}

\begin{proof}
$ Q_L $ is the projection onto $ H_{\cM_L} $, where $ \cM_L $ is the
set of all nonempty $ M \in \cC $ contained in one of the $ n+1 $
intervals. The intersection of subspaces corresponds to the
intersection of subsets.
\qqed\end{proof}

Stochastic analysis gives us another useful tool for calculating the
first chaos, pioneered by Jon Warren \cite[Th.~12]{Wa}. Let $
(B_{s,t})_{s\le t} $ be a decomposable Brownian motion, that is, a
decomposable process such that $ B_{s,t} \sim \N(0,t-s) $. One says
that $ B $ has the \emph{representation property,} if every $ X \in
L_2 (\Om,\F,P) $ such that $ \Ex X = 0 $ is equal to a stochastic
integral,
\[
X = \int_{-\infty}^{+\infty} H(t) \, \D B_{0,t} \, ,
\]
where $ H $ is a predictable \emph{process} w.r.t.\ the filtration $
(\F_{-\infty,t})_{t\in\R} $.

\begin{lemma}

If $ B $ has the representation property then the first chaos is equal
to the set of all \emph{linear} stochastic integrals
\[
\int_{-\infty}^{+\infty} \phi(t) \, \D B_{0,t} \, , \qquad \phi \in
L_2(\R) \, .
\]

\end{lemma}

\begin{proof}
\begin{sloppypar}
Linear stochastic integrals evidently belong to the first chaos. Let $
X $ belong to the first chaos. Consider martingales $ B(t) = B_{0,t}
$, $ X(t) = \cE{ X }{ \F_{-\infty,t} } = \int_{-\infty}^t H(s) \,
\D B(s) $ and their bracket process $ {\ip X B}_t = \int_{-\infty}^t
H(s) \, \D s $. The two-dimensional process $ ( B(\cdot), X(\cdot) ) $
has independent increments; therefore the bracket process has
independent increments as well. On the other hand, the bracket process
is a continuous process of finite variation. Therefore it is
degenerate (non-random), and $ H(\cdot) $ is also non-random.
\end{sloppypar}\qqed
\end{proof}

It follows that $ \F_\stable $ is generated by $ B $.

\begin{example}\label{6a4}

For the noise of stickiness (see Sect.~4), the process $ \( a(s,t)
\)_{s\le t} $ is a decomposable Brownian motion having the
representation property. Therefore it generates $ \F_\stable $. On the
other hand we know (recall \ref{4h3}) that $ a(\cdot,\cdot) $
does not generate the whole \sif. So, the sticky noise is not
classical (Warren \cite{Wa}).

\begin{sloppypar}
The approach of Theorem \ref{6a2} is also applicable. Let $ \phi : G_3
\to [-1,+1] $ be a Borel function, and $ 0 < t-\eps < t < 1 $. We
consider $ \phi(\xi_{0,1}) = \phi ( \xi_{0,t-\eps} \xi_{t-\eps,t}
\xi_{t,1} ) $ (you know, $ \xi_{t-\eps,t} =
f_{a(t-\eps,t),b(t-\eps,t),c(t-\eps,t)} $), and compare it with $
\phi ( \xi_{0,t-\eps} \ti\xi_{t-\eps,t} \xi_{t,1} )
$, where $ \ti\xi_{t-\eps,t} = f_{a(t-\eps,t),b(t-\eps,t),0} $.
\end{sloppypar}
\[
\begin{gathered}\includegraphics{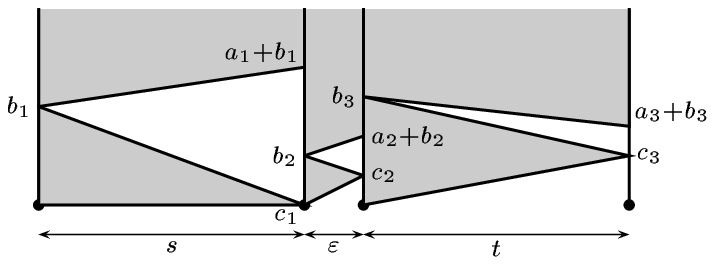}\end{gathered}
\]
It
appears that
\[
\| \phi ( \xi_{0,t-\eps} \xi_{t-\eps,t} \xi_{t,1} ) - \phi (
\xi_{0,t-\eps} \ti\xi_{t-\eps,t} \xi_{t,1} ) \|_{L_2} = O(\eps^{3/4})
= o(\sqrt\eps) \, ,
\]
provided that $ t $ is bounded away from $ 1 $ (otherwise we get $ O (
\eps^{3/4} (1-t)^{-1/2} ) $ with an absolute constant). Taking into
account that $ \ti\xi_{t-\eps,t} $ is measurable w.r.t.\ the \sif\
generated by $ a(\cdot,\cdot) $ we conclude that the projection of $
\phi (\xi_{0,1}) $ onto the first chaos is measurable w.r.t.\ the \sif\
generated by $ a(\cdot,\cdot) $. See \ref{sec:7.2} for the rest.

\end{example}

\subsection{Higher levels of chaos}
\label{sec:6.2}

We still consider an arbitrary continuous factorization. Any Borel
subset $ \cM \subset \cC $ determines a subspace $ H_\cM \subset
L_2 (\Om,\F,P)
$. However, the subset $ \Cfinite \subset \cC $ is special; the
corresponding subspace, being equal to $ L_2 (\F_\stable) $ by Theorem
\ref{5.9}, is of the form $ L_2(\F_1) $ for a sub-\sif\ $ \F_1
\subset \F $.

Another interesting subset is $ \cC_\countable $, the set of all at
most countable compact subsets of $ \R $. It is not a Borel subset of
$ \cC $ \cite[Th.~27.5]{Ke} but still, it is universally measurable
\cite[Th.~21.10]{Ke} (that is, measurable w.r.t.\ every Borel
measure), since its complement is analytic \cite[Th.~27.5]{Ke}. The
Cantor-Bendixson derivative $ M' $ of $ M \in \cC $ is, by
definition, the set of all limit points of $ M $. Clearly, $ M' \in
\cC $, $ M' \subset M $, and $ M' = \emptyset $ if and only if $ M $
is finite. The iterated Cantor-Bendixson derivative $ M^{(\al)} $ is
defined for every ordinal $ \al $ by transfinite recursion: $ M^{(0)}
= M $; $ M^{(\al+1)} = ( M^{(\al)} )' $; and $ M^{(\al)} =
\cap_{\be<\al} M^{(\be)} $ if $ \al $ is a limit ordinal; see
\cite[Sect.~6.C]{Ke}. If $ M \notin \cC_\countable $ then $ M^{(\al)}
\ne \emptyset $ for all $ \al $. If $ M \in \cC_\countable $ then $
M^{(\al)} = \emptyset $ for some finite or countable ordinal $ \al $;
the least $ \al $ such that $ M^{(\al)} = \emptyset $ is called the
Cantor-Bendixson rank of $ M \in \cC_\countable $. It is always of the
form $ \be+1 $, and $ M^{(\be)} $ is a finite set.

Recall the proof of Theorem \ref{5.9}: the indicator of $ \Cfinite
$ belongs to the set $ S $ introduced in \ref{sec:5.2}. Here is
a more general fact.

\begin{lemma}\label{6.5}
Let $ \al $ be an at most countable ordinal, and $ \cM_\al $ the set
of all $ M \in \cC $ such that $ M^{(\al)} = \emptyset $. Then the
indicator function of $ \cM_\al $ belongs to the set $ S $.
\end{lemma}

\begin{proof}
Transfinite induction in $ \al $. For $ \al = 0 $ the claim is
trivial. Let $ \al $ be a limit ordinal. We take $ \al_k \uparrow \al
$, $ \al_k < \al $, and note that $ \cM_\al = \cM_{\al_1} \cup
\cM_{\al_2} \cup \dots $ (indeed, $ M^{(\al_k)} \downarrow M^{(\al)}
$, and $ M^{(\al_k)} $ are compact). Thus, indicators of $ \cM_{\al_k}
$ converge to the indicator of $ \cM_\al $.

The transition from $ \al $ to $ \al + 1 $ needs the following
property of $ S $: for every $ \phi \in S $ and a closed elementary
set $ E $, the function $ M \mapsto \phi ( M \cap E ) $ belongs to $ S
$. Proof: In the first step of constructing $ S $, $ \phi $ is the
indicator of some $ \{ M : M \subset E_1 \} $; thus $ M \mapsto \phi (
M \cap E ) $ is the indicator of $ \{ M : M \subset E_1 \cup
(\R\setminus E) \} $. The second and third steps preserve the
property.

Assume that the indicator function of $ \cM_\al $ belongs to $ S $; we
have to prove the same for $ \al+1 $. The indicator of $ \cM_{\al+1} $
is $ M \mapsto \phi ( M^{(\al)} ) $, where $ \phi $ is the indicator
of $ \Cfinite $. Taking into account that $ \phi \in S $ (see the
proof of Theorem \ref{5.9}), we will prove a more general fact: the
function $
M \mapsto \phi ( M^{(\al)} ) $ belongs to $ S $ for every $ \phi \in S
$ (not just the indicator of $ \Cfinite $). The property is
evidently preserved by the second and third steps of constructing $ S
$; it remains to prove it in the first step. Here $ \phi $ is the
indicator of $ \{ M : M \subset E \} $ for an elementary $ E $. We
have to express the set $ \{ M : M^{(\al)} \subset E \} $ as a limit
of sets of the form $ \{ M : (M\cap E_1)^{(\al)} = \emptyset \} $
where $ E_1 $ is a closed elementary set. The indicator of $ \{ M :
(M\cap E_1)^{(\al)} = \emptyset \} $ belongs to $ S $, since it is $
\One_{\cM_\al} ( M \cap E_1 ) $. We note that, for $ \eps \to 0 $,
\begin{align*}
\{ M : ( M \cap (-\infty,\eps] )^{(\al)} = \emptyset \} &\uparrow \{ M
 : M^{(\al)} \subset (0,\infty) \} \, , \\
\{ M : ( M \cap (-\infty,-\eps] )^{(\al)} = \emptyset \} &\downarrow
\{ M : M^{(\al)} \subset [0,\infty) \} \, ,
\end{align*}
which does the job for two special cases, $ E = (0,\infty) $ and $ E =
[0,\infty) $, and shows how to deal with a boundary point, belonging to
$ E $ or not. The general case is left to the reader.
\qqed\end{proof}

\begin{theorem}\label{6.7}
Let $ \( (\Om,\F,P), (\F_{s,t})_{s\le t} \) $ be a continuous
factorization.

\textup{(a)}
There exists a sub-\sif\ $ \Ec $ of $ \F $ such that for all $ f \in
L_2(\Om,\F,P) $, $ f $ is $ \Ec $-measurable if and only if $ \mu_f $
is concentrated on $ \cC_\countable $.

\textup{(b)}
For every at most countable ordinal $ \al $ there exists a sub-\sif\ $
\Ec_\al $ of $ \F $ such that for all $ f \in L_2(\Om,\F,P) $, $ f $ is
$ \Ec_\al $-measurable if and only if $ \mu_f $ is concentrated on the
set of $ M \in \cC $ such that $ M^{(\al)} = \emptyset $
\textup{(}that is, of Cantor-Bendixson rank less than or equal to $ \al
$\textup{)}.
\end{theorem}

\begin{proof}
Item (a) follows from (b), since $ \Ec_\al = \Ec_{\al+1} $ for
countable $ \al $ large enough (see \cite[Th.~6.9]{Ke}), and $ \mu_f
(\cC_\countable) = \sup_\al \mu_f \{ M : M^{(\al)} = \emptyset \} $
(see \cite{Ke}, the proof of Th.~21.10, and Th.~35.23).

Item (b) follows from \ref{6.5}, \ref{5.7}.
\qqed\end{proof}

Let us concentrate on Item (b) for $ \al = 0,1,2 $. The case $ \al = 0
$ is trivial: only the empty set $ M $, and only constant functions $
f $. The case $ \al = 1 $ was discussed before: finite sets $ M $ and
stable functions $ f $. The case $ \al = 2 $ means that $ M' $ is
finite.

We define the $ n $-th \emph{superchaos}\index{superchaos, $ n $-th}
as the subspace $ H_\cM
\subset L_2 (\Om,\F,P) $ corresponding to $ \{ M \in \cC : |M'| = n \}
$. These subspaces are orthogonal. The $ 0 $-th superchaos is the
stable subspace, while for $ n = 1,2,\dots $ the $ n $-th superchaos
consists of (some) sensitive functions. By Theorem \ref{6.7}(b), the
subspace spanned by $ n $-th superchaos spaces for all $ n =
0,1,2,\dots $ is of the form $ L_2 (\Om,\Ec_2,P) $ where $ \Ec_2 $ is a
sub-\sif\ of $ \F $. Similarly to Theorem \ref{6a1}, the sub-\sif\
generated by the first superchaos and $ \F_\stable $ is equal to $
\Ec_2 $.

Similarly to \eqref{5b2} and \eqref{5b7} we may `count' points of $ M'
$ by the operator
\begin{multline*}
\bN'_{ \{s_1,\dots,s_n\} } = \sum_{j=1}^{n-1} \( \One - \cE{ \cdot }{
 \F_{-\infty,s_j} \otimes \F^\stable_{s_j,s_{j+1}} \otimes
 \F_{s_{j+1},\infty} } \) \\
= \( \One - U_{0+}^{(s_1,s_2)} \) + \dots + \( \One -
 U_{0+}^{(s_{n-1},s_n)} \) \, ,
\end{multline*}
or rather its limit $ \bN' = \lim_n \bN'_{L_n} $. Further, similarly
to \ref{5b3}, we may define
\[
V_\la = \lim_n \exp ( -\la \bN'_{L_n} ) \, .
\]
This way, an ordinal hierarchy of operators may be constructed. It
corresponds to the Cantor-Bendixson hierarchy of countable compact
sets.

Introducing
\begin{multline*}
Q'_{ \{s_1,\dots,s_n\} } X =
 \cE{ X }{ \F_{-\infty,s_1} \otimes \F^\stable_{s_1,\infty} } + 
 \cE{ X }{ \F^\stable_{-\infty,s_1} \otimes \F_{s_1,s_2} \otimes
 \F^\stable_{s_2,\infty} } \\
 + \dots + \cE{ X }{ \F^\stable_{-\infty,s_{n-1}} \otimes \F_{s_{n-1,s_n}} \otimes
 \F^\stable_{s_n,\infty} } + \cE{ X }{ \F^\stable_{-\infty,s_n}
 \otimes \F_{s_n,\infty} }
\end{multline*}
for $ X \in L_2 (\Om,\F,P) $ such that $ \cE X {\F_\stable} = 0 $, we
get such a counterpart of Theorem \ref{6a2}.

\begin{theorem}\label{6b2}

If finite sets $ L_1 \subset L_2 \subset \dots $ are such that their
union is dense in $ \R $, then operators $ Q'_{L_n} $ converge in
the strong operator topology to the orthogonal projection from the
sensitive subspace onto the first superchaos.

\end{theorem}

\begin{proof}
$ Q'_L $ is the projection onto $ H_{\cM_L} $, where $
\cM_L $ is the set of all nonempty $ M \in \cC $ such that $ M' $ is
contained in one of the $ n+1 $ intervals. The intersection of
subspaces corresponds to the intersection of subsets.
\qqed\end{proof}

\begin{example}\label{6b3}

For the sticky noise, consider such a random variable $ X $: the
number of random chords $ [s,t] \times \{ x \} $ such that $ s > 0 $
and $ t > 1 $. In other words (see \ref{sec:4.9}),
\[
X = | \{ x : \si_1(x) \in \Pi \cap (0,\infty) \} | \, .
\]
The conditional distribution of $ X $ given the Brownian path $
B(\cdot) = a(0,\cdot) $ is $ \Poisson(\la) $ with $ \la = a(0,1) +
b(0,1) = B(1) - \min_{[0,1]} B(\cdot) $, which is easy to guess from
the discrete counterpart (see \eqref{4c4d}). That is a generalization
of a claim from \ref{4h3}. In fact, the conditional distribution of the
set $ \{ x : \si_1(x) \in \Pi \cap (0,\infty) \} $, given the Brownian
path, is the Poisson point process of intensity $ 1 $ on $ [-b(0,1),
a(0,1)] $, which is a result of Warren \cite{Wa}. Taking into
account that the \sif\ generated by $ B(\cdot) $ is $ \F_\stable $
(recall \ref{6a4}), we get $ \cE{ X }{ \F_\stable } = a(0,1) + b(0,1)
$. The random variable
\[
Y = X - \cE{ X }{ \F_\stable } = X - a(0,1) - b(0,1)
\]
is sensitive, that is, $ \cE{ Y }{ \F_\stable } = 0 $. I claim that $
Y $ belongs to the first superchaos.

The proof is based on Theorem \ref{6b2}. Given $ 0 < s_1 < \dots < s_n < 1 $,
we have to check that $ Y $ can be decomposed into a sum $ Y_0 + \dots
+ Y_n $ such that each $ Y_j $ is measurable w.r.t.\ $
\F^\stable_{0,s_j} \otimes \F_{s_j,s_{j+1}} \otimes
\F^\stable_{s_{j+1},1} $. Here is the needed decomposition:
\begin{gather*}
X_j = | \{ x : \si_1(x) \in \Pi \cap (s_j,s_{j+1}) \} | \, , \\
Y_j = X_j - \cE{ X_j }{ \F_\stable } \, .
\end{gather*}
We apply a small perturbation on $ (0,s_j) $ and $
(s_{j+1},1) $ but not on $ (s_j,s_{j+1}) $. The set $ \Pi \cap
(s_j,s_{j+1}) $ remains unperturbed. The function $ \si_1 $ is
perturbed, but only a little; being a function of $ B(\cdot) $, it is
stable.

So, $ Y $ belongs to the first superchaos, and $ X $ belongs to the
first superchaos plus $ L_2(\F_\stable) $. It means that $ \mu_X $ is
concentrated on sets $ M $ such that $ |M'| \le 1 $.

The same holds for random variables $ X_u = | \{ x : x \le u , \,
\si_1(x) \in \Pi \cap (0,\infty) \} | $, for any $ u $. They all are
measurable w.r.t.\ the \sif\ generated by the first superchaos and $
\F_\stable $. The random variable $ c(0,1) $ is a (nonlinear!)
function of these $ X_u $ (recall \ref{sec:4.9}). We see that the first
superchaos and $ \F_\stable $ generate the whole \sif\ $ \F $. Every
spectral set (of every random variable) has only a finite number of
limit points.

\end{example}

\begin{example}\label{6b4}

Another nonclassical noise, discovered and investigated by Warren
\cite{Wa1}, see also Watanabe \cite{Wat0}, may be called the noise of
splitting. It is the scaling limit of the model of \ref{1e1};
see also \ref{sec:8.3}. Spectral measures of the most
interesting random variables are described explicitly! A spectral set
contains a single limit point, and two sequences converging to the
point from the left and from the right.

Again, every spectral set (of every random variable) has only a finite
number of limit points.

\end{example}

\begin{question}

We have no example of a noise whose spectral sets $ M $ are at most
countable, and $ M' $ is not always finite. Can it happen at all? Can
it happen for the refinement of a dyadic coarse factorization
satisfying \eqref{3e1}?

\end{question}

Beyond $ \cC_\countable $ it is natural to use the Hausdorff
dimension, $ \dim M $, of compact sets $ M \in \cC $. The set $ S $
used in Theorems \ref{5.9} and \ref{6.7} helps again. First, a general
lemma.

\begin{lemma}\label{6.12}
For every probability measure $ \mu $ on $ \cC $ the function $ \phi :
\cC \to [0,1] $ defined by $ \phi(M) = \mu \{ M_1 \in \cC : M \cap M_1
= \emptyset \} $, belongs to the set $ S $.
\end{lemma}

\begin{proof}
We may restrict ourselves to compact subsets of a bounded interval;
let it be just $ [0,1] $. For any such set $ M $ let $ M^{(n)} $
denote the union of intervals $ [\frac k n, \frac{k+1}n] $ ($ k =
0,\dots,n-1 $) that intersect $ M $. The sequence $
(M^{(n)})_{n=1}^\infty $ decreases and converges to $ M $ (in the
Hausdorff metric). For every $ n $, the function $ \phi_n (M) = \mu \{
M_1 : M \cap M_1^{(n)} = \emptyset \} $ belongs to $ S $, since it is
the convex combination of indicators of $ \{ M : M \subset E \} $ with
coefficients $ \mu \{ M_1 : M_1^{(n)} = [0,1] \setminus E \} $, where
$ E $ runs over $ 2^n $ elementary sets. It remains to note that $
\phi_n (M) \uparrow \phi(M) $, since $ M \cap M_1 = \emptyset $ if and
only if $ M \cap M_1^{(n)} = \emptyset $ for some $ n $.
\qqed\end{proof}

\begin{lemma}\label{6.13}
For every $ \al \in (0,1) $ there exists a function $ \phi \in S $
such that $ \phi(M) = 1 $ for all $ M $ satisfying $ \dim M < \al $,
and $ \phi(M) = 0 $ for all $ M $ satisfying $ \dim M > \al $.
\end{lemma}

\begin{proof}
We may restrict ourselves to the space $ \cC_{0,1} $ of all compact
subsets of $ (0,1) $. There exists a probability measure $ \mu $ on $
\cC_{0,1} $ such that the function $ \phi(M) = \mu \{ M_1 : M_1 \cap M
= \emptyset \} $ satisfies two conditions: $ \phi(M) = 1 $ for all
$ M $ such that $ \dim M < \al $, and $ \phi(M) < 1 $ for all $ M $
such that $ \dim M > \al $. That is a result of J.~Hawkes, see
\cite[Th.~6]{Ha}, \cite[Lemma 5.1]{Pe}. By \ref{6.12}, $ \phi \in S
$. By multiplicativity (of $ S $), also $ \phi^n \in S $ for all $ n
$. The function $ \lim_n \phi^n $ satisfies the required conditions.
\qqed\end{proof}

As a by-product we see that the Hausdorff dimension is a \emph{Borel}
function $ \cC \to \R $. (To this end we use an additional limiting
procedure, as in the proof of Theorem \ref{6.14}.)

\begin{theorem}\label{6.14}
Let $ \( (\Om,\F,P), (\F_{s,t})_{s\le t} \) $ be a continuous
factorization, and $ \al \in (0,1) $ a number. Then there exist
sub-\sif s $ \Ec_{\al-}, \Ec_{\al+} $ of $ \F $ such that for all $ f
\in  L_2(\Om,\F,P) $,

\textup{(a)}
$ f $ is measurable w.r.t.\ $ \Ec_{\al-} $ if and only if $ \mu_f $
is concentrated on the set of $ M \in \cC $ such that $ \dim M < \al
$;

\textup{(b)}
$ f $ is measurable w.r.t.\ $ \Ec_{\al+} $ if and only if $ \mu_f $
is concentrated on the set of $ M \in \cC $ such that $ \dim M \le \al
$.
\end{theorem}

\begin{proof}
We choose $ \al_k \to \al $, apply \ref{6.13} for each $ k $,
consider the limit $ \phi $ of corresponding functions $ \phi_k $, and
use \ref{5.7}. The case $ \al_k < \al $ leads to (a), the case $
\al_k > \al $ leads to (b).
\qqed\end{proof}

A more general notion behind Theorems \ref{5.9}, \ref{6.7} and
\ref{6.14} is an ideal. Recall that a subset $ I $ of $ \cC $ is
called an ideal, if
\[
\begin{gathered}
M_1 \subset M_2, \; M_2 \in I \imply M_1 \in I \, , \\
M_1, M_2 \in I \imply ( M_1 \cup M_2 ) \in I \, .
\end{gathered}
\]
In particular, $ \Cfinite $ and $ \cC_\countable $ are ideals. For
every finite or countable ordinal $ \al $, all $ M \in \cC $ such that
$ M^{(\al)} = \emptyset $ are an ideal. For every $ \al \in (0,1) $,
all $ M \in \cC $ such that $ \dim M < \al $ are an ideal. The same
holds for `$ \dim M \le \al $'. All these ideals are shift-invariant:
\[
\begin{gathered}
M \in I \imply (M+t) \in I \quad \text{for all } t \, , \\
M + t = \{ m + t : m \in M \} \, ,
\end{gathered}
\]
but in general, an ideal need not be shift-invariant. Also, all ideals
mentioned above are Borel subsets of $ \cC $, except for $
\cC_\countable $; the latter is universally measurable, but not
Borel. The following theorem is formulated for Borel ideals, but holds
also for universally measurable ideals. Conditions \ref{6.15} (a,b,c)
parallel \ref{3d1} (a,b,c), which means that sub-\sif s $ \Ec_{s,t} $
form a continuous factorization of the quotient probability space $
(\Om,\F,P) / \Ec $.

\begin{theorem}\label{6.15}
Let $ \( (\Om,\F,P), (\F_{s,t})_{s\le t} \) $ be a continuous
factorization, $ I \subset \cC $ a Borel ideal, $ \Ec \subset \F $ a
sub-\sif, and for every $ f \in L_2 (\Om,\F,P) $, $ f $ be $ \Ec
$-measurable if and only if $ \mu_f $ is concentrated on $ I $. Then
sub-\sif s $ \Ec_{s,t} = \Ec \cap \F_{s,t} $ satisfy the conditions
\begin{gather}
\Ec_{r,t} = \Ec_{r,s} \otimes \Ec_{s,t} \quad  \text{whenever } r \le s
 \le t \, , \tag*{\textup{(a)}} \\
\bigcup_{\eps>0} \Ec_{s+\eps,t-\eps} \text{ generates $ \Ec_{s,t} $
 whenever $ s < t $,} \tag*{\textup{(b)}} \\
\bigcup_{n=1}^\infty \Ec_{-n,n}\text{ generates } \Ec \,
 . \tag*{\textup{(c)}}
\end{gather}
\end{theorem}

\begin{proof}
(a) We introduce Borel subsets $ I_{s,t} = \{ M \in I : M \subset
(s,t) \} $ of $ \cC $ and the corresponding subspaces $ H_{s,t} =
H_{I_{s,t}} $ of $ L_2 (\Om,\F,P) $. The equality $ I_{r,t} = I_{r,s}
\times I_{s,t} $ (treated according to \eqref{495}) follows
easily from the fact that $ I $ is an ideal. Lemma \ref{6.1} (or
rather, its evident generalization) states that $ H_{r,t} = H_{r,s}
\otimes H_{s,t} $. On the other hand,
\[
L_2 (\Ec_{s,t}) = L_2 ( \Ec \cap \F_{s,t} ) = L_2 (\Ec) \cap L_2
(\F_{s,t}) = H_I \cap H_{\cC_{s,t}} = H_{I\cap\cC_{s,t}} = H_{s,t} \,
.
\]
So, $ L_2 (\Ec_{r,t}) = L_2 (\Ec_{r,s}) \otimes L_2 (\Ec_{s,t}) $,
therefore $ \Ec_{r,t} = \Ec_{r,s} \otimes \Ec_{s,t} $.

(c) $ \cup_n I_{-n,n} = I $, therefore $ \cup_n H_{I_{-n,n}} $ is
dense in $ H_I $; that is, $ \cup_n L_2 (\Ec_{-n,n}) $ is dense in $
L_2 (\Ec) $, therefore $ \cup_n \Ec_{-n,n} $ generates $ \Ec $.

(b): similarly to (c).
\qqed\end{proof}

\begin{remark}
If the ideal $ I $ is shift-invariant and the given object is a noise
(not only a factorization), then the sub-factorization $
(\Ec_{s,t}) $ becomes a sub-noise. In particular, every nonclassical
noise has its classical (in other words, stable)
sub-noise.\index{classical sub-noise}
\end{remark}

\begin{question}
Does every Borel ideal correspond to a sub-\sif? (For an arbitrary
continuous factorization, I mean. Though, the question is also open for
noises and shift-invariant ideals.)
\end{question}

\subsection{An old question of Jacob Feldman}
\label{sec:6.3}

Let $ \( (\Om,\F,P), (\F_{s,t})_{s\le t} \) $ be a continuous
factorization. Sub-\sif s $ \F_E $ correspond to elementary sets $ E
\subset \R $ (recall \ref{sec:3.4}) and satisfy
\begin{equation}\label{eq50}
\F_{E_1 \cup E_2} = \F_{E_1} \otimes \F_{E_2} \quad \text{whenever }
E_1 \cap E_2 = \emptyset \, .
\end{equation}
It is natural to ask whether or not the map $ E \mapsto \F_E $ can be
extended to all Borel sets $ E \subset \R $ in such a way that
\eqref{eq50} is still satisfied and in addition,
\begin{equation}\label{eq51}
\F_{E_n} \uparrow \F_E \quad \text{whenever } E_n \uparrow E \, .
\end{equation}
The answer is positive if and only if the given continuous
factorization is classical (Theorem \ref{6.21} below, see also
\cite{Ts99}), which solves a question of Feldman \cite{Fe}.

Note that \eqref{eq51} implies
\begin{equation}\label{eq52}
\F_{E_n} \downarrow \F_E \quad \text{whenever } E_n \downarrow E \, .
\end{equation}
Proof: Let $ E_n \downarrow E $, then $ \F_{\R\setminus E_n} \uparrow
\F_{\R\setminus E} $ by \eqref{eq51}, and so $ \F_{\R\setminus E} $
is independent of $ \cap \F_{E_n} $. If $ \F_E $ is strictly less than
$ \cap \F_{E_n} $, then $ \F_E \otimes \F_{\R\setminus E} $ is
strictly less than $ (\cap \F_{E_n}) \otimes \F_{\R\setminus E} $,
which cannot happen, since $ \F_E \otimes \F_{\R\setminus E} = \F $ by
\eqref{eq50}.

An extension satisfying \eqref{eq51}, \eqref{eq52} is unique (if it
exists) by the monotone class theorem. Therefore an extension (of $
(\F_E) $ to the Borel \sif) satisfying \eqref{eq50}, \eqref{eq51} is
unique (if it exists).

\begin{lemma}\label{5.125}
If the factorization is classical then an extension satisfying
\eqref{eq50}, \eqref{eq51} exists.
\end{lemma}

\begin{proof}
By (slightly generalized) Theorem \ref{6a1}, for every elementary $ E
$, the \sif\ $ \F_E = \F_E^\stable $ is generated by the corresponding
portion $ H_E^{(1)} = L_2(\F_E) \cap H^{(1)} $ of the first chaos $
H^{(1)} $. The space $ H_E^{(1)} $ corresponds (in the sense of
Theorem \ref{3.25}) to the subset $ \cM_E^{(1)} \subset \cC $ of all
single-point subsets of $ E $.

Given an arbitrary Borel set $ E \subset \R $, we define the subset $
\cM_E^{(1)} \subset \cC $ as above (that is, all single-point subsets
of $ E $), consider the corresponding subspace $ H_E^{(1)} \subset
H^{(1)} $, and introduce the sub-\sif\ $ \F_E \subset \F $ generated
by $ H_E^{(1)} $.

Given $ f \in H^{(1)} $, we denote by $ f_E $ the orthogonal
projection of $ f $ to $ H_E^{(1)} $; here $ E $ is an arbitrary Borel
set. If $ E_n \uparrow E $ (or $ E_n \downarrow E $) then $ f_{E_n}
\to f $ in $ L_2 $. If $ E $ is elementary then
\[
\Ex \E^{\I f} = \( \Ex \E^{\I f_E} \) \( \Ex \E^{\I f_{\R\setminus E}} \) 
\]
due to independence. The monotone class theorem extends the equality to
all Borel sets $ E $. We conclude that $ f_E $ and $ f_{\R\setminus E}
$ are independent. Therefore \sif s $ \F_E $ and $ \F_{\R\setminus E} $
are independent for every Borel set $ E $. Taking into account that $
H_{E_1\cup E_2}^{(1)} = H_{E_1}^{(1)} \oplus H_{E_2}^{(1)} $ whenever
$ E_1 \cap E_2 = \emptyset $ we get \eqref{eq50}.

If $ E_n \uparrow E $ then $ H_{E_n}^{(1)} \uparrow H_{E}^{(1)} $,
which ensures \eqref{eq51}.
\qqed\end{proof}

Condition (a) of the next lemma is evidently necessary for the
extension to exist. In more topological language, for every open set
$ G \subset \R $ the corresponding \sif\ $ \F_G $ is naturally defined
by approximation (of $ G $ by elementary sets) from within, while a
closed set is approximated from the outside. The necessary condition, $
\F_G \otimes \F_{\R\setminus G} = \F $, appears to be equivalent to
the following (see \ref{5.13}(b)): the set $ M \cap G $ is compact,
for almost all $ M \in \cC $.

\begin{lemma}\label{5.13}
For all elementary sets $ E_1 \subset E_2 \subset \dots $ the
following two conditions are equivalent:

\textup{(a)} $\displaystyle \Big( \bigvee_n \F_{E_n} \Big) \otimes \Big(
\bigwedge_n \F_{\R\setminus E_n} \Big) = \F \, ; $

\textup{(b)} the set $ \{ M \in \cC : \forall n \> M \cap \( (\cup
E_k) \setminus E_n \) \ne \emptyset \} $ is negligible w.r.t.\ the
spectral measure $ \mu_f $ for every $ f \in L_2 (\Om,\F,P) $.
\end{lemma}

\begin{proof}
Denote $ F_n = \R \setminus E_n $, $ \Ec_n = \F_{E_n} $, $ \F_n =
\F_{\R\setminus E_n} $, $ \Ec_\infty = \vee_n \Ec_n $, $ \F_\infty =
\wedge_n \F_n $. Clearly, $ \Ec_\infty $ and $ \F_\infty $ are
independent, and (a) becomes $ \Ec_\infty \vee \F_\infty = \F $. Denote
also $ \cM_n = \{ M \in \cC : M \subset E_n \} $, $ \cN_n = \{ M \in
\cC : M \subset F_n \} $, $ \cM_\infty = \cup_n \cM_n = \{ M \in \cC :
\exists n \> M \subset E_n \} $, $ \cN_\infty = \cap_n \cN_n = \{ M
\in \cC : M \subset \cap F_n \} $; then $ H_{\cM_n} = L_2 (\Ec_n) $, $
H_{\cN_n} = L_2 (\F_n) $. We have $ \cM_n \uparrow \cM_\infty $ and $
\cN_n \downarrow \cN_\infty $; therefore $ L_2 (\Ec_n) = H_{\cM_n}
\uparrow H_{\cM_\infty} $ and $ L_2 (\F_n) = H_{\cN_n} \downarrow
H_{\cN_\infty} $. On the other hand, $ \Ec_n \uparrow \Ec_\infty $ and $
\F_n \downarrow \F_\infty $; therefore $ L_2 (\Ec_n) \uparrow L_2
(\Ec_\infty) $ and $ L_2 (\F_n) \downarrow L_2 (\F_\infty) $. So,
\[
H_{\cM_\infty} = L_2 (\Ec_\infty) \, , \quad H_{\cN_\infty} = L_2
(\F_\infty) \, .
\]
Denote $ \cM_\infty \vee \cN_\infty = \{ M_1 \cup M_2 : M_1 \in
\cM_\infty, M_2 \in \cN_\infty \} $; the same for $ \cM_1 \vee
\cN_\infty $ etc. We have $ H_{\cM_1 \vee \cN_n} = H_{\cM_1} \otimes
H_{\cN_n} $ and $ \cM_1 \vee \cN_n \downarrow \cM_1 \vee \cN_\infty $;
thus $ H_{\cM_1 \vee \cN_\infty} = H_{\cM_1} \otimes H_{\cN_\infty} $
(note a relation to \ref{6.1}). Similarly, $ H_{\cM_n \vee
\cN_\infty} = H_{\cM_n} \otimes H_{\cN_\infty} $. However, $ \cM_n
\vee \cN_\infty \uparrow \cM_\infty \vee \cN_\infty $, and we get $
H_{\cM_\infty \vee \cN_\infty} = H_{\cM_\infty} \otimes H_{\cN_\infty}
$, that is,
\[
H_{\cM_\infty \vee \cN_\infty} = L_2 (\Ec_\infty) \otimes L_2
(\F_\infty) \, .
\]
Now (a) becomes $ H_{\cM_\infty \vee \cN_\infty} = H $, which means
negligibility of the set $ \cC \setminus ( \cM_\infty \vee \cN_\infty
) = \{ M : \forall n \> M \cap \( (\cup E_k) \setminus E_n \) \ne
\emptyset \} $, that is, (b).
\qqed\end{proof}

Every classical factorization satisfies \ref{5.13}(b), since a finite
set $ M $ cannot intersect $ (\cup E_k) \setminus E_n $ for all $ n $.

\begin{lemma}\label{5.14}
If Condition \textup{\ref{5.13}(b)} is satisfied for every $ (E_n) $
then the factorization is classical.
\end{lemma}

\begin{proof}
Let the factorization be not classical. Then we can choose a sensitive
$ f \in L_2 (\Om,\F,P) $, $ \| f \| = 1 $. Assume for convenience that
$ f \in L_2(\F_{0,1}) $, and consider the spectral measure $ \mu_f $;
$ \mu_f $-almost all $ M $ are infinite subsets of $ (0,1) $. We
choose $ p_1,p_2,\dots \in (0,1) $ such that $ \sum p_k \le 1/3 $ (say,
$ p_k = 2^{-k}/3 $). Integer parameters $ n_1 < n_2 < \dots $ will be
chosen later. We introduce independent random elementary sets $ B_1,
B_2, \dots \subset [0,1] $ as follows:
\[
\Prob \bigg\{ B_k = \Big( \frac{l_1-1}{n_k}, \frac{l_1}{n_k} \Big)
\cup \dots \cup \Big( \frac{l_m-1}{n_k}, \frac{l_m}{n_k} \Big) \bigg\}
= p_k^m (1-p_k)^{n_k-m}
\]
whenever $ 1 \le l_1 < \dots < l_m \le n_k $, $ m \in \{ 0,\dots,n_k
\} $. That is, we have a two-parameter family of independent events, $
\( \frac{l-1}{n_k}, \frac{l}{n_k} \) \subset B_k $, where $ l \in \{
1,\dots,n_k \} $, $ k \in \{ 1,2,\dots \} $. The probability of such
an event is equal to $ p_k $. We define $ E_k = B_1 \cup \dots \cup
B_k $; thus $ E_1 \subset E_2 \subset \dots $ is a (random) increasing
sequence of elementary subsets of $ [0,1] $.

We treat $ M $ as a random compact subset of $ (0,1) $, distributed $
\mu_f $ and independent of $ B_1, B_2, \dots \, $ Let $ \ti P $ be the
corresponding probability measure (in fact, product measure) on the
space $ \ti\Om $ of sequences (of sets) $ (M,B_1,B_2,\dots) $. For
each $ k = 0,1,2,\dots $ we define an event $ A_k $, that is, a
measurable subset of $ \ti\Om $, by the following condition on $
(M,B_1,B_2,\dots) $:
\[
M \setminus E_k \text{ is infinite and does not intersect } B_{k+1} \,
;
\]
of course, $ E_0 = \emptyset $.

We can choose $ n_1, n_2, \dots $ such that $ \sum_k \ti P(A_k) \le
1/3 $. Proof: $ \ti P(A_k) $ is a function of $ n_1, \dots, n_k,
n_{k+1} $ that converges to $ 0 $ when $ n_{k+1} \to \infty $ (while $
n_1,\dots,n_k $ are fixed).

The probability of the event
\[
M \setminus E_k \text{ is infinite for all } k
\]
is no less than $ 1 - \sum p_k \ge 2/3 $. Proof: Each $ M $ has a
limit point (at least one), and the point is covered by (the closure
of) $ B_1 \cup B_2 \cup \dots $ with probability $ \le \sum p_k $.

So, there is a positive probability ($ \ge 1/3 $) to such an event:
\[
\text{for each $ k $, the set } M \setminus E_k \text{ is infinite and
intersects } B_{k+1} \, .
\]
However, the conditional probability, given $ B_1, B_2, \dots $ (but
not $ M $) of the event
\[
\text{for each $ k $, the set } M \setminus E_k \text{ intersects }
B_{k+1}
\]
must vanish according to \ref{5.13}(b).
\qqed\end{proof}

\begin{theorem}\label{6.21}
A continuous factorization is classical if and only if the map $ E
\mapsto \F_E $ can be extended from the algebra of elementary sets to
the Borel \sif, satisfying \eqref{eq50} and \eqref{eq51}.
\end{theorem}

\begin{proof}
If the factorization is classical then the extension exists by
\ref{5.125}. Let the extension exist; then \ref{5.13}(a) is satisfied
for all $ (E_k) $, therefore \ref{5.13}(b) is also satisfied, and the
factorization is classical by \ref{5.14}.
\qqed\end{proof}

\subsection{Black noise}
\label{sec:6.4}

\begin{definition}

A noise is \emph{black,}\index{black noise} if its stable \sif\ $
\F_\stable $ is degenerate. In other words: its first chaos contains
only $ 0 $.

\end{definition}

Why `black'? Well, the white noise is called `white' since its
spectral density is constant. It excites harmonic oscillators of all
frequencies to the same extent. For a black noise, however, the
response of any linear sensor is zero!

What could be a physically reasonable nonlinear sensor able to sense a
black noise? Maybe a fluid could do it, which is hinted at by the
following words of Shnirelman \cite[p.~1263]{Shn} about the paradoxical
motion of an ideal incompressible fluid: `\dots\ very strong external
forces are present, but they are infinitely fast oscillating in space
and therefore are indistinguishable from zero in the sense of
distributions. The smooth test functions are not ``sensitive'' enough
to ``feel'' these forces.'

The very idea of black noises, nonclassical factorizations, etc.\ was
suggested to me by Anatoly Vershik in 1994.

\begin{lemma}\label{6.19}
Let $ \( (\Om,\F,P), (\F_{s,t})_{s\le t} \) $ be a continuous
factorization, $ a < b $, $ \cM $ a Borel subset of $ \cC_{a,b} = \{
M \in \cC : M \subset (a,b) \} $, and $ \ti \cM = \{ M \in \cC : M \cap
(a,b) \in \cM \} $. If $ \mu_f (\cM) = 0 $ for all $ f \in L_2 (\Om,\F,P)
$ then $ \mu_f (\ti \cM) = 0 $ for all $ f \in L_2 (\Om,\F,P) $.
\end{lemma}

\begin{proof}
I prove it for $ (a,b) = (0,\infty) $, leaving the general case to the
reader. We have $ \cC = \cC_{-\infty,0} \times \cC_{0,\infty} $, $ \cM
\subset \cC_{0,\infty} $ and $ \ti \cM = \cC_{-\infty,0} \times \cM $ (in
the sense of \eqref{495}). By \ref{6.1}, $
H_{\ti \cM} =  H_{\cC_{-\infty,0} \times \cM} = H_{\cC_{-\infty,0}}
\otimes H_\cM $. By \eqref{326}, the space $ H_\cM $ is trivial (that is,
$ \{ 0 \} $). Therefore $ H_{\ti \cM} $ is also trivial; it remains to
use \eqref{326} again.
\qqed\end{proof}

Recall that a compact set $ M $ is called perfect, if it has no
isolated points. (The empty set is also perfect.) The set $
\cC_\perfect $ of all perfect compact subsets of $ \R $ is a Borel
set in $ \cC $, see \cite[proof of Th.~27.5]{Ke}.

\begin{theorem}\label{6.20}
For every continuous factorization $ \( (\Om,\F,P), (\F_{s,t})_{s\le
t} \) $ the following two conditions are equivalent:

\textup{(a)}
the first chaos space is trivial \textup{(}contains only $ 0
$\textup{)};

\textup{(b)}
for every $ f \in L_2 (\Om,\F,P) $ the spectral measure $ \mu_f $ is
concentrated on $ \cC_\perfect $.
\end{theorem}

\begin{proof}
(b) implies (a) evidently (a single-point set cannot be
perfect). Assume (a). Applying \ref{6.19} to the set $ \cM $ of
all single-point subsets of $ (a,b) $ we see that $ \mu_f $-almost all
$ M \in \cC $ are such that $ M \cap (a,b) $ is not a single-point
set, for all rational $ a < b $. It means that $ M $ is perfect.
\qqed\end{proof}

So, a noise is black if and only if spectral measures are concentrated
on (the set of all) perfect sets.

Existence of black noises was proven first by Tsirelson and Vershik
\cite[Sect.~5]{TV}. A simpler and more natural example is described in
the next section. Another example is found by Watanabe \cite{Wat}.

If all spectral sets are finite or countable (as in
\ref{6b3}, \ref{6b4}), such a noise cannot contain a black sub-noise.

\begin{question}

If a noise contains no black sub-noise, does it follow that all
spectral sets are at most countable?

\end{question}

Perfect sets may be classified, say, by Hausdorff dimension. For any $
\al \in (0,1) $, sets $ M \in \cC $ of Hausdorff dimension $ \le \al $
are a shift invariant ideal, corresponding to a sub-noise. Also, all $
M \in \cC $ of Hausdorff dimension $ \al $ correspond to a `chaos
subspace number $ \al $'. A continuum of such chaos subspaces (not in
a single noise, of course) could occur, describing different `levels
of sensitivity'. For now, however, I know of perfect spectral sets of
Hausdorff dimension $ 1/2 $ only.

\begin{question}\label{6c4}

Can a noise have perfect spectral sets of Hausdorff dimension other
than $ 1/2 \, $? (See also the end of \ref{sec:8.3}.)

\end{question}

\begin{question}\label{6.23}

Can a black noise emerge as the refinement of a dyadic coarse
factorization satisfying \eqref{3e1}?

\end{question}

The following results (especially \ref{6.31}) may be treated
as continuous-time counterparts of Theorem \ref{thBKS} (of Benjamini,
Kalai and Schramm). Given a continuous factorization $ \( (\Om,\F,P),
(\F_{s,t})_{s\le t} \) $ and a function $ f \in L_2 (\Om,\F,P) $, we
define
\[
\bH (f) = \limsup_{\{t_1,\dots,t_n\}\uparrow\,} \sum_{k=1}^{n+1} \Big(
\Ex \sqrt{ \cVar{ f }{ \F_{\R\setminus(t_{k-1},t_k)} } } \, \Big)^2 \,
;
\]
here $ t_0 = -\infty $, $ t_{n+1} = +\infty $, and the `$ \limsup $'
is taken over all finite sets $ L = \{ t_1, \dots, t_n \} \subset \R
$, $ t_1 < \dots < t_n $, ordered by inclusion. That is, `for every $
\eps $ there exists $ L_\eps $ such that for all $ L \supset L_\eps $
\dots' and so on. We also introduce
\[
\bH_1 (f) = \lim_{\{t_1,\dots,t_n\}\uparrow\,} \sum_{k=1}^{n+1} \Var \(
\cE{ f }{ \F_{t_{k-1},t_k} } \) \, .
\]
This time we may write `$ \lim $' (or `$ \inf $') instead of `$
\limsup $' due to monotonicity (w.r.t.\ inclusion); the more $ L = \{
t_1, \dots, t_n \} $ the less the sum.

\begin{lemma}\label{6.24}
\begin{sloppypar}
$ \sqrt{ \Var \( \cE f {\F_{s,t}} \) } \le \Ex \sqrt{ \cVar{ f }{
\F_{\R\setminus(s,t)} } } $ for all $ f \in L_2 (\Om,\F,P) $ and $ s <
t $.
\end{sloppypar}
\end{lemma}

\begin{proof}
The space $ L_2 (\Om,\F,P) = L_2(\F) = L_2 ( \F_{s,t}
\otimes \F_{\R\setminus(s,t)} ) = L_2 ( \F_{s,t} ) \otimes L_2 (
\F_{\R\setminus(s,t)} ) $ may also be thought of as the space $ L_2
\( \F_{\R\setminus(s,t)}, L_2 ( \F_{s,t} ) \) $ consisting of $
\F_{\R\setminus(s,t)} $-measurable square integrable
vector-functions, taking on values in $ L_2 ( \F_{s,t} ) $. We
consider the element $ \ti f \in L_2 \( \F_{\R\setminus(s,t)}, L_2
( \F_{s,t} ) \) $ corresponding to $ f \in L_2(\F) $ (according
to the canonical isomorphism of these two spaces). The mean value of
the vector-function is $ \Ex \ti f = \cE f {\F_{s,t}} $ (these two `$
\Ex $' act on different spaces). Convexity of
the seminorm $ \sqrt{ \Var(\cdot) } $ on $ L_2 (\F_{s,t}) $ gives $
\sqrt{ \Var(\Ex\ti f) } \le \Ex \sqrt{ \Var(\ti f) } $, where $
\Var(\ti f) $ means the pointwise variance (each value of $ \ti f $ is
a random variable; the latter has its variance), basically the same as
$ \cVar{ f }{ \F_{\R\setminus(s,t)} } $.
\qqed\end{proof}

\begin{corollary}\label{6.25}
$ \bH_1(f) \le \bH(f) $.
\end{corollary}

\begin{lemma}\label{6.26}
$ \bH_1(f) = \| Q_1 f \| $ for all $ f \in L_2 (\Om,\F,P) $; here $
Q_1 $ is the orthogonal projection onto the first chaos.
\end{lemma}

\begin{proof}
Follows immediately from Theorem \ref{6a2}.
\qqed\end{proof}

\begin{corollary}\label{6.27}
Every $ f \in L_2 (\Om,\F,P) $ such that $ \bH (f) = 0 $ is
orthogonal to the first chaos.
\end{corollary}

\begin{corollary}\label{6.28}
If a noise is such that $ \bH (f) = 0 $ for all $ f \in L_2
(\Om,\F,P) $, then the noise is black.
\end{corollary}

\begin{lemma}\label{6.29}
Let $ g \in L_2 (\F) $, $ h \in L_\infty (\F_{0,\infty})
$, and $ f = \cE{ gh }{ \F_{-\infty,0} } $. Then $ \bH (f) \le \| h
\|_\infty^2 \bH (g) $.
\end{lemma}

\begin{proof}
\begin{sloppypar}
It is sufficient to prove the inequality for the influence, $ \Ex
\sqrt{ \cVar{ f }{ \F_{\R\setminus(s,t)} } } \le \| h \|_\infty \Ex
\sqrt{ \cVar{ g }{ \F_{\R\setminus(s,t)} } } $ for any $ (s,t)
\subset (-\infty,0) $. Similarly to the proof of \ref{6.24}, we
consider $ \ti g \in L_2 \( \F_{0,\infty}, L_2(\F_{-\infty,0}) \) $
corresponding to $ g \in L_2 ( \F_{-\infty,0} \otimes \F_{0,\infty} )
$. We have $ \ti g h \in L_2 \( \F_{0,\infty}, L_2(\F_{-\infty,0} \) $, $
\Ex ( \ti g h ) = f $. Convexity of the seminorm $ \Ex \sqrt{ \cVar
\cdot {\F_{(-\infty,0)\setminus(s,t)}} } $ on $ L_2 ( \F_{-\infty,0} ) $ gives
$ \Ex \sqrt{ \cVar f {\F_{(-\infty,0)\setminus(s,t)}} } \le \Ex \Ex \sqrt{
\cVar{ \ti g h }{\F_{(-\infty,0)\setminus(s,t)}} } $, where `$ \Var $' and the
internal `$ \Ex $' act on $ L_2 ( \F_{-\infty,0} ) $, while the outer
`$ \Ex $' acts on $ L_2 ( \F_{0,\infty} ) $. The right-hand side is
equal to $ \Ex \Big( |h| \Ex \sqrt{ \cVar{ \ti g
}{\F_{(-\infty,0)\setminus(s,t)}} } \Big) $ and so, cannot exceed $ \|
h \|_\infty \Ex \Ex \sqrt{ \cVar{ \ti g
}{\F_{(-\infty,0)\setminus(s,t)}} } = \| h \|_\infty \Ex \sqrt{ \cVar{
g }{\F_{\R\setminus(s,t)}} } $.
\end{sloppypar}
\end{proof}

\begin{lemma}\label{6.30}
If $ f \in L_2 (\Om,\F,P) $ is such that $ \bH (f) = 0 $,
then $ \mu_f $ is concentrated on $ \cC_\perfect $.
\end{lemma}

\begin{proof}
Similarly to the proof of Theorem \ref{6.20}, it is sufficient to
prove, for every $ (a,b) \subset \R $, that $ \mu_f $-almost all $
M \in \cC $ are such that $ M \cap (a,b) $ is not a single-point
set. Lemma \ref{6.1} shows that the subspace corresponding to $ \{ M
\in \cC : | M \cap (a,b) | = 1 \} $ is $ H_{-\infty,a} \otimes
H^{(1)}_{a,b} \otimes H_{b,\infty} $, where $ H^{(1)}_{a,b} $ is the first
chaos intersected with $ H_{a,b} $. We have to prove that $ f $ is
orthogonal to $ H_{-\infty,a} \otimes H^{(1)}_{a,b} \otimes
H_{b,\infty} $, that is, to $ gh $ for every $ g \in H^{(1)}_{a,b} $,
$ h \in H_{-\infty,a} \otimes H_{b,\infty} = L_2 (
\F_{\R\setminus(a,b)} ) $, and we may assume that $ h \in  L_\infty (
\F_{\R\setminus(a,b)} ) $.

We have $ \Ex (fgh) = \Ex \( g \cE{ fh }{ \F_{a,b} } \) $. Lemma
\ref{6.29} (slightly generalized) shows that $ \bH \( \cE{ fh }{
\F_{a,b} } \) \le \| h \|_\infty^2 \bH (f) $. Thus, $ \bH
\( \cE{ fh }{ \F_{a,b} } \) = 0 $; by \ref{6.27}, $ \Ex \( g
\cE{ fh }{ \F_{a,b} } \) = 0 $.
\qqed\end{proof}

\begin{corollary}\label{6.31}
Let $ \( (\Om,\F,P), (\F_{s,t})_{s\le t} \) $ be a continuous
factorization. If $ f \in L_2 (\Om,\F,P) $ satisfies $ \bH (f) =
0 $ and $ \Ex f = 0 $, then $ f $ is sensitive.
\end{corollary}

Here are counterparts of \ref{5.5} and Theorem \ref{5.9}
inspired by the work \cite{LJR2} of Le Jan and Raimond.

\begin{lemma}\label{6.36}
Let $ f \in L_2 (\Om,\F,P) $, and $ g = \eta \circ f $
where $ \eta : \R \to \R $ satisfies $ | \eta(x) - \eta(y) | \le |x-y|
$ for all $ x,y \in \R $. Then
\[
\bH (g) \le \bH (f) \, .
\]
\end{lemma}

\begin{proof}
\begin{sloppypar}
It is sufficient to prove the inequality for the influence, $ \Ex
\sqrt{ \cVar{ g }{ \F_{\R\setminus(s,t)} } } \le \Ex \sqrt{ \cVar{ f
}{ \F_{\R\setminus(s,t)} } } $, or a stronger inequality $ \cVar{ g }{
\F_E } \le \cVar{ f }{ \F_E } $ a.s., for an arbitrary elementary set
$ E $. It is a conditional counterpart of the inequality $ \Var ( \eta
\circ X ) \le \Var (X) $ for any random variable $ X $. A proof of the
latter: $ \Var ( \eta \circ X ) = \frac12 \Ex ( \eta \circ X_1 - \eta
\circ X_2 )^2 \le  \frac12 \Ex ( X_1 - X_2 )^2 = \Var (X) $, where $
X_1, X_2 $ are independent copies of $ X $.
\end{sloppypar}
\end{proof}

\begin{theorem}\label{6.37}
\begin{sloppypar}
For every continuous factorization $ \( (\Om,\F,P), (\F_{s,t})_{s\le
t} \) $ there exists a sub-\sif\
$ \F_\jetblack $ of $ \F $ such that $ L_2(\Om,\F_\jetblack,P) $ is the
closure \textup{(}in $ L_2(\Om,\F,P) $\textup{)} of $ \{ f \in
L_2(\Om,\F,P) : \bH(f) = 0 \} $.
\end{sloppypar}
\end{theorem}

\begin{proof}
The set $ \{ f : \bH(f) = 0 \} $ is closed under linear operations,
and also under the nonlinear operation $ f \mapsto |f| $, therefore
its closure is of the form $ L_2(\F_\jetblack) $.
\qqed\end{proof}

\begin{corollary}
$ L_2 (\F_\jetblack) \subset H_{\cC_\perfect} $.
\end{corollary}

\begin{question}
Whether $ \F_\jetblack $ is nontrivial for every black noise, or not?
\end{question}

\section{Example: The Brownian Web as a Black Noise}
\label{sec:7}
\subsection{Convolution semigroup of the Brownian web}
\label{sec:7.1}

A one-dimensional array of random signs can produce some classical and
nonclassical noises in the scaling limit, but I still do not
know\pagebreak[1]
whether it can produce a black noise, or not (see \ref{6.23}).
\[
\setlength{\unitlength}{0.7cm}\nopagebreak[4]
\begin{picture}(12,3)
\put(0,0.5){\includegraphics[scale=0.7]{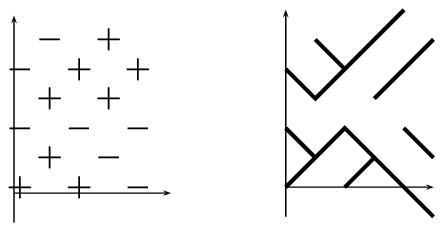}}
\put(1,0){\makebox(0,0){(a)}}
\put(3.8,0){\makebox(0,0){(b)}}
\put(6,0.5){\includegraphics[scale=0.7]{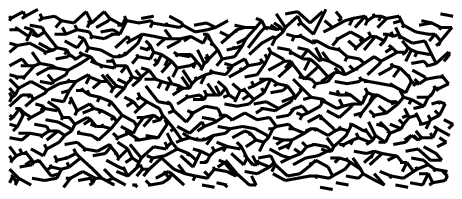}}
\put(8.5,0){\makebox(0,0){(c)}}
\end{picture}
\]
This is why I turn to a two-dimensional array of random
signs (a). It produces a system of coalescing random walks (b) that
converges to the so-called
\emph{Brownian web}\index{Brownian web}
(c), consisting of infinitely many coalescing Brownian motions
(independent before coalescence).

The Brownian web was investigated by Arratia, Toth, Werner, Soucaliuc,
and recently by Fontes, Isopi, Newman and Ravishankar \cite{FINR}
(other references may be found therein). The scaling limit may be
interpreted in several ways, depending on the choice of
`observables', and may involve delicate points, because of complicated
topological properties of the Brownian web as a random geometric
configuration on the plane. However, we avoid these delicate points by
treating the Brownian web as a stochastic flow in the sense of
Sect.~4, that is, a two-parameter family of random variables in a
semigroup.

In order to keep finite everything that can be kept finite, we
consider Brownian motions in the circle $ \T = \R/\Z $ rather than the
line $ \R $.

It is well-known that a countable dense set of coalescing `particles',
given at the initial instant, becomes finite, due to coalescence,
after any positive time. Moreover, the finite number is of finite
expectation. Thus, for any given $ t > 0 $, the Brownian
web on the time interval $ (0,t) $ gives us a random map $ \T \to \T $
of the following elementary form (a step function):
\[
\begin{gathered}\includegraphics{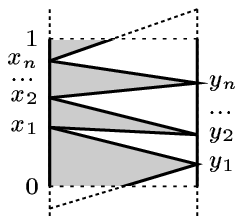}\end{gathered} \qquad
\begin{gathered}
f_{x_1,\dots,x_n}^{y_1,\dots,y_n} : \T \to \T \, , \\
x_1 < \dots < x_n < x_1, \> y_1 < \dots < y_n < y_1 \text{
 (cyclically),} \\
f_{x_1,\dots,x_n}^{y_1,\dots,y_n} (x) = y_{k+1} \text{ for } x \in
 (x_k,x_{k+1}] \, .
\end{gathered}
\]
Of course, $ n $ is random, as well as $ x_1,\dots,x_n $ and $
y_1,\dots,y_n $. The value at $ x_k $ does not matter; we let it be $ y_k $
for convenience, but it could equally well be $ y_{k+1} $, or remain
undefined. Points $ x_1,\dots,x_n $ will be called left critical
points of the map, while $ y_1,\dots,y_n $ are right critical points.

We introduce the set $ G_\infty $ consisting of all step functions $
\T \to \T $ and, in addition, the identity function. If $ f,g \in
G_\infty $ then their composition $ fg $ belongs to $ G_\infty $; thus
$ G_\infty $ is a
semigroup. It consists of pieces of dimensions $ 2,4,6,\dots $ and the
identity. Similarly to $ G_3 $ (recall \eqref{4d2}), $ G_\infty $ is
not a topological semigroup, since the composition is discontinuous.

The distribution of the random map is a probability measure $ \mu_t $
on $ G_\infty $. These maps form a convolution semigroup, $ \mu_s *
\mu_t = \mu_{s+t} $. Similarly to \ref{sec:4.5}, discontinuity
of composition
does not harm, since the composition is continuous almost everywhere
(w.r.t.\ $ \mu_s \otimes \mu_t $). Left and right critical points do not
meet.\footnote{%
 They meet with probability $ 0 $, as long as $ s $ and $ t $ are
 fixed. Otherwise, delicate points are involved\dots}

Having the convolution semigroup, we can construct the stochastic
flow, that is, a family of $ G_\infty $-valued random variables
$(\xi_{s,t})_{s\le t} $ such that
\begin{gather*}
\xi_{s,t} \sim \mu_{t-s} \, , \\
\xi_{r,s} \xi_{s,t} = \xi_{r,t} \quad \text{a.s.}
\end{gather*}
whenever $ -\infty < r < s < t < \infty $, and
\[
\xi_{t_1,t_2}, \dots, \xi_{t_{n-1},t_n} \quad \text{are independent}
\]
whenever $ -\infty < t_1 < \dots < t_n < \infty $.

Indeed, for each $ i $, we can take independent $ \xi_{k/i,(k+1)/i} :
\Om[i] \to G_\infty $ for $ k \in \Z $ according to the discrete
model, and define $ \xi_{k/i,l/i} =
\xi_{k/i,(k+1)/i} \dots \xi_{(l-1)/i,l/i} $. For any two coarse
instants $ s \le t $, the distribution of $ \xi_{s[i],t[i]} $
converges weakly (for $ i \to \infty $) to $ \mu_{t[\infty]-s[\infty]}
$. The refinement gives us
\[
\xi_{s,t} : \Om \to G_\infty \, , \qquad
\xi_{s,t} =
f_{x_1(s,t),\dots,x_{n(s,t)}(s,t)}^{y_1(s,t),\dots,y_{n(s,t)}(s,t)} \,
;
\]
$ x_k(\cdot,\cdot) $ and $ y_k(\cdot,\cdot) $ are continuous a.s. Also,
\begin{equation}\label{7a1}
\Ex n(s,t) < \infty \, .
\end{equation}
We consider the sub-\sif\ $ \F_{s,t} $ generated by all $ \xi_{u,v} $
for $ (u,v) \subset (s,t) $ and get a continuous factorization. Time
shifts are evidently introduced, and so, we get a noise --- the
\emph{noise of coalescence.}\index{noise of coalescence}

\subsection{Some general arguments}
\label{sec:7.2}

Probably we could use $ \bH $ and Theorem \ref{6.37} in order to prove
that the noise of coalescence is black (see also
\cite{LJR2}). However, I choose another way (via $ \bH_1 $ rather than
$ \bH $).

Random variables of the form $ \phi(\xi_{s,t}) $ for arbitrary $ s <
t $ and arbitrary bounded Borel function $ \phi : G_\infty \to \R $
generate the whole \sif\ $ \F $. Products of the form $ \phi_1
(\xi_{t_0,t_1}) \dots \phi_n (\xi_{t_{n-1},t_n}) $ for $ t_0 < \dots <
t_n $ span $ L_2 $ (as a closed subspace); however, we cannot expect
that \emph{linear} combinations of such $ \phi (\xi_{s,t}) $ are dense
in $ L_2 $.

Denote by $ Q_1 $ the orthogonal projection of $ L_2(\Om,\F,P) $ onto
the first chaos.

\begin{lemma}\label{7.1}

Linear combinations of all $ Q_1 \phi(\xi_{s,t}) $ are dense in the
first chaos.

\end{lemma}

Proof: Follows easily from the next (quite general) result, or rather,
its evident generalization to $ n $ factors.

\begin{lemma}

Let $ r \le s \le t $, $ X \in L_2 (\F_{r,s}) $, $ Y \in L_2
(\F_{s,t}) $. Then $ Q_1 (XY) = Q_1 (X) \Ex(Y) + \Ex(X) Q_1 (Y) $.

\end{lemma}

\begin{proof}
In terms of operators $ R_\phi $ given by \ref{5.3} we have $
Q_1 (XY) = R_{\phi_{r,t}} (XY) $, where $ \phi_{r,t} : \cC_{r,t} \to
\R $ is the indicator of $ \{ M \in \cC : |M\cap(r,t)|=1 \}
$. Similarly, $ Q_1 (X) = R_{\phi_{r,s}} (X) $, and $ \Ex (X) =
R_{\psi_{r,s}} (X) $, where $ \psi_{r,s} $ is the indicator of $ \{ M
\in \cC : |M\cap(r,s)| = 0 \} $. However, $ \phi_{r,t} = \phi_{r,s}
\psi_{s,t} + \psi_{r,s} \phi_{s,t} $ almost everywhere on $ \cC_{r,t}
$ (w.r.t.\ every spectral measure).
\qqed\end{proof}

In order to prove that the noise (of coalescence) is black, it
suffices to prove that $ Q \phi(\xi_{s,t}) = 0 $ for all $ s,t,\phi
$. We'll prove that $ Q \phi(\xi_{0,1}) = 0 $; the general case is
similar. According to \ref{6.26} we have to prove that $ \bH_1 (
\phi(\xi_{0,1}) ) = 0 $. Assuming that $ \Ex \phi(\xi_{0,1}) = 0 $ we
will check the sufficient condition:
\[
\| \cE{ \phi(\xi_{0,1}) }{ \F_{t-\eps,t} } \| = o (\sqrt\eps) \quad
\text{for } \eps \to 0 \, ,
\]
uniformly in $ t $. When doing so, we may assume that $ t $ is bounded
away from $ 0 $ and $ 1 $. Indeed, $ \| \cE{ \phi(\xi_{0,1}) }{
\F_{t,1} } \| \to 0 $ for $ t \to 1- $, due to continuity of the
factorization (recall \ref{3d1}(b)).

\begin{lemma}

$ \cE{ \phi(\xi_{0,1}) }{ \F_{t-\eps,t} } = \cE{ \phi(\xi_{0,1}) }{
\xi_{t-\eps,t} } $.

\end{lemma}

The proof is left to the reader; a hint:
\begin{multline*}
\cE{ \phi ( \xi_{t_1,t_5} ) }{ \xi_{t_2,t_3}, \xi_{t_3,t_4} } = \iint
 \phi ( \xi_{12} \xi_{23} \xi_{34} \xi_{45} ) \,
 \D\mu_{t_2-t_1}(\xi_{12}) \D\mu_{t_5-t_4}(\xi_{45}) \\
= \cE{ \phi ( \xi_{t_1,t_5} ) }{ \xi_{t_2,t_4} } \, .
\end{multline*}

\subsection{The key argument}
\label{sec:7.3}

Similarly to \ref{6a4}, we consider $ X = \phi(\xi_{0,1}) = \phi (
\xi_{0,t-\eps} \xi_{t-\eps,t} \xi_{t,1} ) $, $ \Ex X = 0 $, $ |X| \le
1 $ a.s. We have to prove that $ \| \cE{ X }{ \xi_{t-\eps,t} } \| =
o(\sqrt\eps) $ for $ \eps \to 0 $, uniformly in $ t $, when $ t $ is
bounded away from $ 0 $ and $ 1 $. Clearly,
\[
\cE{ X }{ \xi_{t-\eps,t} } = \iint \phi(fgh) \, \D\mu_{t-\eps}(f)
\D\mu_{1-t}(h) \, ,
\]
where $ g = \xi_{t-\eps,t} $.
\[
\begin{gathered}\includegraphics{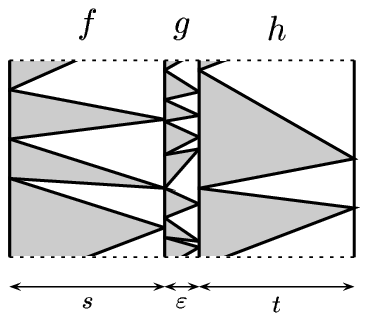}\end{gathered}
\]
We choose $ \ga \in \( \frac13, \frac12 \) $ and divide the strip $
(t-\eps,t) \times \T $ into $ \sim \eps^{-\ga} $ `cells' $ (t-\eps,t)
\times (z_k,z_{k+1}) $ of height $ z_{k+1} - z_k \sim \eps^\ga $.
\[
\begin{gathered}\includegraphics{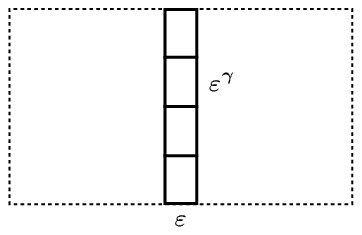}\end{gathered}
\]
We want to think of $ g $ as consisting of independent cells. Probably
it can be done in continuous time, but we have no such
technique for now. Instead, we retreat to the discrete-time model. The needed
inequality for continuous time results in the scaling limit $ i \to
\infty $ provided that in discrete time our estimations are uniform in
$ i $ (for $ i $ large enough).

So, random signs that produce $ g $ are divided into cells. Cells are
independent and, taken together, they determine $ g $ uniquely.

However, a path may cross many cells. This is rather improbable, since
$ \ga < 1/2 $, but it may happen. We enforce locality by a forgery!
Namely, if the path starting at the middle of a cell reaches the
bottom or the top edge of the cell, we replace the whole cell with some
other cell (it may be chosen once and for all) where it does not
happen.
\[
\begin{gathered}\includegraphics[scale=0.8]{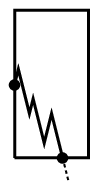}\end{gathered}\qquad\mapsto
\qquad\begin{gathered}\includegraphics[scale=0.8]{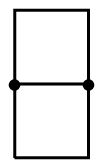}\end{gathered}
\]
Now cells are `local'; a path cannot cross more than two cells, but of
course, the stochastic flow is changed. Namely, $ g $ is changed with
an exponentially small (for $ \eps \to 0 $) probability, which changes
$ \cE{ X }{ \xi_{t-\eps,t} } $ by $ o(\sqrt\eps) $ (much less, in
fact). Still, cells are independent.

Does a cell (of $ g $) influence the composition, $ fgh \, $? It
depends on $ f $ and $ h $. If the left edge $ \{ t-\eps \} \times
[z_k,z_{k+1}] $ of the cell contains no right critical point of $ f $,
the cell can influence, since a path starting in an adjacent cell can
cross the boundary between cells. However, if the enlarged left edge $
\{ t-\eps \} \times [z_k-\eps^\ga,z_{k+1}+\eps^\ga] $ contains no
right critical point of $ f $ (in which case we say `the cell is
blocked by $ f $'), then the cell cannot influence, because of the
enforced locality. Similarly, if the enlarged right edge $
\{ t \} \times [z_k-\eps^\ga,z_{k+1}+\eps^\ga] $ contains no
left critical point of $ h $ (in which case we say `the cell is
blocked by $ h $'), the cell cannot influence.

The probability of being not blocked by $ f $ is the same for all
cells, since the distribution of $ f $ is invariant under rotations of
$ \T $ (discretized as needed). The sum of these probabilities does not
exceed $ 3 \Ex n(0,t-\eps) $ (recall \eqref{7a1}), which is $ O(1) $
when $ \eps \to 0 $. (Here we need $ t $ to be bounded away from $ 0
$.) Thus,
\begin{gather*}
\Pr{ \text{a given cell is not blocked by $ f $} } = O(\eps^\ga) \, ;
 \\
\Pr{ \text{a given cell is not blocked by $ h $} } = O(\eps^\ga) \, ;
 \\
\Pr{ \text{a given cell is not blocked} } = O(\eps^{2\ga}) \, ;
 \\
\Pr{ \text{at least one cell is not blocked} } = O(\eps^\ga) \, .
\end{gather*}
In the latter case we may say that $ g $ is not blocked (by $ f,h $).

Denote by $ A $ the event ``$ g $ is not blocked by $ f,h $'' (it is
determined by $ f $ and $ h $, not $ g $); $ \Pr{ A } = O(\eps^\ga)
$. Taking into account that
\begin{gather*}
X = X - \Ex X = \( X \cdot \One_A - \Ex ( X \cdot \One_A ) \) + \( X
 \cdot (\One-\One_A) - \Ex ( X \cdot (\One-\One_A) ) \) \, , \\
\cE{ X \cdot (\One-\One_A) }{ g } = \Ex ( X \cdot (\One-\One_A) ) \, ,
 \\
\cE X g = \cE{ X \cdot \One_A }{ g } - \Ex ( X \cdot \One_A ) \, ,
\end{gather*}
we have to prove that $ \| \cE{ X \cdot \One_A }{ g } - \Ex ( X \cdot
\One_A ) \| = o(\sqrt\eps) $. Note that it does not result from the
trivial estimation $ \| X \cdot \One_A \| \le \| \One_A \| = \sqrt{\Pr
A} = O(\eps^{\ga/2}) $, $ \ga \in \(\frac13,\frac12\) $. Note also
that, when $ g $ influences $ X $, its influence is usually not small
(irrespective of $ \eps $) because of the stepwise nature of $ f $ and
$ h $.

We express the norm in terms of covariance,
\[
\| \cE{ X \cdot \One_A }{ g } - \Ex ( X \cdot \One_A ) \| = \sup_\psi
\Cov \( X \cdot \One_A, \psi(g) \) \, ,
\]
where the supremum is taken over all Borel functions $ \psi : G_\infty
\to \R $ such that $ \Var \( \psi(g) \) \le 1 $. In terms of the
correlation coefficient
\[
\Corr \( X \cdot \One_A, \psi(g) \) = \frac{ \Cov \( X \cdot \One_A,
\psi(g) \) }{ \sqrt{ \Var ( X \cdot \One_A ) } \sqrt{ \Var ( \psi(g) )
} } \, ,
\]
it is enough to prove that
\[
\Corr \( X \cdot \One_A, \psi(g) \) = o ( \eps^{(1-\ga)/2} ) \, ,
\]
since it implies $ \Cov(\dots) = o ( \eps^{(1-\ga)/2} ) \cdot \| X
\cdot \One_A \| = o ( \eps^{(1-\ga)/2} \eps^{\ga/2} ) = o(\sqrt\eps)
$. Instead of $ o ( \eps^{(1-\ga)/2} ) $ we will get $ O(\eps^\ga) $,
which is also enough since $ \ga > 1/3 $.

It remains to apply the quite general lemma given below, interpreting
its $ Y_k $ as the whole $ k $-th cell (of $ g $), $ X_k $ as the
indicator of the event ``the $ k $-th cell is not blocked'' ($ k =
1,\dots,n $), $ X_0 $ as the pair $ (f,h) $, and $ \phi(\dots) $ as $
X \cdot \One_A $. The lemma is
formulated for real-valued random variables $ Y_k $, but this does not
matter; the same clearly holds for arbitrary spaces, and in fact, we
need only finite spaces. The product $ X_k Y_k $ is a trick for
`blocking' $ Y_k $ when $ X_k = 0 $. Note that dependence between $
X_0, X_1, \dots, X_n $ is allowed.

\begin{lemma}\label{7.4}

Let $ (X_0,X_1,\dots,X_n) $ and $ (Y_1,\dots,Y_n) $ be two independent
random vectors, $ Y_k : \Om \to \R $, $ X_k : \Om \to \{0,1\} $ for $
k = 1,\dots,n $, $ X_0 : \Om \to \R $, and random variables $
Y_1,\dots,Y_n $ be independent. Then
\[
\Corr \( \phi (X_0, X_1 Y_1,\dots,X_n Y_n), \, \psi
 (Y_1,\dots,Y_n) \) \le \sqrt{ \max_{k=1,\dots,n} \Pr{ X_k=1 } }
\]
for all Borel functions $ \phi : \R^{n+1} \to \R $, $ \psi : \R^n \to
\R $ such that the correlation is well-defined \textup{(}that is, $ 0
< \Var \phi(\dots) < \infty $, $ 0 < \Var \psi(\dots) < \infty
$\textup{).}

\end{lemma}

\begin{proof}

\begin{sloppypar}
We may assume that $ X_1, \dots, X_n $ are functions of $ X_0 $.
Consider the orthogonal (in $ L_2(\Om) $) projection $ Q $ from the
space of all random variables of the form $ \psi (Y_1,\dots,Y_n) $ to
the space of all random variables of the form $ \phi (X_0, X_1
Y_1,\dots,X_n Y_n) $, that is, $ Q \psi (Y_1,\dots,Y_n) = \cE{ \psi
(Y_1,\dots,Y_n) }{ X_0, X_1 Y_1,\dots, X_n Y_n } $. We have to prove
that $ \| Q \psi (Y_1,\dots,Y_n) \|^2 \le \( \max_k \Pr{ X_k=1 } \) \|
\psi (Y_1,\dots,Y_n) \|^2 $ whenever $ \Ex \psi (Y_1,\dots,Y_n) = 0 $.
The space of all $ \psi (Y_1,\dots,Y_n) $ is spanned by
\emph{factorizable} random variables $ \psi (Y_1,\dots,Y_n) = \psi_1
(Y_1) \dots \psi_n (Y_n) $. For such a $ \psi $ we have
\begin{multline*}
Q \psi (Y_1,\dots,Y_n) = \cE{ \psi_1 (Y_1) \dots \psi_n (Y_n) }{ X_0,
 X_1 Y_1,\dots,X_n Y_n } \\
= \bigg( \prod_{k:X_k=0} \Ex
 \psi_k(Y_k) \bigg) \bigg( \prod_{k:X_k=1} \psi_k(Y_k) \bigg) \, ;
\end{multline*}
\vspace{-3mm}
\begin{multline*}
\| Q \psi (Y_1,\dots,Y_n) \|^2 = \Ex \Big( \cE{ | Q \psi
 (Y_1,\dots,Y_n) |^2 }{ X_0 } \Big) \\
= \Ex \left( \bigg(
 \prod_{k:X_k=0} | \Ex \psi_k(Y_k) |^2 \bigg) \bigg( \prod_{k:X_k=1}
 \Ex | \psi_k(Y_k) |^2 \bigg) \right) \, .
\end{multline*}
If, in addition, $ \Ex \psi_1 (Y_1) = 0 $ then $ \| Q \psi
(Y_1,\dots,Y_n) \|^2 \le \linebreak
\Pr{ X_1=1 } \| \psi (Y_1,\dots,Y_n) \|^2 $. Similarly,
\[
\| Q \psi (Y_1,\dots,Y_n) \|^2 \le \Big( \max_k \Pr{ X_k=1 } \Big)
\| \psi (Y_1,\dots,Y_n) \|^2 
\]
if $ \Ex \psi (Y_1,\dots,Y_n) = 0 $ and, of course, $ \psi $ is
factorizable, that is, $ \psi (Y_1,\dots,Y_n) = \linebreak[0]
\psi_1 (Y_1) \dots
\psi_n (Y_n) $. The latter assumption cannot be eliminated just by
saying that factorizable random variables of zero mean span all random
variables of zero mean. Instead, we use two facts.
\end{sloppypar}

The first fact. The
space of all random variables $ \psi(\dots) $ has an
\emph{orthogonal} basis consisting of factorizable random
variables satisfying an additional condition: each factor $ \psi_k
(Y_k) $ is either of zero mean, or equal to $ 1 $. (For a proof, start
with an orthogonal basis for functions of $ Y_1 $ only, the first
basis function being constant; do the same for $ Y_2 $; take all
products; and so on.)

\begin{sloppypar}
The second fact. The operator $ Q $ maps \emph{orthogonal}
factorizable random variables, satisfying the additional condition,
into \emph{orthogonal}
random variables. Indeed, let $ \psi (Y_1,\dots,Y_n) = \psi_1 (Y_1)
\dots \psi_n (Y_n) $, $ \psi' (Y_1,\dots,Y_n) = \psi'_1 (Y_1) \dots
\psi'_n (Y_n) $, and each $ \psi_k (Y_k) $ be either of zero mean, or
equal to $ 1 $; the same for each $ \psi'_k(Y_k) $. If $ \Ex \( \psi 
(Y_1,\dots,Y_n) \psi' (Y_1,\dots,Y_n) \) = 0 $ then $ \Ex \( \psi_k
(Y_k) \psi'_k (Y_k) \) = 0 $ for at least one $ k $; let it happen for
$ k=1 $. We have not only $ \Ex \( \psi_1 (Y_1) \psi'_1 (Y_1) \) = 0 $
but also $ \( \Ex \psi_1 (Y_1) \) \( \Ex \psi'_1 (Y_1) \) = 0 $, since
$ \psi_1 $ and $ \psi'_1 $ cannot both be equal to $ 1 $. Therefore
\begin{multline*}
\Ex \( Q \psi (Y_1,\dots,Y_n) \) \( Q \psi' (Y_1,\dots,Y_n) \) = \\
= \Ex \left( \bigg( \prod_{k:X_k=0} \( \Ex \psi_k (Y_k) \) \( \Ex
 \psi'_k (Y_k) \) \bigg) \bigg( \prod_{k:X_k=1} \psi_k (Y_k) \psi'_k
 (Y_k) \bigg) \right) = 0 \, ,
\end{multline*}
since the first term vanishes whenever $ X_1 = 0 $, and the second
term vanishes whenever $ X_1 = 1 $.
\qqed\end{sloppypar}
\end{proof}

Combining all together, we get the conclusion.

\begin{theorem}\label{7c2}
The noise of coalescence is black.
\end{theorem}

\subsection{Remarks}
\label{sec:7.4}

Another proof of Theorem \ref{7c2} should be possible, by showing that
all (zero mean) random variables are sensitive. To this end, we
divide the time axis $ \R $ into intervals of small length $ \eps $, and
choose a random subset of intervals such that each interval is chosen
with a small probability $ 1 - \rho = 1 - \E^{-\la} \sim \la $,
independently of others. On each chosen interval we replace local
random data with fresh (independent) data.

Consider the path $ X(\cdot) $ of the Brownian web, starting at the
origin, $ X(t) = \xi_{0,t} (0) $ for $ t \in [0,\infty) $; it behaves
like a Brownian motion. After the replacement we get another path $
Y(\cdot) $. Their difference, $ \( X(t) - Y(t) \) / \sqrt2 $, behaves
like another Brownian motion when outside $ 0 $, but is somewhat
sticky at $ 0 $. Namely, during each chosen (to the random set) time
interval, the point $ 0 $ has nothing special; however, outside these
time intervals, the point $ 0 $ is absorbing. In this sense, chosen
time intervals act like factors $ f_* $ in the random product of
factors $ f_-, f_+, f_* $ studied in Sect.~4. There, $ f_* $ occurs
with a small probability $ 1/(2\sqrt i) \to 0 $ (recall \ref{4ee3}),
which produces a non-degenerate stickiness in the scaling limit. Here,
in contrast, a time interval is chosen with probability $ 1-\rho \sim
\la $ that does not tend to $ 0 $ when the interval length $ \eps $
tends to $ 0 $. Naturally, stickiness disappears in the limit $ \eps
\to 0 $ (a proof uses the idea of \eqref{4c4c}). That is, interaction
between $ X(\cdot) $ and $ Y(\cdot) $ disappears in the limit $ \eps
\to 0 $. They become independent, no matter how small $ 1 - \rho $ is.

Probably, the same argument works for any finite number of paths $ X_k
(t) = \xi_{0,t} (x_k) $; they should be asymptotically independent of
$ Y_k(\cdot) $ for $ \eps \to 0 $, but I did not prove it.

\smallskip

The spectral measure $ \mu_X $ of the random variable $ X = \xi_{0,1}
(0) $ is written down explicitly in \cite{TsFW}. Or rather, its
discrete counterpart is found; the scaling limit follows by (a
generalization of) Theorem \ref{3c5} (see also \cite{TsSc}). The
measure $ \mu_X $ is a probability measure (since $ \| X \| = 1 $), it
may be thought of as the distribution of a random perfect subset of $
(0,1) $. Note that the random subset is not at all a function on the
probability space $ (\Om,\F,P) $ that carries the Brownian web. There
is no sense in speaking about `the joint distribution of the random
set and the Brownian web'. In fact, they may be treated as
incompatible (non-commuting) measurements in the framework of quantum
probability, see
\cite{Ts98}.

A wonder: $ \mu_X $ is the distribution of $ ( \theta - M ) \cap (0,1)
$, where $ M $ is the set of zeros of the usual Brownian motion, and $
\theta $ is independent of $ M $ and distributed uniformly on $ (0,1)
$.

Moreover, the corresponding equality holds exactly (not only
asymptotically) in the discrete-time model. Strangely enough, the
Brownian motion (or rather, random walk) does not appear in the
calculation of the spectral measure. The relation to Brownian motion
is observed at the end, as a surprise!

\begin{question}

Can $ \mu_X $ (for $ X = \xi_{0,1}(0) $) be found via some natural
construction of a Brownian motion whose zeros form the spectral set
(after the transformation $ x \mapsto \theta - x $)? (See
\cite[Problem 1.5]{TsFW}.)

\end{question}

We see that $ \mu_X $ (for $ X = \xi_{0,1}(0) $)  is concentrated on
sets of Hausdorff dimension $ 1/2 $.

\begin{question}

Is $ \mu_X $ concentrated on sets of Hausdorff dimension $ 1/2 $ for
an arbitrary random variable $ X $ such that $ \Ex X = 0 $ (over the
noise of coalescence)?

\end{question}

An affirmative answer would probably give us another proof that the
noise is black. A stronger conjecture may be made.

\begin{question}

Is $ \mu_X $ for an arbitrary $ \F_{0,1} $-measurable $ X $ (over the
noise of coalescence), satisfying $ \Ex X = 0 $, absolutely continuous
w.r.t.\ $ \mu_{\xi_{0,1}(0)} \, $?

\end{question}

\subsection{A combinatorial by-product}

Consider a Markov chain $ X = (X_k)_{k=0}^\infty $ (a
half-difference of two independent simple random walks, or a
double-speed simple random walk divided by two): $ X_0 = 0 $ and
\[
\cP{ X_{k+1}=X_k+\De x }{ X_k } =
\begin{cases}
 1/4 &\text{for $ \De x = -1 $},\\
 1/2 &\text{for $ \De x = 0 $},\\
 1/4 &\text{for $ \De x = +1 $}
\end{cases}
\]
for each $ k = 0,1,2,\dots $

Let $ Z $ be the (random) set of zeros of $ X $, that is,
\[
Z = \{ k=0,1,\dots\ : X_k = 0 \} \, .
\]
Given a set $ S \subset \{ 0,1,2,\dots \} $ and a number $ k =
0,1,2,\dots $, we consider the event $ Z \cap [0,k] \subset k - S $,
that is, $ \forall l = 0,\dots,k \; \( l \in Z \impl k-l \in S \) $,
and its probability. We define
\[
p_{n,S} = \frac1n \sum_{k=0}^{n-1} \Pr{ Z \cap [0,k] \subset k - S }
\, ;
\]
of course, only $ k \in S $ can contribute (since $ 0 \in Z $).

On the other hand, we may trap $ X $ at $ 0 $ on $ S $; that is, given a
set $ S \subset \{ 0,1,2,\dots \} $, we introduce another Markov chain
$ X^{(S)} = \(X^{(S)}_k\)_{k=0}^\infty $ such that $ X^{(S)}_0 = 0 $
and for each $ k = 0,1,2,\dots $
\[
\cP{ X^{(S)}_{k+1}=x+\De x }{ X^{(S)}_k = x } =
\begin{cases}
 1/4 &\text{for $ \De x = -1 $},\\
 1/2 &\text{for $ \De x = 0 $},\\
 1/4 &\text{for $ \De x = +1 $}
\end{cases}
\]
except for the case $ k \in S $, $ x = 0 $,
\[
\cP{ X^{(S)}_{k+1} = 0 }{ X^{(S)}_k = 0 } = 1 \quad \text{if } k \in S
\, .
\]

\begin{theorem}\label{7.9}
$ p_{n,S} = \frac1n \sum_{k\in S} \Pr{ X^{(S)}_k = 0 } $ for every $ n
= 1,2,\dots $ and $ S \subset \{ 0,1,\dots,n-1 \} $.
\end{theorem}

\begin{example}
Before proving the theorem, consider a special case; namely, let
$ S $ consist of just a single number $ s $. Then $ \Pr{ Z \cap
[0,k] \subset k - S } = \Pr{ Z \cap [0,k] \subset \{ k - s \} } $
vanishes for $ k \ne s $. For $ k = s $ it becomes $ \Pr{ Z \cap [0,s]
= \{0\} } = 2^{-(2s-1)} \( \binom{2s-2}{s-1} + \binom{2s-2}{s} \) $.
Therefore $ p_{n,\{s\}} = \frac1n 2^{-(2s-1)} \( \binom{2s-2}{s-1} +
\binom{2s-2}{s} \) $, assuming $ s \ge 2 $; also, $ p_{n,\{0\}} =
\frac1n $ and $ p_{n,\{1\}} = \frac1{2n} $. On the other hand, $
\frac1n \sum_{k\in S} \Pr{ X^{(S)}_k = 0 } = \frac1n \Pr{ X_s = 0 } =
\frac1n \cdot 2^{-2s} \binom{2s}{s} $. The equality becomes $
\binom{2s-2}{s-1} + \binom{2s-2}{s} = \frac12 \binom{2s}{s} $ (for $ s
\ge 2 $).
\end{example}

\begin{proof}[Proof (sketch).]
We use the discrete-time counterpart of the Brownian web (see
\ref{sec:7.1} and \cite[Sect.~1]{TsFW}) and consider $ \xi_{0,n}(0) $,
the value at time $ n $ of the path starting at the origin. At every
instant $ k \notin S $ we replace the corresponding random signs with
fresh (independent) copies, which leads to another random variable $
\xi'_{0,n}(0) $. We calculate the covariance $ \Ex \( \xi_{0,n}(0)
\xi'_{0,n}(0) \) $ in two ways, and compare the results.

\emph{The first way.}
The difference process $ \xi_{0,\cdot}(0) - \xi'_{0,\cdot}(0) $ is
distributed like the process $ 2X^{(S)} $ (similarly to
\ref{sec:7.4}). Thus
\[
4 \Ex \( X^{(S)}_n \)^2 = \Ex \( \xi_{0,n}(0) - \xi'_{0,n}(0) \)^2 =
2n - 2 \Ex \( \xi_{0,n}(0) \xi'_{0,n}(0) \) \, .
\]
On the other hand, $ \frac12 - \Ex \( X^{(S)}_{k+1} \)^2 + \Ex \(
X^{(S)}_k \)^2 = \frac12 \Pr{ X^{(S)}_k = 0 } $ if $ k \in S $,
otherwise $ 0 $. Therefore $ n - 2 \Ex \( X^{(S)}_n \)^2 = \sum_{k\in
S} \Pr{ X^{(S)}_k = 0 } $. So,
\[
\Ex \( \xi_{0,n}(0) \xi'_{0,n}(0) \) = \sum_{k\in S} \Pr{ X^{(S)}_k
= 0 } \, .
\]

\emph{The second way.}
In terms of the spectral measure $ \mu $ of the random variable $
\xi_{0,n}(0) $ we have $ \Ex \( \xi_{0,n}(0) \xi'_{0,n}(0) \) = \mu \{
M : M \subset S \} $. However, the probability measure $ \frac1n \mu $
is equal to the distribution of $ (\theta-Z) \cap [0,\infty) $; here $
Z $ is (as before) the set of zeros of $ X $, and $ \theta $ is a
random variable independent of $ Z $ and distributed uniformly on $ \{
0,1,\dots,n-1 \} $. (See \cite[Prop.~1.3]{TsFW}, see also \cite{WW}.)
Therefore $ \frac1n \mu \{ M : M \subset S \} = \Pr{ (\theta-Z) \cap
[0,\infty) \subset S } = \Pr{ Z \cap [0,\theta] \subset \theta - S } =
p_{n,S} $. So,
\[
\Ex \( \xi_{0,n}(0) \xi'_{0,n}(0) \) = n p_{n,S} \, .
\]
\qqed\end{proof}

\begin{question}
Is there a simpler proof of Theorem \ref{7.9}? Namely, can we avoid
the spectral measure and its relation to the set of zeros?
\end{question}

A continuous-time counterpart of Theorem \ref{7.9} is left to the
reader.

\section{Miscellany}
\label{sec:8}
\subsection{Beyond the one-dimensional time}
\label{sec:8.1}

Scaling limits of models driven by \emph{two-dimensional} arrays of
random signs are evidently important. The best examples appear in
percolation theory. Also the Brownian web is an example and, after
all, it may be treated as an oriented percolation.

In such cases, independent sub-\sif s should correspond to disjoint
regions of $ \R^2 $, not only of the form $ (s,t) \times \R $. In
fact, a rudimentary use of these can be found in Sect.~\ref{sec:7}
(recall `cells' in \ref{sec:7.3}). In general it is unclear
what kind of regions can be used; probably, regions with piecewise
smooth boundaries always fit, while arbitrary open sets do not fit
unless the two-dimensional noise is classical (recall \ref{sec:6.3}).

In spite of the great and spectacular progress of the percolation
theory (see for instance \cite{SW} and references therein), `the noise
of percolation'\index{noise of percolation} is still a dream.

\begin{question}
For the critical site percolation on the triangular lattice, invent an
appropriate coarse \sif, and check two-dimensional counterparts of the
two conditions of \ref{3b1} for an appropriate class of
two-dimensional domains. Is it possible?
\end{question}

\begin{remark}
Hopefully, the answer is affirmative, that is, the two-dimensional
noise of percolation will be defined. Then it should appear to be a
(two-dimensional) black noise, due to (appropriately adapted)
\ref{6.28}, \ref{7.1} and (most important) the
critical exponent for a small cell of size $ \eps \times \eps $
being pivotal \cite[Sect.~5, Item 2]{SW}. The probability is $
O(\eps^{5/4}) $, therefore $ o(\eps) $. The sum for $ \bH (f) $
contains $ O(1/\eps^2) $ terms, $ o(\eps^2) $ each.\footnote{%
 Different arguments (especially, \ref{7.4}) are used in Sect.\
 \ref{sec:7}, since an infinite two-dimensional spectral set could
 have a finite one-dimensional projection.}

Sensitivity of percolation events, disclosed in \cite{BKS}, is
micro-sensitivity (recall \ref{sec:5.3}). Existence of the black
noise of percolation would mean a stronger property: block
sensitivity. (See also \cite[Problem~5.4]{BKS}.)

It would be the most important example of a black noise!
\end{remark}

For the \emph{general} theory of stability, spectral measures,
decomposable processes etc., the dimension of the underlying space is
of little importance. Basically, regions must form a Boolean
algebra. Such a general approach is used in \cite{TV}, \cite{Ts99}.

Nonclassical factorizations appear already in zero-dimensional `time',
be it a Cantor set, or even a convergent sequence with limit point. For
Cantor sets, see \cite[Sect.~4]{TV}; some interesting models of
combinatorial nature, with large symmetry groups (instead of `time
shifts' of a noise) are examined there. For a convergent sequence with
limit point, see Chapter 1 here (namely, \ref{1b1}), and
\cite[Appendix]{Ts99}.

\subsection{The `wave noise' approach}
\label{sec:8.2}

A completely different way of constructing noises is sketched here.

Consider the linear wave equation in dimension $ 1+1 $,
\begin{equation}\label{eq7.1}
\bigg( \frac{\pd^2}{\pd t^2} - \frac{\pd^2}{\pd x^2} \bigg) u(x,t) = 0
\, ,
\end{equation}
with initial conditions $ u(x,0) = 0 $, $ u_t (x,0) = f(x) $. Its
solution is well-known:
\begin{equation*}
u(x,t) = \frac12 \int_{x-t}^{x+t} f(y) \, dy = \frac12 F(x+t) -
\frac12 F(x-t) \, ,
\end{equation*}
where $ F $ is defined by $ F'(x) = f(x) $. The formula holds in a
generalized sense for nonsmooth $ F $, which covers the following
case: $ F(x) = B(x) = $ Brownian motion (combined out of two
independent branches, on $ [0,+\infty) $ and on $ (-\infty,0] $); $
f(x) = B'(x) $ is the white noise. The random field on $ (-\infty,\infty)
\times [0,\infty) $,
\begin{equation*}
u(x,t) = \frac12 B(x+t) - \frac12 B(x-t) \, , \qquad B = \text{
Brownian motion,}
\end{equation*}
is continuous, stationary in $ x $, scaling invariant (for any $ c $
the random field $ u(cx,ct) / \sqrt c $ has the same distribution as $
u(x,t) $), satisfies the wave equation \eqref{eq7.1} and the following
independence condition:
\begin{equation}\label{eq7.2}
\begin{aligned}
& u \big|_L \;\; \text{ and } \;\; u \big|_R \;\; \text{ are
  independent,} \\
& \quad \text{where } \;\; L = \{ (x,t) : x < -t < 0 \}, \; R = \{
  (x,t) : x > t > 0 \}.
\end{aligned}
\begin{aligned}\includegraphics{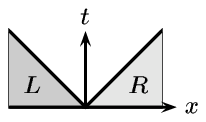}\end{aligned}
\end{equation}
The independence is a manifestation of: (1) the independence inherent
to the white noise (its integrals over disjoint segments are
independent), and (2) the hyperbolicity of the wave equation
(propagation speed does not exceed $ 1 $).

A solution with such properties is essentially unique. That is, if $
u(x,t) $ is a continuous random field on $ (-\infty,\infty) \times
(0,\infty) $, stationary in $ x $, satisfying the wave equation \eqref{eq7.1}
and the independence condition \eqref{eq7.2}, then necessarily $ u(x,t) =
\mu_0 + \mu_1 t + \sigma \( B(x+t) - B(x-t) \) $ for a Brownian motion
$ B $. Scaling invariance forces $ \mu_0 = \mu_1 = 0 $.

It is instructive that a wave equation may be used in a
non-traditional way. Traditionally, a solution is determined by its
initial values. In contrast, the independence condition \eqref{eq7.2},
combined with some more conditions, determines a random solution with
no help of initial conditions! Not an individual sample function is
determined, of course, but its distribution (a probability measure on
the space of solutions of the wave equation).

Somebody with no preexisting idea of white noise or Brownian motion
can, in principle, use the above approach. Observing that $ u(x,0) = 0
$ but $ u_t (x,0) $ does not exist (in the classical sense), he may
investigate $ u(x,t) / t $ for $ t \to 0 $ as a way toward the white
noise.

\begin{question}

Can we construct a nonclassical (especially, black) noise, using a
nonlinear hyperbolic equation?

\end{question}

I once tried the nonlinear wave equation
\begin{equation}\label{eq7.3}
\bigg( \frac{\pd^2}{\pd t^2} - \frac{\pd^2}{\pd x^2} \bigg) u(x,t) =
\eps t^{-(3-\eps)/2} \sin \( t^{-(1+\eps)/2} u(x,t) \) \, ,
\end{equation}
$ \eps $ being a small positive parameter. The equation is
scaling-invariant: if $ u(x,t) $ is a solution, then $ u(cx,ct) /
c^{(1+\eps)/2} $ is also a solution. We search for a random field $
u(t,x) $, continuous, stationary in $ x $, scaling invariant,
satisfying \eqref{eq7.3} and the independence condition \eqref{eq7.2}. Its behavior
for $ t \to 0 $ should give us a new noise. Does such a random field
exist? Is it unique (in distribution)? If the answers are affirmative,
then we get a noise,
\[
\F_{x,y} \;\; \text{ is the \sif\ generated by } \{ u(z,t) : x+t < z <
y-t \} \, , \;\; \begin{aligned}\includegraphics{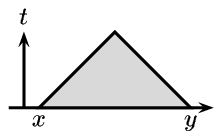}\end{aligned}
\]
and maybe it is black. However, I did not succeed with it.

A modified `waive noise' approach was used successfully in
\cite[Sect.~5]{TV}, proving, for the first time, the existence of a
black noise. The modification is to keep the auxiliary dimension, but
make it discrete rather than continuous:
\[
\begin{aligned}\includegraphics[scale=0.7]{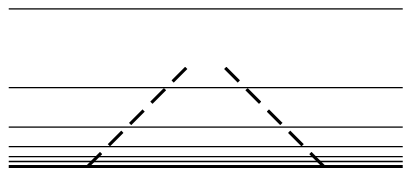}\end{aligned}
\]
More specifically, consider a sequence of stationary random processes
$ u_k (\cdot) $ on $ \R $ such that
\begin{myitemize}
\item $ u_k $ is $ 2\eps_k $-dependent (for some $ \eps_k \to 0 $); it
means that $ u_k \big|_{(-\infty,-\eps_k]} $ and $ u_k
\big|_{[\eps_k,+\infty)} $ are independent;
\item $ u_{k-1} (x) $ is uniquely determined by $ u_k \big|_{
[x-(\eps_{k-1}-\eps_k), \, x+(\eps_{k-1}-\eps_k)] } \, $.
\end{myitemize}
Such a sequence $ (u_k) $ determines a noise; namely, $ \F_{x,y} $ is
generated by all $ u_k (z) $ such that $ x+\eps_k \le z \le y-\eps_k
$. White noise can be obtained by a linear system of Gaussian
processes:
\[
u_{k-1} (x) = \int_{x-(\eps_{k-1}-\eps_k)}^{x+(\eps_{k-1}-\eps_k)} V_k
(y-x) u_k (y) \, dy \, ,
\]
where kernels $ V_k $, concentrated on $
[-(\eps_{k-1}-\eps_k), (\eps_{k-1}-\eps_k)] $, are chosen
appropriately. A nonlinear system (of quite non-Gaussian processes) of
the form
\[
u_{k-1} (x) = \phi \bigg( \frac{\const}{\eps_{k-1}-\eps_k}
\int_{x-(\eps_{k-1}-\eps_k)}^{x+(\eps_{k-1}-\eps_k)} u_k (y) \, dy
\bigg)
\]
was used for constructing a black noise. But, it is not really a
\emph{construction} of a specific noise. Existence of $ (u_k) $ is
proven, but uniqueness (in distribution) is not. True, every such $
(u_k) $ determines a black noise. However, none of them is singled
out.

\subsection{Groups, semigroups, kernels}
\label{sec:8.3}

A Brownian motion $ X $ in a topological group $ G $ is defined as a
continuous $ G $-valued random process with stationary independent
increments, starting from the unit of $ G $. For example, if $ G $ is
the additive group of reals, then the general form of a Brownian
motion in $ G $ is $ X(t) = \si B(t) + vt $, where $ B(\cdot) $ is the
standard Brownian motion, $ \si \in [0,\infty) $ and $ v \in \R $ are
parameters. If $ G $ is a Lie group, then Brownian motions $ X $ in $
G $ correspond to Brownian motions $ Y $ in the tangent space of $ G $
(at the unit) via the stochastic differential equation $ (dX) \cdot
X^{-1} = dY $ (in the sense of Stratonovich).

A noise corresponds to every Brownian motion in a topological group,
just as the white noise corresponds to $ B(\cdot) $. If the noise is
classical, it is the white noise of some dimension ($ 0,1,2,\dots $ or
$ \infty $). If this is the case for all Brownian motions in $ G $,
we call $ G $ a \emph{white group.}\index{white group}
Thus, $ \R $ is white, and
every Lie group is white. Every commutative topological group is white
(see \cite[Th.~1.8]{Ts98}). The group of all unitary operators in $
l_2 $ (equipped with the strong operator topology) is white (see
\cite[Th.~1.6]{Ts98}). Many other groups are white since they are
embeddable into a group known to be white; for example, the group of
diffeomorphisms is white (an old result of Baxendale).

\begin{question}

Is the group of all homeomorphisms of (say) $ [0,1] $ white?

\end{question}

In a topological group, Brownian motions $ X $ and continuous abstract
stochastic flows $ \xi $ are basically the same:
\[
X(t) = \xi_{0,t} \, ; \qquad \xi_{s,t} = X^{-1}(s) X(t) \, .
\]
In a semigroup, however, a noise corresponds to a flow, not to a
Brownian motion (see also \ref{4ce2}).

A nonclassical noise (of stickiness) was constructed in Sect.~4 out of
an abstract flow in a $ 3 $\nobreakdash-\hspace{0pt}dimensional
semigroup $ G_3 $; however, $ G_3 $ is not a topological semigroup,
since composition is discontinuous.

\begin{question}

Can a nonclassical noise arise from an abstract stochastic flow in a
finite-dimensional topological semigroup?

\end{question}

The continuous (but not topological) semigroup $ G_3 $ emerged in
Sect.~4 from the discrete semigroup $ G_3^\discrete $ via the scaling
limit. Or rather, a flow in $ G_3 $ emerged from a flow in $
G_3^\discrete $ via the scaling limit. A similar approach to the discrete
model of \ref{1e1} gives something unexpected. The continuous
semigroup that emerges is $ G_2 $, the two-dimensional topological
semigroup described in \eqref{4d1}. However, its representation is not
single-valued:
\[
\begin{gathered}\includegraphics{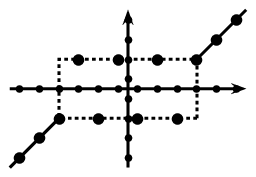}\end{gathered}\quad\mapsto\quad
\begin{gathered}\includegraphics{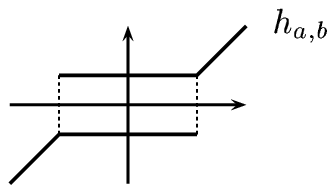}\end{gathered}
\]
Namely, $ h_{a,b} (x) $ for $ x \in (-b,b) $ is $ \pm(a+b) $, that
is, either $ a+b $ or $ -(a+b) $ with probabilities $ 0.5,0.5
$. Such $ h $ is not a function, of course. Rather, it is a
\emph{kernel,} that is, a measurable map from $ \R $ into the space of
probability measures on $ \R $. Composition of kernels is
well-defined, thus, a representation (of a semigroup) by kernels
(rather than functions) is also well-defined.

The stochastic flow in $ G_2 $, resulting from \ref{1e1} via
the scaling limit, is identical to the flow $ ( \xi_{s,t}^{(2)} ) $ of
\ref{sec:4.7}. Its noise is the usual (one-dimensional) white
noise. The representation of $ G_2 $ by kernels turns the abstract
flow into a \emph{stochastic flow of kernels} as defined by Le Jan and
Raimond \cite[Def.~1.1.3]{LJR}. However, a kernel (unlike a function)
introduces an additional level of randomness. When the kernel says
that $ h_{a,b} (x) = \pm (a+b) $, someone has to choose at random one
of the two possibilities. Who makes the decision?

One may treat a point as a macroscopically small collection of many
microscopic atoms, and $ \om \in \Om $ as a macroscopic flow (on the
whole space-time); given $ \om $, atoms are (conditionally)
independent, ``which means that two points\footnote{%
 Or rather, atoms.}
thrown initially at the same place separate'' \cite[p.~4]{LJR}. No
need to deal explicitly with a continuum of independent
choices. ``Turbulent evolutions [are represented] by flows of
probability kernels obtained by dividing infinitely the initial
point'' \cite[p.~4]{LJR}.

Alternatively, one can postulate that if two atoms meet at a
(macroscopic!) point,
they must coalesce. In one-dimensional space (and \emph{sometimes} in
higher dimensions) such a postulate itself prevents a continuum of
independent choices and leads to a flow of maps (the Brownian web is
an example). A countable dense set of atoms makes decisions; others
must obey. A flow of maps is a (degenerate) special case of a flow of
kernels. However, coalescence can produce a flow of maps out of a
non-degenerate flow of kernels, as explained in
\cite[Sect.~2.3]{LJR}.

Conversely, a coalescent flow can produce a non-degenerate flow of
kernels via ``filtering by a sub-noise'' \cite[Sect.~2.3]{LJR}. In the
simplest case (filtering by a trivial sub-noise), we just retain the
one-particle motion of the given coalescent flow, forget the rest of
the flow, and let atoms perform the motion independently.

A large class of flows on $ \R^n $ (and other homogeneous spaces) is
investigated in \cite{LJR}. Some of these flows are shown to be
coalescent and to generate nonclassical noises (neither white nor
black). Flows are homogeneous in space (and isotropic). Thus, we have
a hierarchy of nonclassical models. First, toy models
(recall \ref{1b1}, \ref{1b3}) having a singular time
point. Second, `simple' models (\ref{sec:1.5}, \ref{sec:4.9})
homogeneous in time but having a singular spatial point. Third,
`serious' models (the Brownian web, and Le Jan-Raimond's isotropic
Brownian flows), homogeneous in space and time.

Noises generated by one-dimensional flows (also homogeneous in space
and time) are investigated by Warren and Watanabe \cite{WW}. Spectral
sets of Hausdorff dimension other than $ 0 $ and $ 1/2 $ are found!
Roughly, it answers Question \ref{6c4}; however, these spectral sets
are not perfect --- they have isolated points.

\subsection{Abstract nonsense of Le Jan-Raimond's theory}
\label{sec:8.4}

A new semigroup, introduced recently by Le Jan and Raimond \cite{LJR},
is quite interesting for the theory of stochastic flows and
noises. Its definition involves some technicalities considered here.

A kernel is defined in \cite{LJR} as a measurable mapping from a
compact metric space $ \cM $ to the (also compact) space $ \P(\cM) $ of
all probability measures on $ \cM $. The space $ E $ of all kernels is
equipped with the \sif\ $ \Ec $ generated by evaluations, $ E \ni K
\mapsto K(x) \in \P(\cM) $, at points $ x \in \cM $. Note that every $ \Ec
$\nobreakdash-measurable function uses the values of $ K(x) $ only for
a countable set of points
$ x $, which is scanty, since $ K(x) $ is just measurable (rather than
continuous) in $ x $. Thus, $ (E,\Ec) $ is not a standard Borel
space,\footnote{%
 For a definition, see \cite[Sect.~12.B]{Ke} or \cite[Def.~7.1]{Al}.}
and the composition of kernels is not a measurable operation, which
obscures the technique and makes proofs more difficult (as noted on
page 11 of \cite{LJR}).

Fortunately, the theory can be reformulated equivalently in terms of
Borel operations on standard Borel spaces, as outlined
below. Additional simplification comes from disentangling space and
time (entangled in Theorem 1.1.4 of \cite{LJR}) and explicit use of
the de Finetti theorem.

The hassle about measurability is another manifestation of the
well-known clash between finite-dimensional distributions and
modifications of a random process. Say, for the usual Poisson process
on $ [0,\infty) $, its finite-dimensional distributions do not tell
us whether sample paths are continuous from the left (right), or
not. A process $ X = X(t,\om) $ has a lot of modifications $ Y(t,\om)
$; these satisfy $ \forall t \; \Pr{ \{ \om : X(t,\om)=Y(t,\om) \} } =
1 $, which does not imply $ \Pr{ \{ \om : \forall t \>
X(t,\om)=Y(t,\om) \} } = 1 $. If a process admits continuous sample
paths (like the Brownian motion), the continuous modification is
preferable. If a process is just continuous in probability (like the
Poisson process, but also, say, some stationary Gaussian processes,
unbounded on every interval), we are unable to prefer one modification
to others, in general.

In order to describe the class of all modifications of a random
process, we have two well-known tools: first, a compatible family of
finite-dimensional distributions, and second, a probability measure on
the (non-standard!) Borel space of all (or only measurable; but
definitely, not only continuous) sample paths, whose \sif\ is
generated by evaluations. Assuming the process to be continuous in
probability, we find the first tool much better; joint distributions
depend on points continuously, and everything is standard.

The same for kernels. These may be thought of as sample paths of a
random process whose `time' runs over $ \cM $, and `values' belong to $
\P(\cM) $. However, the process will appear (implicitly) only in Theorem
\ref{8.8}; its finite-dimensional distributions are $ \nu_n
(x_1,\dots,x_n) $ there.

\begin{definition}\label{8.6}
A \emph{multikernel} from a compact metric space $ \cM_1 $ to a compact
metric space $ \cM_2 $ is a sequence $ (P_n)_{n=1}^\infty $ of
continuous maps $ P_n : \cM_1^n \to \P(\cM_2^n) $, compatible in the
sense that\footnote{%
 Here $ \int g \, \D P_n(x_1,\dots,x_n) $ is not an integral in $
 x_1,\dots,x_n $. Rather, $ x_1,\dots,x_n $ are parameters. The
 integral is taken in other variables (say, $ y_1,\dots,y_n $),
 suppressed in the notation and running over $ \cM_2 $.}
\[
\int_{\cM_2^n} g \, \D P_n (x_1,\dots,x_n) = \int_{\cM_2^m} f \, \D P_m
(x_{i_1},\dots,x_{i_m})
\]
for all $ n $ and $ x_1, \dots, x_n \in \cM_1 $, whenever $ i_1, \dots
\, i_m $ are pairwise distinct elements of $ \{ 1,\dots,n \} $, $ f :
\cM_2^m \to \R $ is a continuous function, and $ g : \cM_2^n \to \R $ is
defined by $ g (y_1,\dots,y_n) = f ( y_{i_1},\dots,y_{i_m}) $ for $
y_1,\dots,y_n \in \cM_2 $.
\end{definition}

We do not assume $ i_1 < \dots < i_m $. For example:
\begin{align*}
g(y_1,y_2) = f(y_1) &\imply \int g \, \D P_2 (x_1,x_2) = \int f \, \D P_1
 (x_1) \, ; \\
g(y_1,y_2) = f(y_2) &\imply \int g \, \D P_2 (x_1,x_2) = \int f \, \D P_1
 (x_2) \, ; \\
g(y_1,y_2) = f(y_2,y_1) &\imply \int g \, \D P_2 (x_1,x_2) = \int f \,
 \D P_2 (x_2,x_1) \, .
\end{align*}
Note also that $ x_1, x_2, \dots $ need not be distinct.

\begin{definition}
A multikernel $ (P_n)_{n=1}^\infty $ is \emph{single-valued,} if
\[
\int_{\cM_2^2} g \, \D P_2 (x,x) = \int_{\cM_2} f \, \D P_1 (x) \quad
\text{for all } x \in \cM_1 \, ,
\]
whenever $ g : \cM_2^2 \to \R $ is a continuous function, and $ f :
\cM_2 \to \R $ is defined by $ f(y) = g(y,y) $ for $ y \in \cM_2 $.
\end{definition}

An equivalent definition: $ (P_n)_{n=1}^\infty $ is single-valued, if
\[
\int_{\cM_2^2} \rho \, \D P_2 (x,x) = 0 \quad \text{for all } x \in
\cM_1 \, ,
\]
where $ \rho  : \cM_2^2 \to \R $ is the metric, $ \rho(y_1,y_2) =
\dist(y_1,y_2) $.

Another equivalent definition:
\[
\sup_{\rho(x_1,x_2)\le\eps} \int_{\cM_2^2} \rho \, \D P_2 (x_1,x_2) \to
0 \quad \text{for } \eps \to 0 \, .
\]
(Compare it with continuity in probability.)

My `multikernel' is a time-free counterpart of a `compatible family of
Feller semigroups' of \cite{LJR}. My `single-valued' corresponds to
their (1.7). What could correspond to their `stochastic convolution
semigroup'? It is a single-valued multikernel from $ \cM_1 $ to $
\P(\cM_2) $. Yes, I mean it: maps from $ \cM_1^n $ to $ \P \(
(\P(\cM_2))^n \) $. It may look frightening, but think what happens
if $ \cM_1 $ contains only one point, and $ \cM_2 $ --- only two
points, say, $ 0 $ and $ 1 $. Then a multi\-kernel from $ \cM_1 $ to $
\cM_2 $ is a law of an exchangeable sequence of events. A
single-valued multikernel from $ \cM_1 $ to $ \cM_2 $ would mean that
all events coincide, but we need rather a single-valued multikernel
from $ \cM_1 $ to $ \P(\cM_2) = [0,1] $; nothing but a probability
measure on $ [0,1] $. The De Finetti theorem (see \cite{Al}, for
instance) tells us that every
exchangeable sequence of events arises from a probability measure on $
[0,1] $. Here is a more general result.

\begin{theorem}\label{8.8}
For every multikernel $ (P_n)_{n=1}^\infty $ from $ \cM_1 $ to $ \cM_2
$ there exists a single-valued multikernel $ (\nu_n)_{n=1}^\infty $
from $ \cM_1 $ to $ \P(\cM_2) $ such that
\[
\int_{\cM_2^n} f \, \D P_n (x_1,\dots,x_n) = \int_{(\P(\cM_2))^n} F \,
\D \nu_n (x_1,\dots,x_n)
\]
for all $ n $ and $ x_1,\dots,x_n \in \cM_1 $, whenever $ f : \cM_2^n
\to \R $ is a continuous function, and $ F : (\P(\cM_2))^n \to \R $ is
defined by $ F (\mu_1,\dots,\mu_n) = \int f \,
d(\mu_1\otimes\dots\otimes\mu_n) $ for $ \mu_1,\dots,\mu_n \in
\P(\cM_2) $.
\end{theorem}

\begin{proof}
We choose a discrete probability measure $ \mu_0 $ on $ \cM_1 $ whose
support is the whole $ \cM_1 $. That is, we choose a countable (or
finite) dense set $ A \subset \cM_1 $, and give a positive probability
to each point of $ A $. For every $ n $ we consider the following
measure $ Q_n $ on $ (\cM_1 \times \cM_2)^n $:
\begin{multline*}
\int f_1 \otimes g_1 \otimes \dots \otimes f_n \otimes g_n \, \D Q_n \\
= \int \Big( \int g_1 \otimes \dots \otimes g_n \, \D P_n(x_1,\dots,x_n)
\Big) f_1(x_1) \dots f_n(x_n) \, \D \mu_0 (x_1) \dots \D \mu_0(x_n) \, .
\end{multline*}
In other words, if $ Q_n $ is the distribution of $ (X_1,Y_1; \dots;
X_n,Y_n) $, then $ X_1,\dots,X_n $ are i.i.d.\ distributed $ \mu_0 $
each, and the conditional distribution of $ (Y_1,\dots,Y_n) $ given $
(X_1,\dots,X_n) $ is $ P_n (X_1,\dots,X_n) $. The measure $ Q_n $ is
invariant under the group of $ n! $ permutations of $ n $ pairs, due
to compatibility of the multikernel $ (P_n)_{n=1}^\infty $. For the
same reason, $ Q_n $ is the marginal of $ Q_{n+1} $. Thus, $
(Q_n)_{n=1}^\infty $ is the distribution of an exchangeable infinite
sequence of $ \cM_1 \times \cM_2 $-valued random variables $
(X_n,Y_n) $.

The De Finetti theorem \cite[Th.~3.1 and Prop.~7.4]{Al} states that
the joint distribution of all $
(X_n,Y_n) $ is a mixture of products, in the sense that there exists a
probability measure $ \nu $ on $ \P(\cM_1\times\cM_2) $ such that for
every $ n $, the joint distribution of $ n $ pairs $ (X_1,Y_1), \dots,
(X_n,Y_n) $ is the mixture of products $ Q ^{\otimes n} = Q \otimes
\dots \otimes Q $, where $ Q \in \P(\cM_1\times\cM_2) $ is distributed
$ \nu $. The first marginal of $ Q $ is equal to $ \mu_0 $ (for $ \nu
$-almost every $ Q $), since $ X_n $ are i.i.d.\ ($ \mu_0 $).

Let $ x_1, \dots, x_n \in A $. The event $ X_1=x_1, \dots, X_n=x_n $
is of positive probability. Given the event, the conditional
distribution $ P_n (x_1,\dots,x_n) $ of $ Y_1,\dots,Y_n $ is the
mixture of products $ Q_{x_1} \otimes \dots \otimes Q_{x_n} $, where $
Q_x $ is the conditional measure on $ \cM_2 $, that corresponds to $ Q
$, and $ Q \in \P(\cM_1\times\cM_2) $ is distributed $ \nu $; indeed,
$ \nu $-almost all $ Q $ ascribe the same probability to the event $
X_1=x_1, \dots, X_n=x_n $.

We define $ \nu_n (x_1,\dots,x_n) $ for $ x_1,\dots,x_n \in A $ as the
joint distribution of $ \P(\cM_2) $-valued random variables $ Q_{x_1},
\dots, Q_{x_n} $, where $ Q $ is distributed $ \nu $; then
\begin{multline}\label{54}
\int_{(\P(\cM_2))^n} F \, \D \nu_n (x_1,\dots,x_n) \\
= \int_{\P(\cM_1\times\cM_2)} \bigg( \int_{\cM_2^n} f \, \D (Q_{x_1}
\otimes\dots\otimes Q_{x_n}) \bigg) \, \D \nu(Q) \\
= \int_{\cM_2^n} f \, \D P_n (x_1,\dots,x_n)
\end{multline}
whenever $ f : \cM_2^n \to \R $ is a continuous function, and $ F :
(\P(\cM_2))^n \to \R $ is defined by $ F (\mu_1,\dots,\mu_n) = \int f
\, d(\mu_1\otimes\dots\otimes\mu_n) $ for $ \mu_1,\dots,\mu_n \in
\P(\cM_2) $.

Till now, $ \nu_n (x_1,\dots,x_n) $ is defined for $ x_1,\dots,x_n
\in A $ (rather than $ \cM_1 $). We want to check that $ \int
\ti\rho_2 \, \D \nu_2 (x_1,x_2) \to 0 $ for $ \rho_1 (x_1,x_2) \to 0 $;
here $ \rho_1 $ is a metric on $ \cM_1 $ conforming to its topology,
and $ \ti\rho_2 $ is a metric on $ \P(\cM_2) $ conforming to its weak
topology. Due to compactness of $ \P(\cM_2) $, it is enough to check
that $ \int h^2 \, \D\nu_2(x_1,x_2) \to 0 $ for $ \rho_1 (x_1,x_2) \to
0 $ whenever $ h : \P(\cM_2) \times \P(\cM_2) \to \R $ is of the form
$ h(Q_1,Q_2) = \int f \, \D Q_1 - \int f \, \D Q_2 $ for a continuous
function $ f : \cM_2 \to \R $. Consider $ \ti f : \P(\cM_2) \to \R $,
$ \ti f (Q) = \int f \, \D Q $ for $ Q \in \P(\cM_2) $. We have
\[
\int_{(\P(\cM_2))^2} \ti f \otimes \ti f \, \D\nu_2 (x_1,x_2) =
\int_{\cM_2^2} f \otimes f \, \D P_2 (x_1,x_2) \, ,
\]
which is a special case of \eqref{54}. It may also be written as
\[
\Ex \ti f (Q_{x_1}) \ti f (Q_{x_2}) = \cE{ f(Y_1) f(Y_2) }{ X_1=x_1,
X_2=x_2 } \, ;
\]
here $ Q_{x_1} $ and $ Q_{x_2} $ are treated as random variables on
the probability space $ \( \P(\cM_1\times\cM_2), \nu \) $ (thus, the
two expectations are taken on different probability spaces). The
right-hand side is a continuous function of $ x_1, x_2 $; denote it $
\phi(x_1,x_2) $. We have
\begin{multline*}
\int h^2 \, \D \nu_2 (x_1,x_2) = \Ex \( \ti f (Q_{x_1}) - \ti f
 (Q_{x_2}) \)^2 \\
= \phi(x_1,x_1) - \phi(x_1,x_2) - \phi(x_2,x_1) + \phi(x_2,x_2) \, ,
\end{multline*}
which tends to $ 0 $ for $ \rho_1 (x_1,x_2) \to 0 $. So,
\[
\int_{(\P(\cM_2))^2} \ti \rho_2 \, \D \nu_2(x_1,x_2) \to 0 \quad
\text{for } \rho_1 (x_1,x_2) \to 0 \, .
\]
It follows easily that each $ \nu_n $ is uniformly continuous on $ A^n
$ and, extending it by continuity to $ \cM_1^n $, we get a
single-valued multikernel.
\qqed\end{proof}

Definition \ref{8.6} may be reformulated as follows.

\begin{definition}\label{8.9}
A multikernel from a compact metric space $ \cM_1 $ to a compact
metric space $ \cM_2 $ is a continuous map $ P_\infty : \cM_1^\infty
\to \P (\cM_2^\infty) $, satisfying conditions \textup{(1)} and
\textup{(2)} below. Here $ \cM^\infty = \cM \times \cM \times \dots $
is the product of an infinite sequence of copies of $ \cM $
\textup{(}still a metrizable compact space\textup{).}

\textup{(1)} $ P_\infty $ intertwines the natural actions of the
permutation group of the index set $ \{ 1,2,3,\dots \} $ on $
\cM_1^\infty $ and $ \P (\cM_2^\infty) $ \textup{(}via $ \cM_2^\infty
$\textup{).}

\textup{(2)} For every $ n $, the projection of the measure $ P_\infty
(m) $ to the product $ \cM_1^n $ of the first $ n $ factors depends
only on the first $ n $ coordinates $ m_1, \dots, m_n $ of the point $
(m_1,m_2,\dots) = m \in \cM_1^\infty $.
\end{definition}

Proof of equivalence between definitions \ref{8.6} and \ref{8.9} is
left to the reader.

It is well-known that a continuous map $ \cM_1 \to \P (\cM_2) $ is
basically the same as a linear operator $ C(\cM_2) \to C(\cM_1) $,
positive and preserving the unit. Thus, a multikernel from $ \cM_1 $
to $ \cM_2 $ may be thought of as a positive unit-preserving linear
operator $ C (\cM_2^\infty) \to C (\cM_1^\infty) $ satisfying two
conditions parallel to \ref{8.9}(1,2).

Given three compact metric spaces $ \cM_1, \cM_2, \cM_3 $, a
multikernel from $ \cM_1 $ to $ \cM_2 $ and a multikernel from $ \cM_2
$ to $ \cM_3 $, we may define their composition, a multikernel from $
\cM_1 $ to $ \cM_3 $. In terms of operators it is just the product of
two operators, $ C (\cM_3^\infty) \to C (\cM_2^\infty) \to C
(\cM_1^\infty) $.

\enlargethispage*{10pt}

The set of all multikernels from $ \cM_1 $ to $ \cM_2 $, treated as
operators $ C (\cM_2^\infty) \to C (\cM_1^\infty) $, is a closed (and
bounded, but not compact) subset of the operator space equipped with
the strong operator topology. Thus, the set of multikernels becomes a
Polish space (that is, a topological space underlying a complete
separable metric space).

Composition of multikernels, $ C (\cM_3^\infty) \to C (\cM_2^\infty)
\to C (\cM_1^\infty) $, is a (jointly) continuous operation. (Indeed,
the product of operators is continuous in the strong operator
topology, as far as all operators are of norm $ \le 1 $.)

So, multikernels from $ \cM $ to $ \cM $ are a Polish semigroup (that
is, a topological semigroup whose topological space is Polish).

\printindex

\bigskip
\filbreak
{
\small
\begin{sc}
\parindent=0pt\baselineskip=12pt
\parbox{2.3in}{
Boris Tsirelson\\
School of Mathematics\\
Tel Aviv University\\
Tel Aviv 69978, Israel
\smallskip
\emailwww{tsirel@tau.ac.il}
{www.tau.ac.il/\textasciitilde tsirel/}
}
\end{sc}
}
\filbreak

\end{document}